\documentclass[a4paper,12pt]{amsbook}
\usepackage{amsmath,amssymb,amscd,amsfonts,amsmath,mathrsfs}

\DeclareFontFamily{U}{mathabx}{\hyphenchar\font45}
\DeclareFontShape{U}{mathabx}{m}{n}{
      <5> <6> <7> <8> <9> <10> gen * matha
      <10.95> <12> <14.4> <17.28> <20.74> <24.88> matha10}{}
\DeclareSymbolFont{mathabx}{U}{mathabx}{m}{n}
\DeclareMathSymbol{\partialslash}  {0}{mathabx}{"43}

\usepackage[mathcal]{euscript}

\usepackage[all]{xy}

\newcommand{\C}{{\mathbf C}}            
\newcommand{\CP}{\mathbf{CP}^1}         
\newcommand{\CPt}{\widehat{\mathbf{CP}}^1}    
\newcommand{\Om}{{\mathbf \Omega}}           
\newcommand{\R}{{\mathbf R}}            
\newcommand{\N}{{\mathbf N}}            
\newcommand{\Z}{{\mathbf Z}}            
\renewcommand{\H}{{\mathbb H}}        
\newcommand{\Cs}{{\mathcal H}}           
\newcommand{\E}{{\mathcal E}}           
\newcommand{\F}{{\mathcal F}}           
\newcommand{\EC}{{\mathcal C}}          
\newcommand{\Nahm}{\mathsf{N}}           
\newcommand{\Sky}{{\mathcal R}}         
\newcommand{\dbar}{\bar{\partial}}             
\renewcommand{\d}{\mbox{d}}               
\newcommand{\D}{\mathrm{D}}             
\newcommand{\Dou}{\mathcal{D}}          
\newcommand{\germ}{(\mathrm{germ}\ \mathcal{D})}          
\newcommand{\germz}{(\mathrm{germ}\ \mathcal{D}_{z_{0}})} 
\newcommand{\Dir}{\partialslash  }      
\renewcommand{\S}{\mathrm{E}}             
\newcommand{\Sg}{(\mathrm{germ}\ \mathrm{E})}             
\newcommand{\Id}{\mbox{Id}}             
\newcommand{\ra}{\to}           
\newcommand{\xra}{\xrightarrow}
\newcommand{\lra}{\longrightarrow}              
\newcommand{\absl}{\left\arrowvert}     
\newcommand{\absr}{\right\arrowvert}    
\newcommand{\Dm }{\mbox{Dom}_{\mbox{max}}}          

\renewcommand{\setminus}{\smallsetminus}
\newcommand{\inc }{\hookrightarrow}   
\newcommand{\contr }{\llcorner }
\renewcommand{\j}{\jmath}  

\newtheorem{prop}{Proposition}[chapter]
\newtheorem{hyp}[prop]{Hypothesis}
\newtheorem{rk}[prop]{Remark}
\newtheorem{clm}[prop]{Claim}
\newtheorem{lem}[prop]{Lemma}
\newtheorem{defn}[prop]{Definition}
\newtheorem{cor}[prop]{Corollary}
\newtheorem{thm}[prop]{Theorem}

\makeatletter

\@addtoreset{equation}{chapter}
\makeatother



\DeclareMathSizes{10}{10}{6}{5}
\DeclareMathSizes{11}{11}{7}{5}
\DeclareMathSizes{12}{12}{7}{5}

\title{
\begin{tabular}{p{-2cm} c} 
& \\
& \\
& \hspace{-3.7cm} {\Huge Nahm transform for integrable connections} \\
& \\
& \hspace{-3.7cm} {\Huge on the Riemann sphere}
\end{tabular}
}

\author{
\begin{tabular}{p{-2cm} c}
& \\
& \\
& \\
& \hspace{-1cm} \bf{\Huge Szil\'ard Szab\'o}\\
& \\
& \hspace{-1cm} szilard.szabo@m4x.org
\end{tabular}
}

\begin{document}
\maketitle
\tableofcontents

\chapter{Introduction}

\section{Historical context and abstract of the thesis}
Nahm transform is a non-linear analog for instantons of the usual Fourier transform on functions. 
It has been extensively studied starting from the beginning of the 1980's, inspired by the seminal
work of M. F. Atiyah, V. Drinfeld, N. J. Hitchin and Yu. I. Manin on a correspondence 
(the \emph{ADHM-transform}) between finite-energy solutions of the \emph{Yang-Mills equations} 
and some algebraic data (see \cite{adhm}, \cite{Don-Kronh}). 
The Yang-Mills equations are the anti-self-duality equations for a unitary connection on a 
Hermitian vector bundle defined over $\R^{4}$; their finite-energy solutions are called \emph{instantons}. 

Since then, it turned out that the general picture concerning this correspondence is as follows: let $X$ be
any manifold obtained as a quotient of $\R^{4}$ by a closed additive subgroup $\Lambda$. The solutions of the 
Yang-Mills equations invariant by $\Lambda$ (that are clearly not of finite energy in the case $\Lambda \neq \{ 0 \}$) can
be identified in an obvious manner to solutions of a system of differential equations on $X$, called the 
\emph{reduction} of the Yang-Mills equations. On the other hand, denoting by $(\R^{4})^{\ast }$ the
dual of the vector space $\R^{4}$, $\Lambda$ determines a closed additive subgroup $\Lambda^{\ast}$ called the 
\emph{dual subgroup} by saying that an element $\xi \in (\R^{4})^{\ast }$ is in $\Lambda^{\ast}$ if and only if $\xi (\lambda ) \in \Z$ 
for all $\lambda \in \Lambda$. Hence, we can form the \emph{dual manifold} $X^{\ast } = (\R^{4})^{\ast }/\Lambda^{\ast }$ of $X$, that also
admits a reduction of the Yang-Mills equations. Nahm transform is then a procedure that maps 
solutions of the reduced equations on $X$ to solutions of the reduced equations on $X^{\ast}$ bijectively up to
overall gauge transformations on both sides. One remarks that there is a canonical isomorphism between 
$((\R^{4})^{\ast })^{\ast }$ and $\R^{4}$, as well as between $(\Lambda^{\ast })^{\ast }$ and $\Lambda$. Therefore, if we start from 
a solution of the reduced equations on $X$ and iterate Nahm transform twice, we again get a solution of the 
reduced equations on $X$. One important property analogous to usual Fourier transform is that in some cases 
the solution we get this way is, up to a coordinate change $x \mapsto -x$, known to be the solution we started
with; that is, Nahm transform is (up to a sign) \emph{involutive}. Moreover, in some cases one knows that 
the moduli spaces of solutions of the reduced equations modulo gauge transformations on $X$ and on $X^{\ast}$ 
are smooth hyper-K\"ahler manifolds with respect to the metric induced by $L^{2}$-norm and the 
complex structures induced by $\R^{4}$; Nahm transform is then a hyper-K\"ahler isometry between these moduli
spaces. This is to be compared with Parseval's theorem which states that usual Fourier transform defines an
isometry between $L^{2}$-spaces of functions. 

Putting $\Lambda = \{ 0 \}$, one gets $X=\R^{4}$ and $\Lambda^{\ast}=\R^{4}$, so $X^{\ast}=\{ 0 \}$. In this case, 
Nahm transform reduces to the ADHM-transform. 
There are several other examples of Nahm transform in the literature for different subgroups of $\R^{4}$; 
for a nice exposition of these, see the survey paper \cite{jarsur} of M. Jardim. In this work, we are
concerned with the case $\Lambda=\R^{2}$. In this case, the base manifold is $X=\R^{2}$, and its dual $X^{\ast}$ is 
another copy of the real plane that we shall denote by $\hat{\R}^{2}$. These are non-compact manifolds,
with compactifications the Riemann spheres $\CP$ and $\CPt$ respectively. The reduction of the original 
(Yang-Mills) equations can be viewed in two different ways depending on the complex structure that we choose: 
they are the equations defining an integrable connection with harmonic metric, or equivalently, those defining
a Higgs bundle with Hermitian-Einstein metric. Now, it turns out that there are no smooth 
solutions on the Riemann sphere of either one of these equations except for the trivial ones (c.f. \cite{Hit}). 
However, there are solutions having prescribed singularities in some points, and the solutions of one equation
are still in correspondence with those of the other: this is proved by O. Biquard and Ph. Boalch in
\cite{Biq-Boa}. We establish, under some hypotheses on the singularity behavior, Nahm transform for singular 
integrable connections (or equivalently, singular Higgs bundles) on the Riemann sphere. 
Note that Nahm transform for singular objects
have already been studied by O. Biquard and M. Jardim in \cite{Biq-Jar} and by S. A. Cherkis and A. Kapustin in
\cite{Cher-Kap}. On the other hand, using different techniques, B. Malgrange has defined in \cite{Mal} a so-called 
Fourier-Laplace transform for integrable connections with singularities on the Riemann sphere behaving in the 
same manner on the level of singularity data as the one we define here. It is therefore very natural to
believe that these two transforms actually agree. One difference between these works is, however, the
transformation of a parabolic structure and an adapted harmonic metric at the singularities in our case; 
for details, see Section \ref{sectharmmetr}. 

The construction follows the main ideas of other Nahm transforms found in literature. Namely, 
in Section \ref{sectFr} we define positive and negative spinor bundles $S^{\pm}$ over $\CP$, as well as a Dirac operator 
$$
    \Dir: S^{+} \otimes E \lra S^{-} \otimes E .
$$ 
We then let $\xi \in \hat{\C} \setminus \hat{P}$ be a parameter, where $\hat{P}$ is the singular locus of the
transformed objects, and for all $\xi$ twist the operator $\Dir$ by some flat connection to 
obtain a family of operators $\Dir_{\xi}$. In Section \ref{sectpfFr} we prove that the kernel of these twisted 
operators vanish and that the cokernels form a finite-dimensional space. Furthermore, this dimension is
independent of $\xi$; we then define the transformed vector bundle $\hat{E}$ on $\hat{\C}$ as the vector bundle 
with fiber over $\xi$ given by $coKer(\Dir_{\xi})$. In Section \ref{l2coh} we carry out an analog of $L^{2}$-Hodge
theory of a compact K\"ahlerian manifold in this case; namely we establish an isomorphism between this cokernel 
and the first $L^{2}$-cohomology of an elliptic complex, as well as harmonic $1$-forms with respect to the 
Laplacian of the Dirac operator. We then go on to define the transformed  flat bundle 
and the transformed Hermitian metric in Section \ref{flattr}, and we extend the flat bundle over
the singularities -- so defining the transformed meromorphic integrable connection -- in Section \ref{flatext}. 
The transformed metric is then shown to be Hermitian-Einstein in Section \ref{sectharm}. 
Next, in Section \ref{ident} we  give a completely explicit description of  the fibers of the transformed
bundle, first in terms of hypercohomology of a sheaf map, then in terms of the corresponding spectral set. 
Then come the constructions of the extensions of the
transformed Higgs bundle to the singular points (Section \ref{ext}). This allows us to 
obtain the singularity data of the transformed Higgs bundle in Sections \ref{sing} and \ref{parweightsect}, 
and we complete the transform by computing the topology of the transformed Higgs bundle in Section
\ref{secttop}. Finally, Chapter \ref{chinv} deals with the involutivity property of the transform.

\section{Integrable connection point of view} \label{flatintro}
Let $\C$ be the complex line, with its natural holomorphic coordinate $z=x+iy$ and 
Euclidean metric $\vert\mbox{d}z \vert^{2}$; and let $\CP$ be the complex projective line. 
Let $E \to \CP$ be a rank $r$ holomorphic vector bundle on the Riemann sphere,
and $D$ be a meromorphic integrable connection on it, with \emph{first order} or \emph{logarithmic singularities}
at the points of a finite set $\{ p_1,\ldots,p_{n} \} =P \subset \C$ and a second
order singularity at infinity. 
In other words,  on a small disk $\Delta(p_{j},\varepsilon)$ centered at  
$p_j \in  P$ in a holomorphic basis $\{ \tau^j_k \}_{k=1,...,r}$ of $E$, $D$ is
of the form  $D^{j}+b^{j}$ where $b^{j}$ is a holomorphic $1$-form on the disk and 
\begin{eqnarray} \label{dpj}
      D^{j}=\mbox{d}+ \frac{A^{j}}{z-p_{j}} \d z \land .
\end{eqnarray}
We suppose furthermore that $A^{j}$ is diagonal:
$$
    A^{j}= \begin{pmatrix}
                  0 & & & & & \\
                  & \ddots & & & & \\
                  & & 0 & & & \\
                  & & & \mu^{j}_{r_{j}+1} & & \\
                  & & & & \ddots & \\ 
                  & & & & & \mu^{j}_{r}
           \end{pmatrix};
$$
it is called the \emph{residue} of $D$ at $p_{j}$, and $1 \leq 
r-r_{j}\leq r$ is the rank of  $A^{j}$. For convenience, we put
$\mu^{j}_{1}=\ldots=\mu^{j}_{r_{j}}=0$, so that
$A^{j}=diag(\mu^{j}_{k})_{k=1,\ldots r}$. We will often make use of the holomorphic local
decomposition 
\begin{eqnarray} \label{regsing}
    E^{j} = E^{j}_{reg} \oplus E^{j}_{sing}, 
\end{eqnarray}
into the \emph{regular} and \emph{singular components} of $E$ near $p_{j}$; here by definition $E^{j}_{reg}$ is the
holomorphic subbundle of $E^{j}= E|_{\Delta(p_{j},\varepsilon)}$  spanned by $\{ \tau^j_k \}_{k=1,...,r_{j}}$, 
and $E^{j}_{sing}$ is the one spanned by $\{ \tau^j_k \}_{k=r_{j}+1,...,r}$. Intrinsically, 
$E^{j}_{sing}$ is the sum of the generalized eigenspaces corresponding to all eigenvalues converging to
infinity of the integrable connection, whereas $E^{j}_{reg}$ is the sum of the generalized eigenspaces 
corresponding to the eigenvalues that remain bounded. 

In a similar manner, at infinity $D$ is supposed to be equal (up to a holomorphic term) 
to a meromorphic local model having a \emph{second order pole}, 
so that in a holomorphic basis $\{ \tau_{k}^{\infty} \}_{k=1,\ldots,r}$ on a disk $\C \setminus \Delta(0,R)$ 
corresponding to a standard neighborhood of infinity in $\CP$, 
it is of the form $D=D^{\infty } + b^{\infty }$ where $b^{\infty }$ is now a holomorphic $1$-form in the given neighborhood
of infinity, and 
\begin{eqnarray} \label{dinf}
    D^{\infty }=  \mbox{d} + \left( A + \frac{C}{z} \right) \d z \land
\end{eqnarray}
is the second order model with diagonal leading term 
$$
     A =  \begin{pmatrix}
                \xi_{1} & & & & & & & \\
                & \ddots & & & & & & \\
                & & \xi_{1}  & & & & & \\
                & & & \ddots & & & & \\
                & & & & \ddots & & & \\
                & & & & & \xi_{n'} & & \\
                & & & & & & \ddots & \\
                & & & & & & & \xi_{n'}
                        \end{pmatrix}
$$ 
and residue 
$$
    C = \begin{pmatrix}
                \mu_1^{\infty} & & \\
                & \ddots & \\
                & & \mu_r^{\infty}
                        \end{pmatrix}.
$$
Here $\{ \xi_{l} \}_{l=1}^{n'}$ are the distinct eigenvalues of $A$. Each $\xi_{l }$ appears in neighboring positions
$k=1+a_{l},\ldots,a_{l+1}$, in particular its multiplicity is $m_{l}=a_{l+1}-a_{l}$. Of course, we must then have 
$a_{1}=0$ and $a_{n'+1}=r$. In line with the above notation, we set $r_{\infty }=0$ and 
$C=diag( \mu_{k}^{\infty})_{k=1,\ldots,r}$. Furthermore, we will write 
$$
      A=diag( \{ \xi_{l}, m_{l} \} )_{l=1,\ldots,n'} 
$$
for the diagonal matrix $A$ as given above, meaning that $A$ is diagonal with $m_{l}$ neighboring eigenvalues
equal to $\xi_{l}$.   
\begin{defn} 
The integrable connections having singularities near the points of $P \cup \{ \infty \}$ as described above will
be called \emph{meromorphic integrable connections with logarithmic singularities in $P$ and a second-order
  singularity at infinity}, or for simplicity \emph{meromorphic integrable connections} although they are by far not
all the meromorphic integrable connections. 
\end{defn}

\section[The transform of the integrable connection]{The transform of the meromorphic integrable connection}\label{transmer}
Let $(E,D)$ be a stable vector bundle with a meromorphic integrable connection on the sphere. 
Our aim in this paper is to define another complex bundle $\hat{E}$ with a meromorphic connection $\hat{D}$ on
the sphere out of $(E,D)$, which we call the \emph{transformed meromorphic integrable connection}. 
Just as the initial connection, the transformed one will also admit a finite number of simple poles
in points of the line and a second-order pole at infinity. 

In order to define the transformed vector bundle $\hat{E}$, first we need to set some notation. 
Let $\hat{\C}$ be another copy of $\C$. (The importance of distinguishing the two copies of $\C$ is to help us
avoid confusions.) For a parameter $\xi \in \hat{\C}$, consider the following deformation of $D$:
\begin{eqnarray} \label{intdef}
     D^{int}_{\xi } = D - \xi \d z \land ,
\end{eqnarray} 
where $\xi :E \ra E$ stands for multiplication by $\xi$.
Since we only change the $(1,0)$-part of $D$, and by an endomorphism that is independent of $z$,  this is then
another meromorphic integrable connection, with the same underlying
holomorphic bundle as for $D$. Furthermore, its unitary and self-adjoint parts are given by 
\begin{align} 
     D^{+}_{\xi } & = D^{+} - \frac{\xi}{2} \d z + \frac{\bar{\xi}}{2} \d \bar{z} \label{intdefunit}\\
     \Phi^{int}_{\xi } & = \Phi - \frac{\xi}{2} \d z - \frac{\bar{\xi}}{2} \d \bar{z} \label{intdefsa}.  
\end{align} 
Consider the following family in $\xi$ of elliptic complexes $\EC^{int}_{\xi }$ over $\C \setminus P$: 
\begin{equation}\label{ellcompl}
        \Omega^{0}\otimes E \xrightarrow{\text{$D^{int}_{\xi}$}}  \Omega^1 \otimes E \xrightarrow{\text{$D^{int}_{\xi}$}} \Omega^2 \otimes E.
\end{equation}
Fix a Hermitian metric $h$ on $E$ for which the holomorphic sections of the extension at the singularities are
bounded (above and below) by a positive constant, and denote by $\hat{E}^{int}_{\xi }$ the first
$L^{2}$-cohomology of the complex (\ref{ellcompl}) for this metric. In Theorems \ref{Fredholm}
and \ref{laplker} we show that there exists a finite set $\hat{P} \subset \hat{\C}$ such that for 
$\xi \in \hat{\C}\setminus \hat{P}$ the first $L^{2}$-cohomologies of this complex are finite-dimensional of the 
same dimension for all $\xi$. 
\begin{defn} The \emph{transformed vector bundle} $\hat{E}$ is then the vector bundle over 
$\hat{\C} \setminus \hat{P}$ whose fiber over $\xi \in \hat{\C} \setminus \hat{P}$ is the first $L^{2}$-cohomology
$L^{2}H^{1}(D^{int}_{\xi})$ of $\EC^{int}_{\xi }$. 
\end{defn}
Let $\xi_{0} \in \hat{\C} \setminus \hat{P}$, and let $f(z) \in \hat{E}_{\xi_{0}}$ be a class in the first cohomology 
of $\EC^{int}_{\xi_{0}}$. 
\begin{defn}
The \emph{transformed flat connection}  $\hat{D}$ is by definition the flat connection whose parallel section 
$f(\xi; z)$ extending $f$ in some neighborhood of $\xi_{0}$ is given by the first $L^{2}$-cohomology classes in 
$\EC^{int}_{\xi }$ of 
$$
   e^{(\xi - \xi_{0})z }f(z). 
$$
\end{defn}
Finally, $h$ induces a natural Hermitian metric $\hat{h}$ on $\hat{E}$ as follows: in Theorem \ref{laplker} 
we show that any class in $L^{2}H^{1}(D_{\xi})$ can be represented by a unique harmonic $1$-form with respect to
the Laplacian of the Dirac operator. 
\begin{defn}
The \emph{transformed Hermitian metric} $\hat{h}$ on $\hat{E}$ is defined by the $L^{2}$-norm of harmonic 
representatives. 
\end{defn}
All this will be explained in more detail in Section \ref{flattr} and in Definition \ref{deftrhm}. 

When one considers an integrable connection, there exists sometimes a privileged fiber metric on the bundle, 
namely a harmonic one. In order to be able to define harmonicity, decompose as usual $D$ into its
unitary and self-adjoint part
\begin{equation}\label{decompsaunit}
           D = D^{+} + \Phi, 
\end{equation}
put $\nabla_{D^{+}}$ or simply $\nabla^{+}$ for the covariant derivative associated
to the connection $D^{+}$ (so that $\nabla^{+}t$ makes sense for a tensor $t$ of
arbitrary type $(T\CP)^{p} \otimes (T^{\ast }\CP)^{q} \otimes E^{r}\otimes (E^{\ast })^{s}$)  
and denote by $(\nabla^{+})_{h}^{\ast}$ the adjoint operator of $\nabla^{+}$ with
respect to $h$. 
\begin{defn}
The Hermitian metric $h$ is called \emph{harmonic}, if it satisfies the equation 
\begin{eqnarray} \label{harm}
     (\nabla^{+})_{h}^{\ast}\Phi =0.
\end{eqnarray}
\end{defn}
This is a second-order non-linear partial differential equation in $h$. 

Here is the main result of this thesis in a special case (the one without parabolic structures, see Definition
\ref{parstr}). 
\begin{thm}\label{naivethm}
Let $(E,D,h)$ be any meromorphic integrable connection with logarithmic singularities in $P$ as in (\ref{dpj}), 
and a double pole (\ref{dinf}) at infinity, endowed with a harmonic metric $h$. 
Suppose that the eigenvalues of the polar part of $D$ in the punctures satisfy the following assumptions: 
\begin{enumerate}
\item{for fixed $j \in \{ 1,\ldots,n\}$, the complex numbers $\mu^{j}_{k}$ for $k=r_{j}+1,\ldots,r$ are all different, 
      and $\Re \mu^{j}_{k} \notin \Z$}
\item{for fixed $l \in \{1,\ldots,n'\}$, the complex numbers  $\mu^{\infty}_{k}$ for $k=1+a_{l},\ldots,a_{l+1}$ are all different, 
      and $\Re \mu^{\infty}_{k} \notin \Z$}
\end{enumerate}
Then the set of punctures $\hat{P}\in \hat{\C}$ of the transformed bundle is the set $\{ \xi_{1},\ldots,\xi_{n'} \}$ 
of distinct eigenvalues of the leading order term $A$ of $D$ at infinity. 
For $\xi \notin \hat{P}$, the first $L^{2}$-cohomologies of (\ref{ellcompl}) are finite dimensional vector
spaces of the same dimension. They match up to define a smooth vector bundle $\hat{E}$ of rank 
\begin{equation} \label{trrank}
      \hat{r}=\sum_{j=1}^{n} (r-r_{j}) 
\end{equation}
over $\hat{\C} \setminus \hat{P}$. $\hat{D}$ is a flat connection on $\hat{E}$. 
It underlies a meromorphic integrable connection (that we continue to denote $(\hat{E},\hat{D})$) of degree 
$deg(\hat{E})=deg(E)$, called the \emph{transformed meromorphic connection}. It has logarithmic singularities 
in $\hat{P}$ and a double pole at infinity. The non-vanishing eigenvalues of the residue in $\xi_{l} \in \hat{P}$ are 
$\{ -\mu_{1+a_{l}}^{\infty},\ldots,-\mu_{a_{l+1}}^{\infty} \}$. The eigenvalues of the second-order term of the transformed 
meromorphic connection are $\{-p_{1},\ldots,-p_{n}\}$, the multiplicity of $-p_{j}$ being $(r-r_{j})$; 
the eigenvalues of its residue at infinity on the eigenspace of the second-order term corresponding
to $-p_{j}$ are $\{-\mu_{r_{j}+1}^{j},\ldots,-\mu_{r}^{j}\}$. Finally, $\hat{h}$ is harmonic for $\hat{D}$. 
\end{thm}

\begin{rk}
The assumptions (1) and (2) of the theorem are clearly generic in the parameter space of all 
possible eigenvalues.  
\end{rk}

This theorem actually follows from the more general statement \ref{mainthm}. In order to understand the 
more general setup, one needs to consider meromorphic connections endowed with a parabolic structure. 

\section{Parabolic structure and adapted harmonic metric}\label{sectharmmetr}

Actually, we can suppose more structure on the integrable connection: namely, that it comes
with a parabolic structure on $P$ and at infinity. 
\begin{defn} \label{parstr}
A \emph{parabolic structure} on $(E,D)$ is the data of a strictly decreasing filtration by vector subspaces
$$
     E_{p} =  F_{0}E_{p} \supset F_{1}E_{p} \supset \ldots \supset F_{b_{p}-1} E_{p}  \supset F_{b_{p}} E_{p} = \{ 0 \}
$$
(where $1\leq b_{p}\leq r$) of the fiber $E_{p}$ of $E$ in each singular point $p \in P \cup \{ \infty \}$, 
called the \emph{parabolic flag}, such that each $F_{m}$ is spanned by some of the restrictions 
$\{ \tau^{j}_{k}(p) \}_{k=1}^{r}$ of the holomorphic basis to the singularity $p=p_{j}$ or $\infty$, 
together with a sequence of corresponding real numbers  
$$
    0\leq \tilde{\beta}^{j}_{1}< \ldots < \tilde{\beta}^{p}_{b_{p}} < 1
$$
called the \emph{parabolic weights}.
\end{defn}
\begin{rk}
All parabolic weights can be assigned a natural multiplicity, namely the dimension of the corresponding 
graded of the filtration: more precisely, the multiplicity of $\tilde{\beta}^{p}_{k}$ for any $p \in P \cup \{ \infty \}$ 
and any $k\in \{1,\ldots,b_{p} \}$ is by definition 
$$
     \dim (F_{k-1}E_{p}/F_{k}E_{p}).
$$
We will write 
$$
    0\leq {\beta}^{p}_{1}\leq \ldots \leq {\beta}^{p}_{r} < 1 
$$
for the parabolic weights repeated according to their multiplicities, and use this numbering of the weights 
throughout the whole paper instead of the one in their definition. Moreover, we write $\beta^{j}_{k}$ instead of 
$\beta^{p_{j}}_{k}$. 
\end{rk}
\begin{rk} The order of the $\tau^{\infty}_{k}$ spanning $F_{m}E_{\infty}$ in the above definition is not necessarily the same 
as the one in which the eigenvalues of the second-order term A at infinity appear in one
group, as supposed in (\ref{dinf}). However, this will not cause any confusion in the sequel, because 
the basis vectors at infinity in this latter order still have well-defined parabolic weights 
(which are then not necessarily increasing). 
\end{rk}
\begin{defn} \label{parintcon}
A meromorphic integrable connection $(E,D)$ with described local models and parabolic 
structures at the punctures will be called \emph{parabolic integrable connection}. The \emph{parabolic degree}
of $E$ with respect to the given parabolic structure is the real number 
\begin{equation} \label{pardeg}
    deg_{par} (E) = deg(E) + \sum_{j \in \{1,\ldots,n,\infty \}} \sum_{k=1}^{r} \beta^{j}_{k}, 
\end{equation}
where $deg(E)$ is the standard (algebraic geometric) degree of $E$, and the sum is taken over all 
parabolic weights for all punctures $p$. The \emph{slope} of the parabolic integrable connection 
is the real number 
\begin{equation} \label{slope}
     \mu_{par}(E) = \frac{deg_{par} (E)}{rk(E)}, 
\end{equation}
and $(E,D)$ is said to be \emph{parabolically stable} (resp. \emph{semi-stable}) if for any subbundle $F$
invariant with respect to $D$ and endowed with the induced parabolic structure over the singularities, 
the inequality 
\begin{equation} \label{stable}
     \mu_{par}(F) < \mu_{par}(E) 
\end{equation}
(respectively $\mu_{par}(F) \leq \mu_{par}(E)$) holds. Finally, $(E,D)$ is said to be \emph{parabolically polystable} 
if it is a direct sum of parabolically stable bundles that are all invariant by $D$ and of the same slope as
$E$. 
\end{defn}
\begin{rk}\label{pardegvanish}
The notions of stability, semi-stability and polystability make sense for meromorphic integrable 
connections without a parabolic structure as well: in the corresponding definitions, 
one only needs to set all parabloic weights equal to $0$. Notice however that by the residue theorem 
we have 
\begin{align*}
     deg(E) & =  - \Re tr(Res(D,\infty)) -\sum_{j \in \{1,\ldots,n \}} \Re tr(Res(D,p_{j})) \\
            & =  \sum_{k=1}^{r} \Re \mu^{\infty}_{k}-\sum_{j \in \{1,\ldots,n\}} \sum_{k=1}^{r} \Re \mu^{j}_{k}, 
\end{align*}
(the change of sign coming from the fact that the eigenvalues of the residue at infinity are $-\mu^{\infty}_{k}$
because in the local coordinate $w=1/z$ we have $\d z/z=-\d w/w$.)
Therefore (\ref{pardeg}) is in fact equal to 
\begin{align*}
     \sum_{k=1}^{r} (\beta^{\infty}_{k} + \Re \mu^{\infty}_{k}) + \sum_{j \in \{1,\ldots,n\}} \sum_{k=1}^{r} (\beta^{j}_{k} - \Re \mu^{j}_{k}) = 
     \sum_{j \in \{1,\ldots,n, \infty\}} \sum_{k=1}^{r} \gamma^{j}_{k}, 
\end{align*}
where $\gamma^{j}_{k}$ are the parabolic weights of the local system at $p_{j}$ (Proposition 11.1, \cite{Biq97}). 
On the other hand, the parabolic degree of an integrable connection is always equal to $0$: this 
follows from the Gauss-Chern formula 2.9 of \cite{Biq91}. 
Therefore, the case where the parabolic weights $\beta^{j}_{k}$ of the integrable connection vanish is not the one 
where the parabolic weights $\gamma^{j}_{k}$ of the representation of the fundamental group vanish, 
and where by Remark 8.2 of \cite{Biq-Boa} stability reduces to irreducibility of the corresponding representation. 
\end{rk}
\begin{defn}\label{adaptdef}
A Hermitian fiber metric $h$ on $E$ is said to be \emph{adapted} to the parabolic structure if 
near the logarithmic punctures in the holomorphic bases $\tau^{j}_{k}$ it is mutually bounded with 
the diagonal model 
\begin{equation} \label{parmatr}
     diag(|z-p_{j}|^{2 \beta^{j}_{k}})_{k=1}^{r}, 
\end{equation}
and at infinity in the holomorphic basis $\tau^{\infty}_{k}$ it is mutually bounded with  
\begin{equation} \label{parmatrinf}
     diag(|z|^{-2 \beta^{\infty}_{k}})_{k=1}^{r}. 
\end{equation}
\end{defn}
\begin{rk}
In general, without the hypothesis of semisimplicity of the residue in the puctures made in Section
\ref{flatintro}, the local models of the metric near the punctures are more complicated than in the 
above definition: e.g. for the regular singularities one has to take into account an extra filtration 
induced by the nilpotent part of the residue, and add a factor $|\ln(r)|^{k}$ on the 
corresponding $k$-th graded, see the Synopsis of \cite{Sim}. 
\end{rk}
Here is the important existence result of the theory: 
 \begin{thm}[Sabbah C. \cite{Sab-harm}; Biquard O., Boalch Ph. \cite{Biq-Boa}] \label{harmexistthm}
Let $(E,D)$ be a parabolically stable parabolic integrable connection. Then there exists a unique
harmonic Hermitian metric $h$ adapted to the parabolic structure. 
\end{thm}
\begin{rk}
Actually, in the above articles this theorem is proved to hold for parabolic integrable connections having 
poles of arbitrary order in the punctures. On the other hand, for integrable connections with only regular 
singularities, it had already been shown by C. Simpson, see \cite{Sim}. 
\end{rk}
We are now ready to describe the more general version of Nahm transform: that for parabolic integrable connections. 
\begin{thm} \label{mainthm} 
Let $(E,D)$ be any parabolic integrable connection on the sphere with logarithmic singularities in $P$ as in 
(\ref{dpj}), and a double pole (\ref{dinf}) at infinity. Suppose that the eigenvalues of its polar parts $\mu$ and 
the parabolic weights $\beta$ in the punctures satisfy the following assumptions: 
\begin{enumerate}
\item{for fixed $j \in \{ 1,\ldots,n\}$, the complex numbers $\mu^{j}_{k}-\beta^{j}_{k}$ for $k=r_{j}+1,\ldots,r$ are distinct
    and different from $0$, the parabolic weights $\beta^{j}_{k}$ for $k=1,\ldots,r_{j}$ are $0$ and finally 
    $\Re \mu^{j}_{k} \notin \Z$ for $k=r_{j}+1,\ldots,r$}
\item{for fixed $l \in \{1,\ldots,n'\}$, the complex numbers  $\mu^{\infty}_{k}-\beta^{\infty}_{k}$ for $k=1+a_{l},\ldots,a_{l+1}$ are
    distinct and different from $0$, and $\Re \mu^{\infty}_{k} \notin \Z$}
\end{enumerate}
Then, in addition to the conclusions of Theorem \ref{naivethm}, the transformed bundle $(\hat{E},\hat{D})$ 
carries a natural parabolic structure in the punctures (that we will call 
\emph{transformed parabolic structure}), such that the transformed metric of the harmonic metric is adapted to
it. Moreover, the set of its non-vanishing parabolic weights in $\xi_{l} \in \hat{P}$ 
is equal to the set of parabolic weights $\{ \beta^{\infty }_{1+a_{l}},\ldots,\beta^{\infty }_{a_{l+1}} \}$ 
of $E $ at infinity, restricted to the eigenspace of $A$ corresponding to the eigenvalue $\xi_{l}$; 
whereas the parabolic weights of $\hat{E}$ at infinity restricted to the eigenspace of the second-order term of $\hat{D}$ 
corresponding to the eigenvalue  $-p_{j}$ are equal to the parabolic weights 
$\{ \beta^{j}_{r_{j}+1},\ldots, \beta^{j}_{r}\}$ of $E $ at $p_{j}$. All these statements are to be understood 
with multiplicities. 
\end{thm}

\begin{rk}
Again, the conditions (1) and (2) of the theorem are generic in the parameter space of all 
possible eigenvalues and parabolic weights. They will regularly appear along this paper, 
both in analytical and geometric arguments. 
\end{rk}

This theorem is a consequence of Theorem \ref{mainthmhiggs}. 
\begin{defn}
The map 
\begin{align} \label{nahm}
       \Nahm: (E,D) \mapsto (\hat{E},\hat{D})
\end{align} 
described in Theorems \ref{naivethm} and \ref{mainthm} will be called \emph{Nahm transform}. 
\end{defn}

Finally, as we have already mentioned, Nahm transform has an involutibility property:
\begin{thm}
Let $(E,D)$ be a parabolic integrable connection on the sphere satisfying the assumptions of Theorem
\ref{mainthm}. Then 
$$
    \Nahm^{2}(E,D) = (-1)^{\ast} (E,D),
$$
where $-1:\C \ra \C$ is the map $z \mapsto -z$, 
and $(-1)^{\ast}$ the induced map on fiber bundles with connection. In particular, Nahm transform is invertible. 
\end{thm}
This will be proved in Theorem \ref{invtr}, using arguments of the same type as S. K. Donaldson and P. B. Kronheimer in 
\cite{Don-Kronh}, namely the study of the spectral sequence of a suitable double complex. 

\section{Local model for parabolic integrable connections}\label{locmodint}
We suppose in this section that near each singularity, $h$ coincides with the diagonal 
models $h^{j}$ and $h^{\infty }$ given in Definition \ref{adaptdef} (that is, without the extra 
$O(|z-p_{j}|^{2 (\beta^{j}_{k} - \beta^{j}_{k'})})$ and $O(|z|^{-2 (\beta^{\infty}_{k} - \beta^{\infty}_{k'})})$ factors in 
(\ref{parmatr}) and (\ref{parmatrinf}); in particular, this metric is then not harmonic). 
For computations, it will be useful to express the local models of the 
integrable connection near the singularities in some orthonormal bases. As in \cite{Biq-Boa}, 
we consider the orthonormal basis defined by 
\begin{equation} \label{ebasisj}
      e^{j}_{k} = |z|^{-\beta^{j}_{k}-i\Im \mu^{j}_{k}} {\tau^{j}_{k}}   \qquad {k=1,\ldots,r}
\end{equation}
around $p_{j}$. The $h$-unitary part $(D^{+})^{j}$ of $D^{j}$ becomes  
\begin{align}
    (D^{+})^{j}  = 
    \mbox{d}+ i \Re (A^{j}) \d \theta  \label{polarunit} 
\end{align} 
where 
$\Re(A^{j}) = \frac{A^{j}+(A^{j})^{\ast}}{2}= diag(\Re \mu^{j}_{k})_{k=1,\ldots,r}$ (and 
$ \Im (A^{j}) = \frac{A^{j}-(A^{j})^{\ast}}{2i} = diag(\Im
\mu^{j}_{k})_{k=1,\ldots,r}$) stands for the real (imaginary) part of
$A^{j}$, and $r$ and $\theta $ are local polar coordinates around
$p_{j}$ such that we have  $z-p_{j}=r e^{i\theta}$. 
For the self-adjoint part $\Phi^{j}$ of $D^{j}$ in this basis we get: 
\begin{align}
    \Phi^{j} & = \frac{A^{j}}{2} \frac{\d z}{z-p_{j}} + \frac{(A^{j})^{\ast }}{2} \frac{\d
      \bar{z}}{\bar{z}-\bar{p_{j}}} - \beta^{j} \frac{\mbox{d}r}{r} \notag \\
      & = [ \Re(A^{j}) - \beta^{j} ] \frac{\mbox{d}r}{r} - \Im(A^{j}) \mbox{d}\theta \label{polarphi} , 
\end{align}
where $\beta^{j}=diag(\beta^{j}_{k})_{k=1,\ldots,r}$. These together imply that with
respect to this basis, the model for the operator $D$ in polar coordinates is 
\begin{eqnarray}
      D^{j} = \d + i A^{j} \d \theta + [ \Re (A^{j}) - \beta^{j} ] \frac{\mbox{d}r}{r}. \label{polard}
\end{eqnarray}

In an analogous way, in the orthonormal basis $  \left\{ e^{\infty}_{k} \right\}_{k=1,\ldots,r} $ given by 
\begin{equation} \label{ebasisinf}
      e^{\infty }_{k} = |z|^{\beta^{\infty }_{k}+i\Im \mu^{\infty }_{k}} \exp{[({\xi_{k}
          z}-{\bar{\xi_{k}}\bar{z}})/2]} \tau_{k}^{\infty}
\end{equation}
near infinity the unitary part of the model connection $D^{\infty }$ is given by   
\begin{eqnarray*}
         (D^{+})^{\infty} = \d + i \Re (C) \d \theta, 
\end{eqnarray*}
where we have put again 
$\Re (C)= \frac{C+C^{\ast }}{2} = diag (\Re \mu^{\infty}_{k})_{k=1,\ldots,r}$ and $z=r e^{i\theta }$. Moreover, putting 
$ \Re (z  A)=  diag \Re (\{ z \xi_{l}, m_{l} \})_{l=1,\ldots n'} $ and 
$\Im (z  A)= diag \Im  (\{ z \xi_{l}, m_{l} \} )_{l=1,\ldots n'} $, the self-adjoint part of $D^{\infty }$ has the form 
\begin{align} \label{polarphiinf}
        \Phi^{\infty } & = \frac{1}{2}\left( A + \frac{C}{z} \right) \d z + \frac{1}{2}\left( A^{\ast} +
          \frac{C^{\ast }}{\bar{z}} \right) \d \bar{z} + \beta^{\infty}
        \frac{\mbox{d}r}{r} \notag \\
        & = [ \Re (z  A + C) + \beta^{\infty}  ] \frac{\mbox{d}r}{r} +
        \Im ( z  A  + C)  \d \theta 
\end{align}
(the inversion of the sign of $\beta$ comes from the fact that if we make a
coordinate change $w=1/z$, $\vert w \vert = \rho = 1/r = 1/ \vert z \vert$,
then $\d \rho / \rho = - \d r / r$). Remark that in these expressions the terms 
in $\d \theta, \d r / r , \d z / z, \d \bar{z}/ \bar{z}$ are of lower order then
the ones in $\d z, \d \bar{z}, z \d r /r, z \d \theta $; hence the leading order
term of the singular part of $D$ in this basis is just 
\begin{eqnarray} \label{singpart}
                  d + \frac{A}{2} \d z + \frac{A^{\ast }}{2} \d \bar{z}.
\end{eqnarray}

\section{Higgs bundle point of view}\label{higgsintro}
The idea of the proofs of Theorems \ref{naivethm} and \ref{mainthm} will be to exploit the 
correspondence known as nonabelian Hodge theory between parabolic integrable connections 
on one side and parabolic Higgs bundles on the other side. 
Let us recall the definition of the latter notion: 
\begin{defn} \label{parhiggsdef}
A \emph{parabolic Higgs bundle} is given by:  
\begin{enumerate}
\item{a holomorphic bundle $\E$ with holomorphic structure $\bar{\partial}^{\E}$ over $\CP$ 
      called the \emph{holomorphic bundle underlying the Higgs bundle}, 
      and with underlying smooth vector bundle $V$;}
\item{in each point $p \in P \cup \{\infty \}$ a strictly decreasing parabolic flag 
$$
     V_{p} =  F_{0}V_{p} \supset F_{1}V_{p} \supset \ldots \supset F_{c_{p}-1} V_{p} \supset F_{c_{p}} V_{p} = \{ 0 \}
$$
      for some $1\leq c_{p}\leq r$, with parabolic weights 
$$
     0\leq \tilde{\alpha}^{p}_{1} <\ldots<\tilde{\alpha}^{p}_{c_{p}}<1;
$$}
\item{a $\bar{\partial}^{\E}$-meromorphic section $\theta \in \Omega^{1,0}(\CP, End(V))$ (called the \emph{Higgs field}), 
     having a simple pole at the points of $P$ with semi-simple residue respecting the parabolic flag 
     (that is, such that $Res(\theta,p_{j})$ leaves $F_{k}V_{p_{j}}$ invariant for each $p_{j}\in P$ and 
     all $0 \leq k \leq c_{p}$), and a second-order pole at infinity, such that there exists a holomorphic 
     basis of $\E$ near infinity compatible with the parabolic structure in which the residue 
     and second-order term are both diagonal.}
\end{enumerate}
Again, we write 
$$
     0\leq {\alpha}^{p}_{1} \leq \ldots\leq {\alpha}^{p}_{r}<1
$$
for the parabolic weights repeated according to their multiplicities 
$$
    \dim (F_{k-1}V_{p}/F_{k}V_{p}),
$$
and we shorten ${\alpha}^{p_{j}}_{k}$ to ${\alpha}^{j}_{k}$. Finally, we set 
\begin{equation} \label{dsecond}
       D''=\dbar^{\E}+ \theta,
\end{equation} 
that we call the \emph{$D''$-operator} associated to the Higgs bundle. 
\end{defn}
The notions of parabolic degree, slope and (poly/semi-)stability of parabolic Higgs bundles are defined 
analogously to the case of integrable connections, see Definition \ref{parintcon}. 
Theorem 6.1 of  \cite{Biq-Boa} then says the following. 
\begin{thm}[Biquard O.-Boalch Ph., 2004] \label{wnht}
There exists an isomorphism between the moduli space of parabolically stable rank $r$ integrable 
connections with fixed diagonal polar part and parabolic structures up to complex holomorphic
gauge transformations respecting the parabolic flags, and the moduli space of parabolically stable rank $r$
Higgs bundles with fixed diagonal polar part and parabolic structures up to complex holomorphic gauge 
transformations respecting the parabolic flags.
\end{thm} 
\begin{rk}
Actually, this is a consequence of the existence of a harmonic metric (Theorem \ref{harmexistthm}), 
and hence also proved for parabolic integrable 
connections with poles of arbitrary fixed order and diagonal polar part in the punctures and parabolic Higgs 
bundles with poles of the same order with diagonal polar part. 
\end{rk}

The transition from integrable connections to Higgs bundles is given as follows: first, 
the underlying smooth vector bundle of the integrable connection and the Higgs bundle are the same. 
Furthermore, suppose $h$ is the harmonic metric, consider the decomposition (\ref{decompsaunit}) 
of the integrable connection into its unitary and self-adjoint part, and decompose the terms further 
according to bidegree 
\begin{align} \label{decompbidegree}
           D^{+} & = (D^{+})^{1,0} + (D^{+})^{0,1} \\
           \Phi & = \Phi^{1,0} + \Phi^{0,1}. \notag
\end{align}
The partial connection $(D^{+})^{0,1}$ defines then the holomorphic structure of $\E$, and $\Phi^{1,0}$ will 
be the Higgs field $\theta$. The $D''$-operator is of course $(D^{+})^{0,1}+\Phi^{1,0}$.
Harmonicity of the metric implies that $\theta$ is holomorphic. 

The transition in the other direction is also established using a privileged metric. 
\begin{defn} \label{HE}
Let $(\E,\theta)$ be a Higgs bundle. We say that $h$ is a \emph{Hermitian-Einstein metric} for $(\E,\theta)$ if, 
denoting by $D^{+}_{h}$ the Chern connection (the unique $h$-unitary 
connection compatible with $\dbar^{\E}$), by $F_{D^{+}_{h}}$ its curvature, and by 
$\theta^{\ast}_{h}$ the adjoint of $\theta$ with respect to $h$, then these objects satisfy the real Hitchin equation
\begin{eqnarray*}
              F_{D^{+}_{h}}+[\theta,\theta^{\ast}_{h}] = 0, 
\end{eqnarray*}
where $[.,.]$ stands for graded commutator of forms.
\end{defn} 
Let $(\E,\theta)$ be a parabolically stable parabolic Higgs bundle. By \cite{Biq-Boa}, there exists 
a unique Hermitian-Einstein metric $h$ adapted to the parabolic structure. The connection 
\begin{equation}
         D =  D^{+}_{h} + ( \theta +  \theta_{h}^{\ast } )
\end{equation}
on $V$ is then integrable, and $h$ is the corresponding harmonic metric adapted to the parabolic structure. 
In what follows, in order to simplify notations, we are often going to omit the subscript $h$ in the 
notation of the Chern connection and adjoints.

Let now $(E,D)$ be a parabolically stable parabolic integrable connection and $(\E,\theta)$ the associated  
parabolic Higgs bundle. An important result we will be constantly using is the following 
\begin{thm}[Simpson C. \cite{Sim}] \label{simpson}  
Suppose the metric $h$ is harmonic. Then, with the previous notations, the Laplace
operators $\Delta_{D}=DD^{\ast}+D^{\ast}D$ and $\Delta_{D''}=D''(D'')^{\ast}+(D'')^{\ast}D''$ satisfy
$$
      \Delta_{D}=2\Delta_{D''}.
$$
In particular, their domain and kernel coincide. 
\end{thm}

\section{Local model for Higgs bundles}
In this section, we give the eigenvalues of the residue of the Higgs field and the parabolic weights of the 
Higgs bundle in the punctures that correspond to those of the integrable connection via the Theorem \ref{wnht}. 
To obtain local models for the operators in this setup, suppose again that
near $p_{j}$ the metric $h$ coincides with the diagonal model $h^{j}$ given by (\ref{parmatr}) 
(without the correcting $O(|z-p_{j}|^{2(\beta^{j}_{k}-\beta^{j}_{k'})})$ term; in particular, it does not satisfy 
Hitchin's equation). Then, according to \cite{Biq-Boa}, in the local $\dbar^{\E}$-holomorphic trivialisation  
\begin{align}\label{holtriv}
     \sigma^{j}_{k} = |z-p_{j}|^{\Re \mu^{j}_{k}} \frac{e^{j}_{k}}{(z-p_{j})^{[\Re \mu^{j}_{k}]}}
     && (k=1,\ldots,r)
\end{align}
around $p_{j}$, the Higgs field is equal up to a perturbation term to the model Higgs field given by 
\begin{align} \label{modelhiggs}
     \theta^{j} & = \frac{A^{j}-\beta^{j}}{2} \frac{\d z}{z-p_{j}}\notag \\
     & = diag\left( \frac{\mu^{j}_{k}-\beta^{j}_{k}}{2} \frac{\d z}{z-p_{j}} \right)_{k=1,\ldots,r}\notag \\
     & = diag\left(\lambda^{j}_{k} \frac{\d z}{z-p_{j}} \right)_{k=1,\ldots,r} ,
\end{align}
where we have put $\lambda^{j}_{k}=(\mu^{j}_{k}-\beta^{j}_{k})/2$. 
Moreover, in the same trivialisation, the parabolic weights are 
\begin{equation} \label{higgsparweights}
    \alpha^{j}_{k}=\Re (\mu^{j}_{k}) -[\Re (\mu^{j}_{k})],
\end{equation}
where $[.]$ denotes integer part. 
\begin{rk} 
In fact, this formula is not completely correct, because the $\alpha^{j}_{k}$ defined by it are not 
necessarily in increasing order, although they should be by definition. 
One should instead write the same formula for $\alpha^{j}_{s(k)}$, 
where $s$ is a permutation of $\{ 1,\ldots,r\}$. However, in the sequel we discard this minor technical 
detail for the sake of simplicity of the notation. 
\end{rk}
\begin{rk} Since the gauge transformations between the bases 
$\{ \tau^{j}_{k} \}_{k=1,\ldots,r}$ and $\{ \sigma^{j}_{k} \}_{k=1,\ldots,r}$ 
are just multiplications by some functions (in particular diagonal
matrices), it follows that the smooth subbundle
spanned by the sections $\{ \sigma^{j}_{k} \}_{k=r_{j}+1,\ldots,r}$ is the same as
the one spanned by $\{ \tau^{j}_{k} \}_{k=r_{j}+1,\ldots,r}$, which is by
definition the underlying smooth vector bundle of $E^{j}_{sing}$; and similarly, the subbundle spanned by 
$\{\sigma^{j}_{k} \}_{k=1,\ldots,r_{j}}$ is equal to the underlying smooth bundle of $E^{j}_{reg}$. The same remark also
holds for the bases $\{ e^{j}_{k} \}$ instead of $\{ \sigma^{j}_{k} \}$. 
In particular, the residue of the model Higgs field $\theta^{j} $ in the
point $p_{j} \in P$ belongs to $End(E^{j}_{sing}|_{p_{j}})$. 
\end{rk}

Near infinity, the situation is slightly different: for $h = h^{\infty }$ the diagonal model, 
in the local $\dbar^{\E}$-holomorphic frame 
\begin{align} \label{holtrivinf}
            \sigma_{k}^{\infty} = |z|^{-\Re \mu^{\infty }_{k}} z^{[\Re \mu^{\infty }_{k}]} e^{\infty }_{k} && (k=1,\ldots,r)
\end{align}
the Higgs field is equal up to a perturbation term to the model Higgs field  given by  
\begin{align} \label{modelhiggsinf}
      \theta^{\infty } & = \frac{1}{2} A\d z + \frac{\mu^{\infty} - \beta^{\infty}}{2} \frac{\d z}{z}\notag \\
      & = \left( \frac{1}{2} diag(\{ \xi_{l}, m_{l} \})_{l=1,\ldots,n'} 
      + \frac{1}{z} diag(\lambda_{k}^{\infty})_{k=1,\ldots,r} \right) \d z,
\end{align}
where we have put again $\lambda_{k}^{\infty}=(\mu_{k}^{\infty}-\beta_{k}^{\infty})/2$, with parabolic weights
being, as in  the case of simple poles, 
\begin{equation} \label{higgsparweightsinf}
\alpha_{k}^{\infty}=\Re (\mu_{k}^{\infty}) -[\Re (\mu_{k}^{\infty})]. 
\end{equation}
From these data, as above, one can form the model $D''$-operator 
\begin{equation} \label{modeldsecond}
   (D'')^{j} = \dbar^{\E} + \theta^{j} \qquad (j \in \{ 1, \ldots n, \infty \} ). 
\end{equation}
Notice that since we considered holomorphic trivialisations of $\E^{j}$, the
partial connection part of the model coincides with the usual $\dbar$-operator. 

We are now ready to write out the assumptions made in Theorem \ref{mainthm} on the parameters of the 
integrable connection, translated to those of the Higgs bundle: 

\begin{hyp} \label{main}
We suppose that $(\E ,\theta )$ is a parabolically stable Higgs bundle with diagonal polar part of the 
Higgs field in some local holomorphic frame near each puncture, satisfying the properties 
\begin{enumerate} 
\item{for fixed $j\in \{ 1,\ldots,n  \}$ the residues $\lambda^{j}_{k} $ for $k \in \{ r_{j}+1,\ldots,r \} $ are
    non-vanishing and distinct, $\lambda^{j}_{k}$ vanish for $k=1,\ldots,r_{j}$ and finally $\alpha_{k}^{j} \neq 0$ if and only
    if $\lambda^{j}_{k} \neq 0$;}\label{maini}
\item{for fixed $l \in \{ 1, \ldots, n' \} $ the complex numbers $\lambda^{\infty }_{k}$ for $k \in \{ 1+a_{l}, \ldots a_{l+1} \}$ are
    non-vanishing and distinct, and $\alpha^{\infty }_{k} \neq 0$.} \label{mainii}
\end{enumerate}
\end{hyp}
Diagonality of the polar parts has already been assumed when writing the local models (\ref{modelhiggs}) and
(\ref{modelhiggsinf}). The first condition says that no parabolic weight and no eigenvalue of the residue of $\theta$ 
vanishes on the singular component at any singularity, and that on the singular component near a puncture all 
eigenvalues are different; whereas the eigenvalues of the residue and parabolic weights vanish on the regular 
component. One more way to say the same thing is: for all $j \in \{ 1,\ldots n \}$, the residue of $\theta $ defines an
automorphism of $E^{j}_{sing}|_{p_{j}}$, and the parabolic weights corresponding to the holomorphic 
trivialisation (\ref{holtriv}) are non-vanishing exactly on this subspace.
The second one imposes that on the eigenspace corresponding to a fixed eigenvalue of the second-order term 
at infinity, all the eigenvalues of the residue be non-vanishing and distinct, furthermore that no parabolic 
weight vanish at infinity. Note that these conditions are generic among all possible choices of
singularity parameters. 

\section{The transformation of the Higgs bundle} \label{transf}
Let $(\E ,\theta )$ be a parabolic Higgs bundle and $\xi \in \hat{\C} \setminus \hat{P}$ a parameter. 
The natural deformation of the Higgs field is 
\begin{align} \label{deformhiggs}
        \theta_{\xi} &= \theta - \frac{\xi}{2}  \mbox{d}z 
\end{align}
with fixed underlying holomorphic bundle $\E$. It is clear that $\theta_{\xi}$ is then also holomorphic with respect 
to $\dbar^{\E}$ with the same local models at the logarithmic punctures as $\theta$, but its local model near 
infinity is different. If moreover a Hermitian metric is fixed, then we also have 
$$
       \theta_{\xi}^{\ast} = \theta^{\ast } - \frac{\bar{\xi}}{2} \mbox{d}\bar{z}. 
$$
Therefore, the integrable connection corresponding to the deformed Higgs bundle is given by 
\begin{eqnarray} \label{higgsdef}
       D^{H}_{\xi} = D - \frac{\xi}{2} \mbox{d} z- \frac{\bar{\xi}}{2} \mbox{d}\bar{z},
\end{eqnarray} 
and the crucial observation is that via the unitary gauge transformation 
\begin{equation} \label{gauge}
     \exp[ ( \bar{\xi } \bar{z} - \xi z )/2]
\end{equation}
on $\C$ this is equivalent to the deformation (\ref{intdef}). The self-dual part of this deformation is  
\begin{align} 
     \Phi^{H}_{\xi } & = \Phi - \frac{{\xi}}{2} \d z - \frac{\bar{\xi}}{2} \d \bar{z} \label{higgsdefsa},  
\end{align} 
the same deformation as in (\ref{intdefsa}). Therefore it will not make any
confusion to refer to $\Phi_{\xi }$ without mentioning the adopted point of
view; consequently, we drop the corresponding upper indices. The connection defined by (\ref{higgsdef})
is still flat, but the underlying holomorphic structure is different from
the one of $D$ (because of the term in $\mbox{d}\bar{z}$). Notice also that the gauge transformation 
(\ref{gauge}) between these deformations has an exponential singularity at infinity. Denote by  $\EC^{H}_{\xi}$ the 
elliptic complex
\begin{equation}\label{merellcompl}
        \Omega^{0}\otimes E \xrightarrow{\text{${D}^{H}_{\xi}$}}  \Omega^1 \otimes E \xrightarrow{\text{${D}^{H}_{\xi}$}} \Omega^2 \otimes E. 
\end{equation}
\begin{defn}
The smooth vector bundle $\hat{V}$ underlying the transformed Higgs bundle is the vector bundle 
whose fiber over $\xi \in \hat{\C} \setminus \hat{P}$ is the first $L^{2}$-cohomology $L^{2}H^{1}(\EC^{H}_{\xi})$
of $\EC^{H}_{\xi}$.
\end{defn}
In Proposition \ref{trvbprop} we prove that these vector spaces indeed define a finite rank smooth bundle. 
Furthermore, by Theorem \ref{laplker}, any class in $L^{2}H^{1}(\EC^{H}_{\xi})$ admits a unique 
${D}^{H}_{\xi}$-harmonic representative. 
\begin{defn}
The \emph{transformed holomorphic structure} on $\hat{V}$ is the one induced by the orthogonal projection
$\dbar^{\E}$ of the trivial partial connection with respect to the variable $\xi$ on the trivial $L^{2}$-bundle 
over $\CPt$ to $D^{H}_{\xi}$-harmonic $1$-forms. 
The \emph{transformed Higgs field} is multiplication by $-z\d \xi /2$ followed by projection onto harmonic 
$1$-forms. Finally, the \emph{transformed Hermitian metric} is the $L^{2}$-metric of the harmonic
representative. 
\end{defn}
By virtue of Theorems \ref{laplker} and \ref{simpson}, the transformed smooth bundle $\hat{V}$ can also be
computed in this case as the first cohomology of the elliptic complex $\EC_{\xi}''$ given by:
$$
     \Omega^{0}\otimes E \xrightarrow{\text{$D_{\xi}''$}}  \Omega^1 \otimes E
     \xrightarrow{\text{$D_{\xi}''$}} \Omega^2 \otimes E, 
$$ 
where the maps are the corresponding deformations of (\ref{dsecond}) in the
Higgs-bundle point of view. Explicitly, $D_{\xi}''$ reads 
$$
      (D^{H}_{\xi })''= \dbar^{\E} + \theta_{\xi }. 
$$
We use this description of the transformed bundle in Section \ref{sectharm} to show the statement 
of Theorem \ref{naivethm} on the transformed metric:
\begin{thm} 
If the original metric harmonic then the same thing holds for the transformed metric. 
\end{thm}
For this purpose, we prove in fact that the candidate Higgs field $\hat{\theta}$ corresponding to 
$\hat{D}$ and $\hat{h}$ is meromorphic with respect to the transformed holomorphic structure.

Furthermore, in this interpretation, the remaining part of Theorems \ref{naivethm} and \ref{mainthm} can be written: 
\begin{thm} \label{mainthmhiggs} Suppose $(\E,\theta )$ is a parabolic Higgs bundle with logarithmic singularities 
in the points of $P$ and a double pole at infinity, as described in Section \ref{higgsintro}, such that 
its singularity parameters satisfy Hypothesis \ref{main}. Then the transformed Higgs bundle 
$(\dbar^{\hat{\E}},\hat{\theta} )$ is of the same type (that is, it has a finite number of logarithmic
singularities in points of $\hat{\C}$ and a double pole at infinity, with a parabolic structure in these 
points). Furthermore, its topological and singularity parameters are as follows: 
\begin{enumerate}
\item{the rank of $\hat{\E}$ is the sum (\ref{trrank}) of the ranks of the residues of $\theta $ in $P$} \label{i}
\item{its degree is the same as that of $\E$} \label{ii}
\item{the logarithmic singularities are located in the set $\hat{P} $, and for all $l \in \{ 1,\ldots,n' \}$ 
    the rank of the transformed Higgs field in the point $\xi_{l }$ is equal to the multiplicity $m_{l}$ 
    of the eigenvalue $\xi_{l }$ of $A$} \label{iii}
\item{the set of non-vanishing eigenvalues of the residue of $\hat{\theta }$ in the point $\xi_{l }$ is  
$\{ -\lambda^{\infty }_{1+a_{l}}, \ldots , -\lambda^{\infty }_{a_{l+1}} \}$, where $\{ \lambda^{\infty }_{a_{l}+1}, \ldots, \lambda^{\infty }_{a_{l+1}} \}$ are the 
eigenvalues  of the residue of the original Higgs field $\theta $ at infinity, 
restricted to the eigenspace of $A$ corresponding to the eigenvalue $\xi_{l }$ } \label{iv}
\item{the non-vanishing parabolic weights of $\hat{\E }$ in $\xi_{l }$ is the set of parabolic weights 
$\{ \alpha^{\infty }_{1+a_{l }},\ldots,\alpha^{\infty }_{a_{l+1 }} \}$ of $\E $ at infinity, restricted to the same subspace} \label{v}
\item{the eigenvalues of the second-order term of $\hat{\theta} $ at infinity are $\{ -p_{1}/2,\ldots,-p_{n}/2 \}$, and
  the multiplicity of $-p_{j}/2$ is equal to the rank $r-r_{j}$ of the residue of $\theta $ in $p_{j}$} \label{vi} 
\item{on the eigenspace corresponding to $-p_{j}/2$ of the second-order term at infinity,
the eigenvalues of the residue of $\hat{\theta} $ are \newline $\{ -\lambda^{j}_{r_{j}+1},\ldots,- \lambda^{j}_{r} \}$} \label{vii}
\item{the parabolic weights on the same eigenspace at infinity are the parabolic weights 
$\{ \alpha^{j}_{r_{j}+1},\ldots, \alpha^{j}_{r} \}$ of $\E $ at $p_{j}$} \label{viii}
\end{enumerate}
\end{thm} 

The proof of this theorem is the object of Chapter \ref{algint}. 



\chapter{Analysis of the Dirac operator} \label{chFr}
In this chapter, we study the analytical theory needed for our construction along the lines of 
Donaldson-Kronheimer \cite{Don-Kronh}, Jardim \cite{Jardim} and others. First, in Section \ref{sectFr} we
define spinor spaces and Dirac operators that permit us to study the problem. 
We also define a suitable functional space $H^{1}$ and state a Fredholm theorem 
valid for all deformations of the initial connection. 
Then it is natural to define the fibers of the transformed bundle 
as the cokernel of the deformed Dirac operator. The Fredholm theorem 
is then proved in Section \ref{sectpfFr}. 
In Section \ref{l2coh}, we carry out an identification of this cokernel with the first $L^{2}$-cohomology 
$L^{2}H^{1}(D^{int}_{\xi})$ of the  complex $\EC_{\xi}^{int}$  given in (\ref{ellcompl}), similar in vein to
the Hodge isomorphism between the cokernel of the operator 
$\d + \d^{\ast }$ on a compact manifold and the $L^{2}$-cohomology of the operator $\d$. 
However, since the manifold we are working on is non-compact, in proving these results we need a careful 
study of the singularities. 

In all what follows, we fix a parabolic integrable connection with adapted metric $(E,D,h)$ and choose to study the analytic
properties of the deformation from the point of view of integrable connections, hence we set 
for simplicity $D_{\xi }= D^{int}_{\xi }$ until further notification. 

\section{Statement of the Fredholm theorem} \label{sectFr}
\begin{defn}
The positive and negative \emph{spinor bundles} are the vector bundles over $\C \setminus P$ given by 
\begin{align*}
    S^{+} &= \Lambda^{0}T^{\ast }(\C \setminus P) \oplus \Lambda^{2}T^{\ast }(\C \setminus P) & S^{-} = \Lambda^{1}T^{\ast }(\C \setminus P)
\end{align*}
\end{defn}
Recall that we have defined $\hat{P}$ as the set $\{ \xi_{1},\ldots,\xi_{r} \}$ of all
eigenvalues of the second order term of $D $ at infinity.
\begin{defn}
For $\xi \in \hat{\C} \setminus \hat{P}$ the \emph{Dirac operator } is the first-order differential operator
\begin{align*}
              \Dir_{\xi} &=D_{\xi} - D_{\xi}^{\ast}: \Gamma (S^{+}\otimes E) \lra  \Gamma (S^{-}\otimes E)
\end{align*}
where $\Gamma$ is used to denote smooth sections with compact support in $\C \setminus P$. 
Its formal adjoint 
\begin{align*}
              \Dir_{\xi}^{\ast} & = D_{\xi}^{\ast} - D_{\xi} : \Gamma (S^{-}\otimes E) \lra \Gamma
              (S^{+}\otimes E),
\end{align*}
is called the \emph{adjoint Dirac operator. }
\end{defn}
For any $\xi\in \hat{\C}$ let us introduce the following norm on sections $f$ of $S^{+}\otimes E$: 
\begin{align} \label{norm}
     \Vert f \Vert_{H^{1}_{\xi}}^{2} = \int_{\C } &  \vert f
     \vert^{2} + \vert \nabla_{\xi }^{+}f \vert^{2} + \vert \Phi_{\xi} \otimes f \vert^{2} , 
\end{align}
where $\nabla_{\xi }^{+}$ and $\Phi_{\xi}$ are defined in (\ref{intdefunit}) and (\ref{intdefsa}).
Here and in all what follows, we integrate with respect to the
Euclidean volume form $|\d z|^{2}$, and $\vert x \vert^{2}$ denotes $h(x,x)$, unless the
contrary is explicitly stated. Our convention is furthermore to write
$(x,y)$ for $h(x,y)$, and for sections $x$ and $y$, we write $\langle x,y\rangle $
instead of $\int_{\C } (x,y) $. 

Define the space of sections  
\begin{equation} \label{h1space}
       H^{1}_{\xi}(S^{+}\otimes E) = \{ f \in L^{2}_{loc}(S^{+}\otimes E): \Vert f \Vert_{H^{1}_{\xi}} <
       \infty \},
\end{equation} 
where in $L^{2}$ we refer to the metric $h$ on the fibers. We will often
write $H^{1}_{\xi}$ instead of $H^{1}_{\xi}(S^{+}\otimes E)$. 
As we will see by the end of this chapter, this is the appropriate space to regard the Dirac
operators. First we establish the simple 
\begin{lem} \label{lemmaone} The norm $\Vert . \Vert_{H^{1}_{\xi}}$ 
  depends (up to equivalence of norms) neither on $\xi \in \hat{\C}$, 
  nor on the particular connection $D$ having behavior as in 
  (\ref{dpj}) and (\ref{dinf}). 
\end{lem}
\begin{proof} We begin by showing that the norm is independent of $\xi $. 
In order to simplify notations, we let $H^{1}$ stand for $H^{1}_{0}$ from now on. 
It is obviously sufficient to prove that for an arbitrary
$\xi \in \hat{\C}$, the $H^{1}_{\xi}$-norm is equivalent to the $H^{1}$-norm. From the point-wise identity 
$$
     |\Phi_{\xi } \otimes f| = 2 |\theta_{\xi } \otimes f| =  2 |\theta_{\xi }^{\ast } \otimes f| ,
$$
and the point-wise estimation 
\begin{align} \label{punctestunit}
   |\nabla^{+}_{\xi }f| & \leq  |\nabla^{+} f| + 2|\xi || f| ,
\end{align}
one sees that for any section $f=(f_{0},f_{2}) \in \Gamma(S^{+} \otimes E)$ the estimates  
$$
     \Vert f \Vert^{2}_{H^{1}_{\xi}} \leq (1+ 8\vert \xi \vert^{2}) \Vert f \Vert^{2}_{H^{1}}
$$
and 
$$
     \Vert f \Vert^{2}_{H^{1}} \leq (1+ 8\vert \xi \vert^{2}) \Vert f \Vert^{2}_{H^{1}_{\xi}}
$$
hold; the first statement of  the Lemma follows at once. 

Now we show independence of the particular connection $D$ with right
singularity behavior. Introduce the model norm 
\begin{eqnarray} \label{modelnorm}
     \Vert f \Vert^{2}_{H^{1}_{mod}(\Delta(p_{j}, \varepsilon))} = 
      \int_{\Delta(p_{j}, \varepsilon)}  \vert f \vert^{2} + \vert (D^{+})^{j}f \vert^{2} + 
      \vert \Phi^{j} f \vert^{2}
\end{eqnarray}
around points of $P$ and the model norm 
\begin{eqnarray} \label{modelnorminf}
   \Vert f \Vert^{2}_{H^{1}_{mod}(\C \setminus \Delta(0, R))} = \int_{\C \setminus \Delta(0,R)}   \vert f \vert^{2} + 
   \vert (D^{+})^{\infty } f \vert^{2} + \vert \Phi^{\infty} f \vert^{2}   
\end{eqnarray}
near infinity. 
Then it is sufficient to prove the following:
\begin{clm} \label{equivclm}
If $\varepsilon >0$ is chosen sufficiently small and $R>0$ sufficiently
large, then for any smooth section $f \in H^{1}$ we have 
\begin{eqnarray} \label{equiv}
     c \Vert f^{j} \Vert^{2}_{H^{1}(\Delta(p_{j}, \varepsilon))} < 
     \Vert f^{j} \Vert^{2}_{H^{1}_{mod}(\Delta(p_{j}, \varepsilon))}
     <  C \Vert f^{j} \Vert^{2}_{H^{1}(\Delta(p_{j}, \varepsilon))}
\end{eqnarray}
and similarly 
\begin{eqnarray} \label{equivinf}
     c \Vert f^{j} \Vert^{2}_{H^{1}(\C \setminus \Delta(0, R))} < 
     \Vert f^{j} \Vert^{2}_{H^{1}_{mod}(\C \setminus \Delta(0, R))}
     <  C \Vert f^{j} \Vert^{2}_{H^{1}(\C \setminus \Delta(0, R))}
\end{eqnarray}
with some constants $0<c<C$ independent of $f$. 
\end{clm}
\begin{proof}
Consider first the case of $p_{j} \in P$. Decompose $f^{j}=f^{j}_{reg}+f^{j}_{sing}$ corresponding to the 
splitting (\ref{regsing}). Write also 
\begin{align}
    f^{j}_{reg} & = \sum_{k=1}^{r_{j}} \phi^{j}_{k} e^{j}_{k} \label{regdecomp} \\
    f^{j}_{sing} & = \sum_{k=r_{j}+1}^{r} \phi^{j}_{k} e^{j}_{k} \label{singdecomp}
\end{align}
with respect to the orthonormal basis $\{ e^{j}_{k} \}$ introduced in (\ref{ebasisj}), 
where the $\phi^{j}_{k}$ are functions. 
Formulas (\ref{polarunit}) and (\ref{polarphi}) and Hypothesis \ref{main} imply that (\ref{modelnorm}) 
is equivalent to the weighted Sobolev space of sections satisfying 
\begin{align} \label{weighted}
       \sum_{k=1}^{r_{j}} \int_{\Delta(p_{j}, \varepsilon)} & \vert \phi^{j}_{k} \vert^{2} + \vert \d \phi^{j}_{k} \vert^{2} \\
       & + \sum_{k=r_{j}+1}^{r}  \int_{\Delta(p_{j}, \varepsilon)} 
      \absl \frac{\phi^{j}_{k}}{r} \absr^{2} + \vert \d \phi^{j}_{k} \vert^{2}  < \infty ,\notag
\end{align}
where $\d $ stands for the trivial connection on functions. Notice that we only add 
weights on the singular component. 
By \cite{Sim}, Theorem 1 it follows that in $\Delta(p_{j}, \varepsilon)$ the difference between 
$(D^{+})^{j}$ and $D^{+}$ is 
\begin{equation} \label{pertj}
     a^{j} = O(r^{-1+\delta} )
\end{equation} 
for some $\delta >0$, and the same estimation holds for the difference between $\Phi^{j}$ and $\Phi$. 
It is then immediate that for any $c>0$, the estimation 
$$
     \int_{\Delta(p_{j}, \varepsilon)} \absl \frac{\phi^{j}_{k}}{r} \absr^{2} > c
     \int_{\Delta(p_{j}, \varepsilon)} | a^{j} \phi^{j}_{k} |^{2}  
$$
holds for $k =r_{j}+1,\ldots,r$  and for $\varepsilon >0$ sufficiently small. 
We therefore have (\ref{equiv}) for $f_{sing}$. 
On the other hand, for a function $g$ defined in $\Delta(0,1)$ and for $\delta >0$ fixed, from the claim in
the proof of Theorem 5.4 in \cite{Biq-Boa} we have 
$$
    \int_{\Delta(p_{j}, 1)} | r^{-1+\delta }g |^{2} \leq c\left( \int_{\Delta(p_{j}, 1)} | \d g |^{2} + 
    \int_{\Delta(p_{j}, 1) \setminus \Delta(p_{j}, 1/2) } |g|^{2} \right). 
$$
Rescaling this inequality to the disk $ \Delta(p_{j}, \varepsilon) $ one easily checks
that it implies 
\begin{align} \label{regineq}
    \varepsilon^{-2\delta } \int_{\Delta(p_{j}, \varepsilon)} & | r^{-1+\delta } g|^{2} \notag \\ 
    \leq c & \left( \int_{\Delta(p_{j}, \varepsilon)} | \d g |^{2} + \varepsilon^{-2} 
    \int_{\Delta(p_{j}, \varepsilon) \setminus \Delta(p_{j}, \varepsilon/2) } |g|^{2} \right).
\end{align} 
Choosing $\varepsilon$ sufficiently small, applying this to $\phi^{j}_{k}$ for $k=1,\ldots,r_{j}$, and 
recalling that on the regular component $(D^{+})^{j}$ is the trivial connection $\d $ and 
$\Phi^{j} =0$, we obtain (\ref{equiv}) for $f_{reg}$ as well. This establishes the
equivalence of the norms $\Vert . \Vert^{2}_{H^{1}_{mod}}$ and 
$ \Vert . \Vert^{2}_{H^{1}} $ around a finite singularity. 

Around infinity, by \cite{Biq-Boa} Lemma 4.6 the difference between $(D^{+})^{\infty}$ and $D^{+}$ 
is bounded above by a term 
\begin{equation} \label{pertinf}
         a^{\infty } = O(r^{-1 - \delta })
\end{equation}
for some $\delta >0$, and again the same holds for $\Phi^{\infty} - \Phi$.  
The equivalence (\ref{equivinf}) follows immediately from the estimation 
$$
     |r^{-1 - \delta }f|\leq c |f|
$$
for any $c>0$, whenever $r>R$ with $R$ sufficiently large. 
\end{proof}

This then finishes the proof of Lemma \ref{lemmaone} as well.  

\end{proof}

From the previous discussion, we bring out as consequence: 

\begin{cor} \label{weightedcor}
 The Hilbert space $H^{1}(E)$ is the set of sections $f \in L^{2,1}_{loc}(E)$ such that near a logarithmic 
 singularity $p_{j}$, in the decompositions (\ref{regdecomp}) and (\ref{singdecomp}) we have 
 $\phi^{j}_{k} \in L^{2,1}$ for $k=1,\ldots,r_{j}$ and $\phi^{j}_{k}/r, \d \phi^{j}_{k} \in L^{2}$ for $k=r_{j}+1,\ldots,r$; 
 whereas at infinity, the coordinates $\phi^{\infty }_{k}$ of $f$ in the basis (\ref{ebasisinf}) are $L^{2,1}$; 
 equipped with the norm 
\begin{align*}
     & \int_{\C \setminus \cup_{j} \Delta(p_{j}, \varepsilon)} |f|^{2} + |\nabla f|^{2} \\ 
     & + \sum_{j=1}^{n} \left\{ \sum_{k=1}^{r_{j}} \int_{\Delta(p_{j}, \varepsilon)}  \vert \phi^{j}_{k} \vert^{2}  + 
       \vert \d \phi^{j}_{k} \vert^{2} + \sum_{k=r_{j}+1}^{r}  \int_{\Delta(p_{j}, \varepsilon)} 
      \absl \frac{\phi^{j}_{k}}{r} \absr^{2} + \vert \d \phi^{j}_{k} \vert^{2} \right\}
\end{align*}
The same result holds for sections of $\Omega^{2} \otimes E$, coordinates being expressed in the basis 
$\d z \land \d \bar{z}$. 
\end{cor}
\begin{proof}
For sections of $\Omega^{0}$, this follows from Claim \ref{equivclm}, (\ref{weighted}) and 
$$
       |\Phi \otimes f| \leq K|f|.
$$
We then obtain the case of $\Omega^{2}$ by duality.  
\end{proof}

We now come back to the analysis of the Dirac operator. 
From the definitions of $H^{1}( S^{+} \otimes E )$ and $\Dir_{\xi}$ we see that
this latter admits a bounded extension 
\begin{eqnarray} \label{Dirac}
    \Dir_{\xi}:  H^{1}( S^{+} \otimes E ) \lra L^{2}(S^{-} \otimes E). 
\end{eqnarray}
We are now able to announce the first main result of this chapter:

\begin{thm} \label{Fredholm} 
The operator (\ref{Dirac}) is Fredholm; if $h$ is harmonic, its kernel vanishes. 
\end{thm}

\begin{cor} \label{cokerbundle}
The bundle over $\hat{C} \setminus \hat{P}$ whose fiber over $\xi$ is the cokernel of (\ref{Dirac}) is 
a smooth vector bundle. 
\end{cor}
\begin{proof}
We recall the well-known fact that the index of a continuous family of Fredholm operators is constant. 
On the other hand, if the kernel of a Fredholm operator vanishes, then its index is equal to the 
opposite of the dimension of its cokernel. It then follows immediately from the Fredholm theorem that if 
the metric is harmonic, then the dimension of the cokernel of the operator $\Dir_{\xi }$ is a finite 
constant independent of $\xi $. Moreover, by standard implicit function theorem arguments in Hilbert space 
it follows that the cokernels of these Dirac operators in $L^{2}(S^{-}\otimes E)$ vary smoothly with $\xi $. 
\end{proof}
Therefore, we may set the following. 
\begin{defn} \label{deftrvb}
The \emph{transformed vector bundle} $\hat{E}$ of $(E,D,h)$ of a singular integrable 
connection with harmonic metric is the smooth vector bundle 
over $\hat{\C } \setminus \hat{P}$ whose fiber over $\xi $ is the finite-dimensional vector space 
$\hat{E}_{\xi } = coKer(\Dir_{\xi } ) \subset L^{2}(S^{-}\otimes E)$.
\end{defn}

In the remaining of this section, we prove vanishing of the kernel. 
The proof of the first statement of Theorem \ref{Fredholm} is left for the
next section. For the rest of the discussion in this section, we drop the index $\xi $.  
\begin{lem} \label{imorth} The subspaces $Im(\Dir  \vert_{H^{1}(\Omega^{0})})$ and 
$Im(\Dir  \vert_{H^{1}(\Omega^{2})})$ of $L^{2}(\Omega^{1})$ are orthogonal. 
\end{lem}
\begin{proof} 
Let $f_{0} \in H^{1}(\Omega^{0})$ and $f_{2}=g \d z \land \d\bar{z} \in H^{1}(\Omega^{2})$. Suppose first 
that $f_{0}$ is smooth and supported on a compact subset of $\C$, and such that near any 
singularity $p_{j } \in P$ its singular part is supported away from $p_{j}$. Then 
in a neighborhood of any $p_{j}$ in a holomorphic basis $Df_{0}=(\d +a)f$ for some bounded section 
$a \in \Om^{1}(End (\E))$, and so we have by partial integration 
\begin{equation}\label{ddstarint}
     \int_{\C  \setminus P} ( Df_{0},D^{\ast}f_{2} ) = 
     \int_{\C  \setminus P} ( D^{2}f_{0},f_{2} ) = 0,
\end{equation}
since $D$ is flat. 
Therefore, in order to finish the proof it is sufficient to show the following: 
\begin{clm} \label{dense} 
The set of compactly supported smooth sections of $S^{+} \otimes E$ on $\C$ 
with singular part compactly supported away from any singularity is dense in $H^{1}$.  
\end{clm}
\begin{proof}
It is sufficient to show the statement for $\Omega^{0}$, the case of $\Omega^{2 }$ being analogous. 
First we concentrate on infinity. 
Let $f \in H^{1}(E)$, and define cut-off functions $\rho_{R } (r) $ supported in $[0,2R]$ 
and equal to $1$ on $[0,R]$, such that $\rho_{R }'$ is supported in $[R,2R]$ with 
$$
    max| \rho_{R }' | \leq 2/R.
$$
Then we need to check that 
$$
    \rho_{R} (r)  f \lra f
$$
in $H^{1}(E)$ as $R \ra \infty $. In view of Corollary \ref{weightedcor}, this boils down 
to the classical calculations 
$$
     \| (1-\rho_{R} (r)) f \| \leq \int_{R\leq r} |f|^{2} 
$$
and
\begin{align*}
     \| \nabla^{+} ((1-\rho_{R} (r)) f) \| & \leq \int_{R \leq  r \leq 2 R} |\rho_{R }'(r)|^{2} |f |^{2} + K \int_{R \leq r} |\nabla^{+} f|^{2} \\
                                 & \leq K' \int_{R \leq  r \leq 2 R} |f |^{2} + K \int_{R \leq r} |\nabla^{+} f|^{2},
\end{align*}
where $K,K'$ are constants independent of $R$ and $f$.  

Next, let us consider a logarithmic singularity $p_{j}$, and define cut-off functions $\rho_{\varepsilon }$ 
supported in $[0,\varepsilon ] $, equal to $1$ on $[0, \varepsilon /2]$, and such that 
$$
      \max |\rho_{\varepsilon }' | \leq \frac{4}{{\varepsilon }}.
$$
We need to show that 
$$
    (1- \rho_{\varepsilon }) f^{sing} \lra f^{sing}
$$
in $H^{1}(E)$ as $\varepsilon \ra 0$. One sees that 
$$
     \int_{\C } |\rho_{\varepsilon }f^{sing}|^{2}  \leq \int_{r < \varepsilon }|f^{sing}|^{2}  \ra 0,
$$
since $f^{sing} \in L^{2}$. In the same way, 
$$
    \int_{\C } \frac{|\rho_{\varepsilon }f^{sing}|^{2}}{r^{2}} \leq \int_{r < \varepsilon } 
    \frac{|f^{sing}|^{2}}{r^{2}}\ra 0,
$$
since $f^{sing}/r \in L^{2}$. Finally, we also see that 
\begin{align*}
    \int_{\C } | \nabla^{+} (\rho_{\varepsilon } f^{sing})|^{2} & \leq \frac{16}{\varepsilon^{2} } 
    \int_{\varepsilon /2 < r < \varepsilon } |f^{sing}|^{2} + \int_{r < \varepsilon } | \nabla^{+} f^{sing}|^{2} \\
    & \leq \int_{\varepsilon /2 < r < \varepsilon } \frac{16 |f^{sing}|^{2}}{r^{2} } 
    + \int_{r < \varepsilon } | \nabla^{+} f^{sing}|^{2}
\end{align*}
and all of these expressions converge to zero as well. 
\end{proof}

Applying the claim to approximate $f_{0}$ and $f_{2}$ in $H^{1}$ by sections with compactly supported singular
component combined with (\ref{ddstarint}), we immediately get the lemma.
\end{proof}

Now we can come to vanishing of the kernel of (\ref{Dirac}): 
by Lemma \ref{imorth}, we have 
$$
     Ker(\Dir_{\xi}) = Ker(D_{\xi} \vert_{H^{1}(\Omega^{0})}) \oplus Ker(D^{\ast }_{\xi}  \vert_{H^{1}(\Omega^{2})}),
$$
it is therefore sufficient to prove vanishing of the kernels of $D$ and of
$D^{\ast }$. By duality, we only need to treat the case of $D$.
Harmonicity of the metric implies the Weitzenb\"ock formula: 
\begin{equation} \label{weitzen}
      \Dir^{\ast }_{\xi} \Dir_{\xi} = (\nabla^{+}_{\xi})^{\ast } (\nabla^{+}_{\xi}) + (\Phi_{\xi} \otimes )^{\ast } \Phi_{\xi} \otimes
\end{equation}
(see \cite{Biq97}, Thm 5.4.), which then gives by partial integration and Claim \ref{dense} the identity 
\begin{eqnarray} \label{trianglezero}
     \Vert \Dir_{\xi} f \Vert_{L^{2}}^{2} = \Vert D^{+}_{\xi}f \Vert_{L^{2}}^{2} + 
     \Vert \Phi_{\xi} f \Vert_{L^{2}}^{2}  
\end{eqnarray}
for any $f \in H^{1}(\Omega^{0})$. Suppose now that $f $ is in the kernel of $\Dir_{\xi}$. 
Then  (\ref{trianglezero}) implies $\Phi_{\xi} f=0$, and
since $\Phi_{\xi}$ is an isomorphism near infinity because of the choice $\xi\notin \hat{P}$, 
we also have there $f=0$. 
Again by (\ref{trianglezero}), $f$ is covariant constant. This gives
the result, since a covariant constant section vanishing on an open set
vanishes everywhere.

\section{Proof of the Fredholm Theorem }\label{sectpfFr}

A modification of the usual gluing argument of Fredholm-type theorems
works in this case as well. One lets $\phi_{1}$ be a cut-off
function supported in a compact region $R$ outside a neighborhood of the
singularities, and puts $\phi_{2}=1-\phi_{1}$. Since  $\Dir $ is a non-singular
first-order elliptic operator in $R$, elliptic theory of a compact
manifold implies that a parametrix $P_{1}$ exists for
$\Dir$ in this region. Next, one considers the problem in neighborhoods of the
singularities. First, one studies the model operators $\Dir^{j} = D^{j} +
(D^{j})^{\ast }$ instead of the Dirac operator itself. There are two different ways of treating these: 
\begin{enumerate}
\item{either one extends the functional spaces and the model Dirac operator onto
     a natural completion of the neighborhood, which can be either a conformal cylinder or
     a complex line (depending on the form of the metric and the functional
     spaces), and defines a two-sided inverse of $\Dir^{j}$ on this completion} \label{compl}
\item{or one finds directly a two-sided inverse of $\Dir^{j}$ on a small disk around the
  singularity, with a boundary condition verified by any section supported
  outside a neighborhood of the boundary.} \label{boundcond}
\end{enumerate}
Let us see how these allow to deduce the Fredholm theorem: 
if we take $R$ sufficiently
large, then on the support of $\phi_{2}$ all of these inverses $(\Dir^{j})^{-1}$ are defined. 
One then sets 
\begin{align*}
    & P : L^{2}(S^{-} \otimes E) \lra H^{1}(S^{+} \otimes E) \\
    & P(u)= \phi_{1}P_{1}(\phi_{1}u) + \sum_{j} \phi_{2} (\Dir^{j})^{-1}(\phi_{2}u),  
\end{align*}
and shows that this operator is a two-sided parametrix of $\Dir$ on all
$\C$. This can be done along classical lines, the only difference being
that near the singularities we have inverses of the local models of the
operator and not inverses of the operator itself. 
Therefore, we proceed as follows: first, we study the local models of the
Dirac operator around the singularities, and establish the isomorphisms as
in (\ref{compl}) or in (\ref{boundcond}). 
Then we prove that the effect of passing to the model operators
from the global ones at the singularities only amounts to
adding a compact operator $H^{1}(S^{+} \otimes E) \ra L^{2}(S^{-} \otimes E)$, which then gives the theorem. 

\subsection{Logarithmic singularities} Let $\Delta(p, \varepsilon )$ be a small 
neighborhood of $p \in P$. Up to a change of coordinates, we may suppose
$ \varepsilon =1$. Identify $ \Delta(p, 1 ) \setminus \{ p \} = S^{1} \times ]0,1]  $ via
polar coordinates $(r, \theta )$. Since the local model (\ref{polard}) is
diagonal in the basis $\{ e^{j}_{k} \}$, we see that the model Dirac operator on this disk  
$$ 
    \Dir^{j}_{0} = D^{j} - ( D^{j} )^{\ast }: (\Omega^{0} \oplus  \Omega^{2} ) \otimes
    E|_{ \Delta(p, 1 ) } \lra  \Omega^{1} \otimes E|_{ \Delta(p, 1)} 
$$ 
splits into the direct sum of its restrictions to the rank-one
components generated by one of the $\{ e^{j}_{k} \}$. Again, we have two
cases: first, $k \in \{ 1, \ldots r_{j} \}$ (regular case) and secondly $k \in \{
r_{j}+1, \ldots r \}$ (singular case). 

In the regular case, by definition the model Dirac operator on a rank-one
component is just the operator 
$$
   \Dir = \d - \d^{\ast} : S^{+} = \Omega^{0} \oplus  \Omega^{2} \lra  \Omega^{1} = S^{-},  
$$
which identifies to a projection of the real part of the usual Dirac operator 
on a product of two disks in $\C^{2}$ given by 
$$
    \dbar - \dbar^{\ast } : \Omega^{0,0} \oplus  \Omega^{0,2} \lra  \Omega^{0,1}. 
$$
Since this is known to have an inverse for the Atiyah-Patodi-Singer
boundary condition, the case of the regular part at a finite singularity
follows. 

On the singular component near a finite singularity, consider again the
coordinate change $t = - \ln r \in \R^{+} $. The local model of $D$ 
with respect to $t$ is given by 
\begin{align*}
     D^{j} & = d + i \bar{\mu}^{j}_{k} \d \theta + [\Re \mu^{j}_{k}-\beta^{j}_{k}] 
     \frac{\d r} {r} 
\end{align*}
(see (\ref{polard})).
Notice that the rank of $S^{+}$ and that of $S^{-}$ are both equal to $2$: 
we trivialize them using the unit-norm sections $(1,r \ \d r \land \d \theta )$ and 
$(\d r , r \d \theta)$ respectively, so that both $S^{+}\otimes E_{sing}$ and $S^{-}\otimes E_{sing}$ become isomorphic 
to $E_{sing}\oplus E_{sing}$ as Hermitian bundles. 
As we have seen in Lemma \ref{lemmaone}, the space $H^{1}( \Delta(p, 1),E_{sing})$ is equal to the model 
space of all sections $\phi$ having 
$$
    \int_{\Delta (p, 1 )} \left( | \nabla \phi |^{2} + \absl \frac{\phi }{r} \absr^{2} \right) \ r \d r
    \d \theta < \infty.
$$
By conformal invariance of the norm of $1$-forms and $\d t=\d r/r$, this is 
$$
    \int_{  S^{1} \times \R^{+} } \left( | \nabla \phi |^{2} + | \phi |^{2} \right)\ \d t \d \theta < \infty,  
$$
with the norm of the $1$-form $\nabla \phi $ measured with respect to the volume form $\d t \d \theta$. 
This latter is just the definition of the weighted Sobolev space $L^{2,1}_{0}(S^{1}\times \R^{+},E_{sing})$ 
with one derivative in $L^{2}$and weight $0$. In a similar way, the usual
$L^{2}$-space of sections of $E_{sing}$ on the disk is identified with the space 
$L^{2}_{-1}(S^{1} \times \R^{+},E_{sing})$ of $L^{2}$-sections with weight $-1$ on the half cylinder, for 
$$
    \int_{\Delta (p, 1 )}   | \phi |^{2} \ r \d r \d \theta =  \int_{S^{1} \times \R^{+}}|\phi e^{-t} |^{2} \d t \d \theta. 
$$
Hence in the trivialisation $(\d r , r \d \theta)$ of $S^{-}$, the usual $L^{2}$-space of $1$-forms 
on the disk is identified with the weighted space 
$$
     L^{2}_{-1}(S^{1} \times \R^{+}, E_{sing}\oplus E_{sing}).
$$
\begin{clm} 
Let $(r, \theta )$ be polar coordinates around $p=p^{j}$ as above. Let $k \in \{ r_{j}+1,\ldots,r \}$ and 
$$
    (f, g (r \d r \land \d \theta) ) \otimes e_{k}^{j} \in C^{\infty}(\Delta \setminus \{ 0 \},S^{+} \otimes E_{sing}).
$$
Then the value of the model Dirac operator $\Dir^{j}$ on this section is 
\begin{align*}
    & \left( \partial_{r}f + \frac{\Re \mu_{k}^{j}-\beta_{k}^{j}}{r} f - 
    \frac{\partial_{\theta }+i \mu_{k}^{j}}{r} g \right) dr \\ & + 
    \left(   \frac{\partial_{\theta }+i \bar{\mu}_{k}^{j}}{r} f + 
    \partial_{r} g - \frac{\Re \mu_{k}^{j} - \beta_{k}^{j}}{r} g \right) r \d \theta .
\end{align*}
In particular, in the unitary trivialisations $(1,r \ \d r \land \d \theta )$ and $(\d r , r \d \theta)$ of $S^{+}$ and
$S^{-}$, the operator 
$$
     r\Dir^{j}=e^{-t}\Dir^{j}
$$
is translation-invariant with respect to the cylindrical coordinate $t$. 
\end{clm}
\begin{proof} 
This is a direct computation: for $f \otimes e_{k}^{j}$ it follows
immediately from (\ref{polard}). For the image of 
$g (r \d r \land \d \theta) \otimes e_{k}^{j}$, consider first the smooth form $\varphi
\d r \otimes  e_{k}^{j} $ supported in a compact region of $\Delta \setminus \{
0 \}$; then by the same formula we have 
\begin{align*}
     \langle \varphi \d r \otimes  e_{k}^{j}, (D^{j})^{\ast } g (r \d r  \land \d \theta) \otimes e_{k}^{j} \rangle & = \langle D^{j} (\varphi \d
     r),g(r\d r \land \d \theta) \rangle  \\ 
     & = -\langle (\partial_{\theta }+ i \bar{\mu}_{k}^{j} ) \varphi \d r
     \land \d \theta , g (r \d r \land \d \theta) \rangle  \\ 
     & = - \frac{1}{r} \langle (\partial_{\theta }+ i
     \bar{\mu}_{k}^{j} ) \varphi , g \rangle \\ 
     & = \frac{1}{r} \langle \varphi , (\partial_{\theta } + i \mu_{k}^{j} )g  \rangle 
\end{align*}
and thus the projection of $(D^{j})^{\ast } g (r \d r \land \d \theta) \otimes e_{k}^{j} $ on the
$\d r$-component is $(\partial_{\theta } + i\mu_{k}^{j} )g \d r \otimes e_{k}^{j}$. 
The other component is obtained taking a compactly supported smooth 
form $\psi r \d \theta  \otimes e_{k}^{j}$: 
\begin{align*}
    \langle \psi r \d \theta \otimes  e_{k}^{j}, (D^{j})^{\ast } g ( r \d r  \land \d \theta )  \otimes  e_{k}^{j} \rangle & = 
    \langle D^{j} (\psi r \d \theta), g(r\d r \land \d \theta) \rangle \\
    & = \left\langle \left( \partial_{r} + \frac{\Re \mu_{k}^{j} - \beta_{k}^{j}}{r} \right) \psi , g \right\rangle \\
    & =\left\langle \psi ,  \left(-\partial_{r} + \frac{\Re \mu_{k}^{j} - \beta_{k}^{j}}{r}
      \right)g \right\rangle,  
\end{align*}
and the formula of the claim follows. It implies that $r\Dir^{j}$ is translation-invariant 
because $\partial_{r}=-\partial_{t}/r$.
\end{proof}

By definition, the weight $0$ is \emph{critical} for $r\Dir^{j}$ if and only if
there exists a non-trivial solution of
$$
    e^{-t}\Dir^{j} (A e^{ -\nu t+in \theta}, B  e^{-\nu t + in \theta } (r \ \d r \land  \d \theta)) = 0
$$
with some constants $A,B \in \C $ and a constant $\nu \in \C$ such that $\Re \nu=0$. 
Turning back to the coordinate $r$ again, this is equivalent to having
\begin{eqnarray} \label{critval}
    r\Dir^{j} (A r^{\nu } e^{in \theta }, B r^{\nu } e^{in \theta } (r \ \d r \land  \d \theta)) = 0.
\end{eqnarray}
By \cite{Lock-McO}, if $0$ is not a critical weight, then the translation-invariant elliptic 
differential operator 
$$
     e^{-t}\Dir^{j} : L^{2,1}_{0}(S^{1} \times \R^{+}, S^{+} ) \lra L^{2}_{0} (S^{1}\times \R^{+},S^{-}) 
$$
is invertible, and thus so is 
$$
     \Dir^{j} : L^{2,1}_{0}(  S^{1} \times \R^{+}, S^{+} ) \lra L^{2}_{-1} (S^{1}\times \R^{+},S^{-}) 
$$
since 
$$
     e^{t}:L^{2}_{0} \lra L^{2}_{-1}
$$ 
is an isomorphism. Therefore, in order to establish the desired isomorphism in the singular
case, we only need to check the weight $0$ is not critical for $r\Dir^{j}$.

Applying the claim to the equation (\ref{critval}), we see that $0$ is a
critical weight if and only if the system of linear equations 
\begin{align*}
        (\nu + \Re \mu -\beta )A - i (n + \mu )B & = 0 \\
        i (n + \bar{\mu} ) A + (\nu + \beta - \Re \mu )B &= 0
\end{align*}
has a non-trivial solution $(A,B)\in \C^{2}$ for some $\nu \in \C$ with $\Re \nu=0$ 
(here we have omitted indices $j$ and $k$ of $\mu$ and $\beta$ for simplicity). 
This system has a non-trivial solution  if and only if the determinant formed by the
coefficients is equal to $0$: 
$$
     \nu ^{2} - (\Re \mu - \beta )^{2} - \vert n + \mu \vert ^{2} =0.
$$
Since $\Re \nu $ must be $0$, this can only be the case if $\nu = \Re \mu - \beta  = n+ \mu
=0$. By assumption $0 \leq \beta < 1$, and $n$ is an integer, therefore the
only case this can hold is when $\beta = \mu = 0$, which is impossible, since we
are looking at the singular component of the bundle. Therefore, there are
no non-trivial solutions to (\ref{critval}), and  $0$ is not a critical
weight.  

\subsection{Singularity at infinity} 
In this section the importance of the condition $\xi \notin \hat{P}$ will come out; 
therefore we write out the index $\xi $ of our operators.
A neighborhood of infinity in $\C \setminus P$ is given by the complementary 
$\C \setminus \Delta(R) $ of a large disk around $0$. A natural choice of
completion of this manifold is of course $\C $, with its standard metric
$|\d z|^{2}$. We choose to study the local model in the orthonormal basis $\{ e^{\infty }_{k}
\}$  defined in (\ref{ebasisinf}). This
allows us to think of $E$ as the trivial bundle $\C^{r}$ over $\C \setminus \Delta(R)$, with
standard hermitian metric on the fibers. By (\ref{holtrivinf}) this basis (up to a
polynomial scaling factor) is a natural one for the Higgs-bundle point of
view, so the deformation is that considered in (\ref{higgsdef}), and 
the operator $D_{\xi }$ near infinity is given (up to terms of order $r^{-1}$) by 
$$ 
   D_{\xi }^{\infty } = \d + \frac{A- \xi \Id }{2} \d z + \frac{(A- \xi \Id )^{\ast }}{2}
   \d \bar{z} 
$$ 
(see (\ref{singpart})), and a natural extension of it to all of $\C $ can be
given by the same formula. This implies immediately that  
$$
    \Phi_{\xi }^{\infty } =  \frac{A- \xi \Id }{2} \d z + \frac{A^{\ast }- \bar{\xi} \Id }{2} \d \bar{z}
$$ 
and $(D^{\infty })^{+} = \nabla $ (the trivial connection) on all of $\C$. 
For a section $\phi \in L^{2}(\Omega^{0})$ supported in $\C \setminus \Delta(R)$, the condition 
$\Phi_{\xi } \phi \in L^{2}( \Omega^{0} )$ then automatically holds, and $(D_{\xi }^{\infty })^{+} \phi \in L^{2}$ is 
equivalent to $\nabla \phi \in L^{2}$. Therefore, on sections of $\Omega^{0}$ supported on the
complementary of $\Delta(R)$, the $H^{1}$-norm is equivalent to the usual
Sobolev $L^{2,1}$-norm. A similar argument shows that for sections of $\Omega^{2}$, the 
$H^{1}$-norm is also equivalent to the usual $L^{2,1}$-norm. 
Therefore, on all of $\C $, we must consider a natural extension of these 
functional spaces, namely $L^{2,1} (\C , \Omega^{0} \oplus \Omega^{2})$. In an analogous manner, on $S^{-}$ we
consider the extension $L^{2}(\C, \Omega^{1})$ of $L^{2}( \C \setminus \Delta(R) ,
\Omega^{1})$. Therefore, we need to prove the 
\begin{lem} On  $\C $, the Dirac operator 
\begin{eqnarray} \label{diracinf} 
    \Dir_{\xi }^{\infty }= D_{\xi }^{\infty } -  (D_{\xi }^{\infty })^{\ast } : L^{2,1}( \Omega^{0} \oplus \Omega^{2} ) \lra
    L^{2} ( \Omega^{1} ) 
\end{eqnarray}
is an isomorphism. 
\end{lem} 
\begin{proof}
Since $A$ is supposed to be diagonal in this basis with eigenvalues
$\xi_{l}  \ (l=1,\ldots , n')$, we may restrict ourselves to the study of the
operator $D^{\infty } = \d + (\xi_{l}-\xi )/2 \d z + (\bar{\xi }_{l}-\bar{\xi })/2 \d \bar{z}$. We need the 
following:  
\begin{clm} Denote by $\Delta$ the plain Laplace operator $\nabla^{\ast } \nabla$ on forms. Then we have
\begin{eqnarray} \label{weitz}
    \Dir_{\xi }^{\infty} (\Dir_{\xi }^{\infty})^{\ast }  = - \Delta - \frac{\vert \xi_{l }-\xi  \vert^{2}}{4}. 
\end{eqnarray}
\end{clm}
\begin{proof} 
This is an easy computation. 
\end{proof}
Now recall that by the classical theory of the
Laplace operator, $ \Delta +  \lambda^{2} $ with $ \lambda > 0$ is an
isomorphism 
\begin{equation} \label{isominf}
    L^{2,2}(\C , \Omega^{j} ) \lra L^{2} (\C , \Omega^{j} ). 
\end{equation}
This statement can be for example obtained passing to the Fourier transform $|\hat{x}|^{2}+ \lambda^{2}$ of this 
operator.

Coming back to our situation, the condition $\xi \notin \hat{P}$ means exactly that $\xi_{l}-\xi \neq 0$ for any $l=1,\ldots,n'$. 
This immediately implies that (\ref{diracinf}) is surjective: indeed, 
clearly $Im( (\Dir_{\xi }^{\infty})^{\ast } ) \subset L^{2,1}( \Omega^{0} \oplus \Omega^{2} )$, and 
$ \Dir_{\xi }^{\infty} (\Dir_{\xi }^{\infty})^{\ast } $ is surjective by the isomorphism 
(\ref{isominf}). For injectivity, note that a formula similar to
(\ref{weitz}) holds for the Laplace operator $ (\Dir_{\xi }^{\infty})^{\ast } \Dir_{\xi }^{\infty} $ as well. 
This in turn implies that the $L^{2,2}$-kernel of $\Dir_{\xi }^{\infty}$ vanishes. Elliptic regularity then shows
that the $L^{2,1}$-kernel vanishes as well. 
\end{proof}

\subsection{Compact perturbation}
We wish to prove that near each one of the singularities the effect of 
passing from the global operator to its local model, i.e. subtracting 
the perturbation term only amounts to a
compact operator $H^{1}(S^{+}\otimes E ) \ra L^{2}(S^{-}\otimes E )$. This then
finishes the proof of the Fredholm theorem, because the sum of a Fredholm operator 
and a finite number of compact operators is Fredholm. 

Consider first the case of a singularity at a finite point. Recall from 
Lemma \ref{lemmaone} that near $p_{j}$ the space $H^{1}(S^{+} \otimes E)$ 
is equal to the sum 
$$
      L^{2,1}_{eucl} (S^{+} \otimes E_{reg}) \oplus L^{2,1}_{0} (S^{+} \otimes E_{sing}),
$$
where $L^{2,1}_{eucl}$ is the usual Sobolev space on the disk of $L^{2}$-functions 
with one derivative in $L^{2}$ with respect to Euclidean metric, whereas $L^{2,1}_{0}$
is the weighted Sobolev space defined by 
$$
    \int_{\Delta(p_{j},\varepsilon )} \left( \absl \frac{\phi }{r} \absr^{2} + |\nabla \phi
      |^{2} \right) |\d z |^{2} \leq \infty . 
$$
Also, the order of growth of the $1$-form perturbation term $a^{j}$ with
respect to Euclidean metric is by (\ref{pertj}) at most $O(r^{-1+\delta })$,
with $\delta >0$. We need to prove that we have compact Sobolev multiplications for 
functions on the disk 
\begin{equation} \label{mult1}
     L^{2,1}_{eucl} \xrightarrow{a^{j}} L^{2}_{eucl}
\end{equation}
and 
\begin{equation} \label{mult2}
     L^{2,1}_{0} \xrightarrow{a^{j}} L^{2}_{eucl}.
\end{equation}
Consider first (\ref{mult1}): since the disk is a compact manifold, for any
$2<p<\infty $ the inclusion $L^{2,1}_{eucl} \inc L^{p}_{eucl}$ is compact. On the other hand,
$O(r^{-1+\delta }) \d r +  O(r^{-1+\delta }) r \d \theta $ is in $L^{2+\varepsilon }_{eucl}$ 
for some $\varepsilon>0$. Choose $p$ such that $1/2=1/(2+\varepsilon) + 1/p$; 
(\ref{mult1}) then follows immediately from the 
continuous multiplication $L^{2+\varepsilon}_{eucl} \times L^{p}_{eucl} \ra  L^{2}_{eucl} $.
Now, we come to (\ref{mult2}): this is an immediate consequence of the
previous, for the weighted norm $L^{2,1}_{0}$ is stronger then
$L^{2,1}_{eucl}$. 

Next, let us treat the case of the singularity at infinity. In the
coordinate $w=1/z$ we have a second-order singularity on the disk 
$\Delta(0,1/R)$. Let $w= \rho e^{i\vartheta }$; by (\ref{pertinf}) the perturbation is $O(\rho^{-1-\delta})$, and 
the $H^{1}$-norm of a function $\phi $ supported near infinity is given by 
$$
    \int_{ \C \setminus \Delta (0,R)} \left( |\phi |^{2} + |\nabla \phi |^{2}  \right) | \d z |^{2} = 
    \int_{ \Delta (0,1/R)} \left( \absl \frac{\phi }{\rho^{2} } \absr^{2} + |\nabla \phi |^{2} \right)  |\d w|^{2}.
$$ 
In particular, in the coordinate $w$ this norm is also stronger then $L^{2,1}_{eucl}$, so we
conclude from (\ref{mult1}). 

\section{$L^{2}$-cohomology and Hodge theory}\label{l2coh}
In this section we keep on supposing that we have on one side an integrable
connection $D$ with singularities in $P \cup \{ \infty \}$, with prescribed behaviors
at these points, given in regular singularities by (\ref{pertj}) and at infinity by (\ref{pertinf}). 
In Theorem \ref{Fredholm} we proved that the deformed operators $\Dir_{\xi}$
are Fredholm between the spaces $H^{1}$ and $L^{2}$; in particular their indices agree. 
We also showed that if the metric is harmonic then the kernel of the Dirac operator vanishes, hence the 
index of $\Dir_{\xi}$ is equal to the opposite of the dimension of the cokernel
$Coker(\Dir_{\xi } )$, this operator being considered between functional
spaces as in (\ref{Dirac}). This dimension is therefore a constant independent of $\xi $, and it follows 
from the implicit function theorem that the spaces $\hat{E}_{\xi }=Coker(\Dir_{\xi } )$ define a finite-rank smooth
vector bundle $\hat{E}$ over $\hat{\C} \setminus \hat{P} $, the rank being
equal to the opposite of the index of (\ref{Dirac}). 
Here we wish to interpret this cokernel as the first cohomology of the elliptic complex 
\begin{eqnarray} \label{ltwocompl}
     L^{2}(\Omega^{0} \otimes E) \xrightarrow{\text{$D_{\xi}$}}  L^{2}(\Omega^{1} \otimes E)   
                            \xrightarrow{\text{$D_{\xi}$}} L^{2}(\Omega^{2} \otimes E), 
\end{eqnarray}
(see Theorem \ref{coker}), and also as the space of harmonic sections with respect to the Laplace operator
of the adjoint Dirac operator $\Dir_{\xi }^{\ast }$ (Theorem \ref{laplker}). 

Since the operators in (\ref{ltwocompl}) are unbounded, we need to define their domains. 
In this chapter $C^{\infty }_{0}$ stands for smooth sections supported in a compact subset of $\C \setminus P$.
\begin{defn} The \emph{maximal domain} of $D|_{\Omega^{i}}$ is  
$$ 
     \Dm(D|_{\Omega^{i}}) = \{ u \in L^{2}(\Omega^{i}): Du \in L^{2}(\Omega^{i+1}) \} ,
$$  
where $Du \in L^{2} $ is understood in the sense of currents, i.e. 
the functional $v \in C^{\infty }_{0}(\Omega^{i+1}) \mapsto \langle u, D^{\ast }v \rangle $ is
continuous in the $L^{2}$-topology. 
\end{defn} 
By local elliptic regularity, this amounts to the same thing as $Du$ being
an $L^{2}$-section. When it does not cause any confusion,
we will simply write  $\Dm (\Omega^{i})$ for $\Dm(D|_{\Omega^{i}})$. 
It is easy to see that if we consider $D$
on its maximal domain, then the kernel is a closed subspace of $L^{2}$, and
the image of $D$ on $\Omega^{i-1} $ is contained in the kernel of $D$ on $\Omega^{i}$. 
The image of a general differential operator is however not always a closed subspace of the kernel. 
\begin{defn} For $i \in \{ 0,1,2 \}$, the $i^{th}$ \emph{$L^{2}$-cohomology} of $D$  is
   $Ker(D\vert _{\Omega^{i} \otimes E})/ {Im( D\vert _{\Omega^{i-1} \otimes E})} $, where both
   of these operators are considered with maximal domain, and the
   operators not shown in (\ref{ltwocompl}) are trivial.  It is denoted by $L^{2}H^{1}(D)$.  
\end{defn}
Our aim is to obtain the following: 
\begin{thm} \label{coker} The cokernel of $\Dir $ defined on $H^{1}(S^{+} \otimes E)$ is equal to the first  
  $L^{2}$-cohomology of $D$. 
\end{thm}
\begin{proof}
Recall that by definition 
\begin{align}
    Coker(\Dir \vert_{H^{1}(S^{+} \otimes E)} ) & = (Im(\Dir \vert_{H^{1}(S^{+} \otimes E)} ))^{\bot }\notag \\
    & = (Im(D \vert_{H^{1}(\Omega^{0} \otimes E)} ))^{\bot } \cap  (Im(D^{\ast } \vert_{H^{1}(\Omega^{2} \otimes E)} ))^{\bot },\label{descrcoker} 
\end{align}
where $A^{\bot }$ stands for the $L^{2}$-orthogonal of the subspace $A \subset
L^{2}$. Therefore, it is sufficient to prove the following lemmas: 
\begin{lem} \label{maxdom} 
The maximal domain of 
$$ 
    D:L^{2}(\Omega ^{0} \otimes E) \lra L^{2}(\Omega ^{1} \otimes E)
$$ 
is $H^{1}(\Omega ^{0} \otimes E)$. Similarly, 
the maximal domain of 
$$ 
    D^{\ast }:L^{2}(\Omega ^{2} \otimes E) \lra L^{2}(\Omega ^{1} \otimes E)
$$
is $H^{1}(\Omega ^{2} \otimes E)$. In particular, the maximal domain of 
$$
    \Dir :L^{2}(S^{+} \otimes E) \lra L^{2}(S^{-} \otimes E)
$$ 
is $H^{1}(S^{+} \otimes E)$. Moreover, if this latter space is equipped with the norm $\| . \|_{H^{1}}$ defined in (\ref{norm}), 
then $\Dir$ is a bounded operator from $H^{1}(S^{+} \otimes E)$ to $L^{2}(S ^{-} \otimes E)$. 
\end{lem}
\begin{lem} \label{imker} We have 
$$
    (Im(D^{\ast } \vert_{H^{1}(\Omega^{2} \otimes E)} ))^{\bot } = Ker( D \vert_{\Dm(\Omega ^{1}
    \otimes E)}).
$$
\end{lem} 
\begin{lem} \label{closed} The image of $D: H^{1}(\Omega ^{0} \otimes E) \ra  L^{2}(\Omega ^{1}  \otimes E)$ is
  closed. 
\end{lem} 
Indeed, Lemmas \ref{maxdom} and \ref{imker} together with (\ref{descrcoker}) imply
that the cokernel is equal to 
$$
   (Im(D \vert_{\Dm(\Omega^{0} \otimes E)} ))^{\bot } \cap  Ker( D \vert_{\Dm(\Omega ^{1}\otimes E)}),  
$$
which in turn is identified to the first reduced $L^{2}$-cohomology of
(\ref{ltwocompl}), i.e. to 
$$
    Ker(D\vert _{\Dm(\Omega^{1} \otimes E)})/ \overline{Im( D\vert _{\Dm(\Omega^{0} \otimes E)})}, 
$$
where the bar over the image stands for the $L^{2}$-closure of that
space. Lemma \ref{closed} now concludes the proof of Theorem \ref{coker}. 

\begin{proof} (Lemma \ref{imker})  We first show the
\begin{clm} The adjoint of the unbounded operator 
\begin{eqnarray} \label{Dstar}
    D^{\ast } : L^{2} (\Omega^{2} \otimes E) \lra L^{2}(\Omega^{1} \otimes E) 
\end{eqnarray}
with domain $H^{1}(\Omega^{2} \otimes E)$ is the unbounded operator 
\begin{eqnarray} \label{Dstaradj}
    D: L^{2}(\Omega^{1} \otimes E) \lra L^{2} (\Omega^{2} \otimes E)
\end{eqnarray}
with domain $\Dm(\Omega^{1} \otimes E)$. 
\end{clm}
\begin{proof} (Claim)
It is clear that the formal adjoint of (\ref{Dstar})
is (\ref{Dstaradj}), we only need to prove its domain is $\Dm$. By
definition, a section $u \in L^{2}( \Omega ^{1})$ 
is in the domain of the adjoint operator $\mbox{Dom}((D^{\ast })^{\ast })$ if and
only if for all $v \in H^{1} (\Omega^{2} \otimes E)$ we have 
$$
    \vert \langle u,D^{\ast }v \rangle \vert \leq K \| v \| 
$$
with a constant $K$ only depending on $u$. Now, since $v \in H^{1}$ and $u \in
L^{2}$, by Claim \ref{dense} we can perform partial integration to the left-hand side
of this formula. Therefore, $u$ is
in the domain of the adjoint operator if and only if the functional 
$$
     v \mapsto  \langle Du, v \rangle 
$$
is bounded in $L^{2}(\Omega^{2} \otimes E)$. But this condition is equivalent to $Du \in
L^{2}(\Omega^{2} \otimes E)$, and the claim follows. 
\end{proof}

Lemma \ref{imker} now directly follows  from the claim and the general fact
that the cokernel of an unbounded operator is equal to the kernel of its
adjoint. 
\end{proof}

\begin{proof} (Lemma \ref{maxdom}) First we need to prove that for
a section $u$ of $ L^{2}(\Omega ^{0} \otimes E) $ we have $Du \in L^{2}$ if and only if
both $D^{+}u \in L^{2}$ and $\Phi u \in L^{2}$. The ``if '' direction being obvious, we
concentrate ourselves on the opposite statement, and suppose in what
follows that $u$ is an $ L^{2}$-function with $Du \in L^{2}$. 

We first study the singularity at infinity. For $|u|$ sufficiently large, we have the point-wise 
estimate 
$$
    \vert \Phi u \vert \leq 2 K \vert u \vert,  
$$ 
where $K$ is the maximal modulus of the eigenvalues of the matrix
$A$. Therefore, $u \in L^{2}$ at infinity implies $\Phi u \in L^{2}$ at infinity, 
and consequently $D^{+}u=Du - \Phi u \in L^{2}$ at infinity, and we are done. 

Next, consider the case of a singularity at a finite point. 
In the orthonormal basis (\ref{ebasisj}), the operators we study are equal, up to a perturbation
term, to the local models (see (\ref{polarunit}), (\ref{polarphi}), (\ref{polard}))  
\begin{align*}
             (D^{+})^{j} \phi  & =(\mbox{d} + i \Re \mu^{j}_{k} \mbox{d} \theta) \phi  \\ 
             \Phi^{j} \phi  & = [(\Re \mu^{j}_{k} - \beta^{j}_{k})  \frac{\d r}{r} + \Im
             \mu^{j}_{k} \d \theta ] \phi \\
             D^{j}\phi  & = [\mbox{d} + i \bar{\mu}^{j}_{k} \d \theta  + (\Re \mu^{j}_{k} -
             \beta^{j}_{k})  \frac{\d r}{r} ]\phi 
\end{align*}
To simplify  notation, from now on we drop the indices $j$ and $k$. Note
that because of Lemma \ref{lemmaone}, it is sufficient to prove that $\Phi^{j}
\phi $ and $(D^{+})^{j} \phi$ are in $L^{2}$. 
Notice also that since the perturbation $a^{j}$ may mix the regular and singular
components, a priori it is not sufficient to prove for example that
$ \phi_{reg } \in L^{2}$ and $D \phi_{reg } \in L^{2} $ imply $(D^{+})^{j} \phi_{reg } \in
L^{2}$, because $D \phi \in L^{2}$ does not imply directly $D \phi_{reg } \in L^{2}$
in the presence of a mixing perturbation term. However, remark that
denoting by $a^{j}_{r,r}$ the part of the endomorphism $a^{j}$ that takes
the regular component into the regular one, and $a^{j}_{r,s},a^{j}_{s,r},a^{j}_{s,s}$
the other parts, we have 
\begin{align} \label{mix}
        \int_{\Delta(p_{j}, \varepsilon ) } |(D^{j}+a^{j})\phi |^{2} = & \int_{\Delta (p_{j}, \varepsilon )} |(D^{j}+a^{j}_{r,r})
        \phi_{reg}+a^{j}_{s,r}\phi_{sing}|^{2}\notag \\
        & +  \int_{\Delta(p_{j}, \varepsilon ) } |(D^{j}+a^{j}_{s,s})
        \phi_{sing}+a^{j}_{r,s}\phi_{reg}|^{2}  \\ 
        \geq & \int_{\Delta(p_{j}, \varepsilon ) } |D^{j} \phi_{reg}|^{2} + |D^{j} \phi_{sing}|^{2}\notag  \\
        & - |a^{j}\phi_{reg}|^{2} - |a^{j}\phi_{sing}|^{2}, \notag 
\end{align}
and this estimate shows that we can treat the two components separately: the left-hand side 
is finite by hypothesis, whereas the integrals of $|a^{j}\phi_{reg}|^{2}$ and $|a^{j}\phi_{sing}|^{2}$ 
by Kato's inequality and (\ref{regineq}); hence the same thing holds for the integrals of 
$|D^{j} \phi_{reg}|^{2}$ and $|D^{j} \phi_{reg}|^{2}$.

On the regular component, the above expressions simplify to $D^{j}=(D^{+})^{j}=\nabla $ (the
trivial connection), and $\Phi^{j} =0$. What we need to show is that $\phi_{reg} ,D \phi_{reg}  \in
L^{2}$ implies $\nabla \phi_{reg}  \in L^{2}$, if $D =\nabla + a^{j}$ with $a^{j} =O(r^{-1+\delta })$. 
Recall that by Kato's inequality and (\ref{regineq}) with $\varepsilon >0$ chosen sufficiently small
we have
$$
     \int_{\Delta(p_{j}, \varepsilon )}  \vert a^{j} \phi_{reg} \vert^{2} \leq   
     \int_{\Delta(p_{j}, \varepsilon )} \vert D \phi_{reg} \vert^{2} + 
     \int_{\Delta(p_{j}, \varepsilon )\setminus \Delta(p_{j}, \varepsilon /2 )} \vert \phi_{reg} \vert^{2}. 
$$
It follows that 
\begin{align*}
    \int_{\Delta(p_{j}, \varepsilon )} \vert \nabla \phi_{reg} \vert^{2} 
    & \leq  \int_{\Delta(p_{j}, \varepsilon )} \vert D \phi_{reg} \vert^{2} + 
    \int_{\Delta(p_{j}, \varepsilon )}  \vert a^{j} \phi_{reg} \vert^{2} \\
    & < 2 \int_{\Delta(p_{j}, \varepsilon )} \vert D \phi_{reg} \vert^{2} + 
    2\int_{\Delta(p_{j}, \varepsilon )} \vert  \phi_{reg} \vert^{2}. 
\end{align*}
Now by the hypothesis $\phi, D \phi \in L^{2}$, the right-hand side
is finite. Therefore $\nabla \phi  \in L^{2}$ as we wished to show. 

Consider now the singular case: again, we need to show that if we have a 
section $\phi \in L^{2}$ such that $D\phi \in L^{2}$, then $D^{+}\phi_{sing}, \Phi \phi_{sing} \in L^{2}$.
Here, usual elliptic regularity does not give the claim,
because we need to deduce that $\phi_{sing} /r \in L^{2}$. From now on, we write
$\phi = \phi_{sing}$ to lighten notation. Decompose $\phi$ into its Fourier-series near $p^{j}$: 
$$
    \phi(r,\theta )= \sum_{n=-\infty }^{\infty } \phi_{n}(r)e^{in\theta }
$$
Choosing $\varepsilon$ sufficiently small, we can make the perturbation
term $a^{j}$ be smaller on $\Delta(p_{j}, \varepsilon )$ then $\nu/r$ for
any $\nu>0$.   
Write first the $\d \theta $-term of $D^{j}\phi$: 
$$
    D_{\theta }^{j}\phi = (\partial_{\theta }+ i \bar{\mu })\phi \d \theta = i \d \theta  \sum_{n=-\infty }^{\infty } (n + \bar{\mu })
    \phi_{n}(r) e^{in\theta }. 
$$ 
By this and the estimate on the perturbation, we infer that 
\begin{align} \label{mixsing}
   \Vert (D^{j}_{\theta } +a^{j}) \phi \Vert_{L^{2}(\Delta(p_{j}, \varepsilon ))}^{2} \geq 
   \Vert D^{j}_{\theta } \phi \Vert_{L^{2}(\Delta(p_{j}, \varepsilon ))}^{2} -  
    \Vert \nu  \phi /r \Vert_{L^{2}(\Delta(p_{j}, \varepsilon ))}^{2}& \notag \\  
    =  \int_{\Delta(p,\varepsilon )} \sum_{n=-\infty }^{\infty } (\vert n + \bar{\mu } \vert^{2} 
      -\nu^{2})\frac{\vert \phi_{n}(r) \vert^{2}}{r^{2}}& \\
      = \int_{\Delta(p,\varepsilon )} \sum_{n=-\infty }^{\infty } (\vert n + \Re \mu  \vert^{2} -
       \nu^{2} + \vert \Im \mu \vert^{2}) \frac{\vert \phi_{n}(r) \vert^{2}}{r^{2}}&. \notag
\end{align}
By Hypothesis \ref{main} we have $\Re \mu \notin \Z$, and so if $\nu $ is sufficiently small, then
the last expression can be bounded from below by 
\begin{align} \label{est}
   \frac{1}{2} \int_{\Delta(p,\varepsilon )}& \sum_{n=-\infty }^{\infty } (\vert n + \Re \mu  \vert^{2} +
    \vert \Im \mu \vert^{2}) \frac{\vert \phi_{n}(r) \vert^{2}}{r^{2}} \\
    = &\frac{1}{2} \int_{\Delta(p,\varepsilon )} |(D_{\theta }^{+})^{j}\phi|^{2}+|\Phi_{\theta }^{j} \phi|^{2}.\notag
\end{align}
As in the regular case, by (\ref{regineq}) the left-hand side of (\ref{mixsing}) is finite, 
so we see that $(D_{\theta }^{+})^{j}\phi \in L^{2}$ and $\Phi_{\theta }^{j} \phi \in L^{2} $.
The $\d r$-part $\Phi_{r}^{j} \phi $ of $\Phi^{j} \phi $ is in $L^{2} $ if and only if
$$
      \int_{\Delta(p,\varepsilon )} \vert \Re \mu - \beta \vert^{2} \frac{\vert \phi (r)
       \vert^{2}}{r^{2}}   < \infty .
$$
Again by our main hypothesis $\Re \mu \notin \Z$ there exists a constant
$K >0$ such that  
$$
     \sum_{n=-\infty }^{\infty } \vert \Re \mu - \beta \vert^{2} \frac{ \vert \phi_{n} (r)
       \vert^{2}}{r^{2}} \leq  
     K \sum_{n=-\infty }^{\infty } \vert n+ \Re \mu \vert^{2} \frac{\vert \phi_{n} (r) \vert^{2}}{r^{2}}.
$$
As we have already seen, this last expression is integrable, therefore
$\Phi^{j} \phi \in L^{2}$. Since the perturbation is negligible compared to the
behavior $O(r^{-1})$ of (\ref{est}), we then also have $\Phi \phi \in L^{2}$. 
We conclude  using $D^{+}\phi = D\phi - \Phi\phi$.  

By duality, the case of a $2$-form $v \d z \land \d \bar{z}$ is settled the same way. 
The general case (that of $S^{+} \otimes E$) then follows from Lemma \ref{imorth}. 
The fact that 
$$
     \Dir : H^{1}(S^{+} \otimes E) \lra L^{2} (S^{-} \otimes E) 
$$
is bounded, is then immediate (and has already been pointed out, see (\ref{Dirac})).
\end{proof}

\begin{proof} (Lemma \ref{closed}) This is immediate 
from Theorem 1 and Claim \ref{imorth}.  \end{proof}
We have established lemmata \ref{imker}, \ref{closed} and \ref{maxdom}, hence we 
finished the proof of Theorem \ref{coker}. 
\end{proof}

\begin{thm} \label{laplker} The first $L^{2}$-cohomology of the complex (\ref{ltwocompl})
is canonically isomorphic to the kernel of the adjoint Dirac operator 
\begin{equation} \label{diradj}
     \Dir_{\xi}^{\ast }: L^{2}(S^{-} \otimes E) \lra L^{2}(S^{+} \otimes E)
\end{equation}
on its domain, or alternatively to the kernel of the Laplace operator 
\begin{equation} \label{lapldirac}
     \Delta_{\xi } = \Dir_{\xi } \Dir_{\xi}^{\ast } = - D_{\xi } D_{\xi }^{\ast } - D_{\xi }^{\ast } D_{\xi } :
     L^{2}(S^{-} \otimes E) \lra L^{2}(S^{-} \otimes E)
\end{equation}
on its domain.  
\end{thm}
\begin{proof}
By duality, we get from Lemma \ref{imker} that 
$$
      (Im(D|_{H^{1}(\Omega^{0} \otimes E)}))^{\bot } = ker(D^{\ast }|_{\Dm(\Omega^{1} \otimes E)}), 
$$
and this implies 
\begin{align*}
    coKer(\Dir|H^{1}(S^{+} \otimes E )) & = ker( D^{\ast }|_{\Dm(\Omega^{1} \otimes E)} ) \cap 
    ker( D|_{\Dm(\Omega^{1} \otimes E)} ) \\
    & = ker( \Dir^{\ast }|_{\Dm(\Omega^{1} \otimes E)} ) .
\end{align*}
It remains to show that this latter is equal to $ker( \Dir \Dir^{\ast }|_{\Dm(\Omega^{1} \otimes E)} )$. 
It is clear that 
$$
     ker( \Dir \Dir^{\ast }|_{\Dm(\Omega^{1} \otimes E)} ) \supseteq  ker( \Dir^{\ast }|_{\Dm(\Omega^{1} \otimes E)}).
$$ 
Suppose now $u \in L^{2}(\Omega^{1} \otimes E)$ satisfies $\Dir \Dir^{\ast }u=0$. This means that 
$$
   \Dir^{\ast }u\in Ker(\Dir)\subset \Dm(\Dir)=H^{1}(S^{+} \otimes E)
$$ 
by Lemma \ref{maxdom}. Vanishing of the $L^{2}$-kernel of $\Dir$ on $H^{1}(S^{+} \otimes E)$ 
(c.f. Theorem \ref{Fredholm}) gives $\Dir^{\ast }u=0$, that is $u \in Ker(\Dir^{\ast })$, whence 
$$
    ker( \Dir \Dir^{\ast }|_{\Dm(\Omega^{1} \otimes E)} ) \subseteq ker( \Dir^{\ast }|_{\Dm(\Omega^{1} \otimes E)}).
$$
\end{proof}

Finally, let us introduce the norm 
$$
     \| f \|_{H^{2} (S^{+} \otimes E) } = \int_{\C } |f|^{2} + |(\nabla^{+})^{\ast } \nabla^{+} f|^{2} + 
      |(\Phi \otimes )^{\ast } \Phi \otimes f|^{2}
$$
and the corresponding function space 
$$
    H^{2} (S^{+} \otimes E) = \{ f : \hspace{5mm}  \| f \|_{H^{2} (S^{+} \otimes E) } < \infty \} 
$$
Then we have the following. 
\begin{thm} \label{Hilbertisom}
The domain of the Laplace operator $\Delta_{\xi }=\Dir_{\xi}^{\ast } \Dir_{\xi }$ is $H^{2} (S^{+} \otimes E)$. It defines 
a Hilbert-space isomorphism 
$$
     H^{2} (S^{+} \otimes E) \lra  L^{2} (S^{+} \otimes E).
$$
\end{thm}
\begin{proof}
The fact that $\Delta_{\xi }$ is a well-defined bounded operator on $H^{2} (S^{+} \otimes E)$ follows from the 
Weitzenb\"ock formula (\ref{weitzen}). Its maximal domain is the set of $u\in L^{2} (S^{+} \otimes E)$ 
such that $\Dir_{\xi }u\in \Dm(\Dir_{\xi}^{\ast })$. This latter is, by computations similar to Lemma \ref{maxdom}, 
the Sobolev space $H^{1}(S^{-} \otimes E)$ is with $1$ derivative in $L^{2}$, and weight $-1$ on the irregular 
component near logarithmic singularities like in Corollary \ref{weightedcor}. 
We deduce that the maximal domain of $\Delta_{\xi }$ is $H^{2}(S^{+} \otimes E)$, and that it splits as 
$$
       H^{2}(S^{+} \otimes E) \xrightarrow{\Dir_{\xi }} H^{1}(S^{-} \otimes E) \xrightarrow{\Dir_{\xi }^{\ast }} L^{2} (S^{+} \otimes E).
$$
Exactly as in Theorem \ref{Fredholm}, the first map is Fredholm with vanishing kernel from the 
Sobolev space $H^{2}(S^{+} \otimes E)$ into $H^{1}(S^{-} \otimes E)$, both space being endowed with the $L^{2}$-inner product. 
This with the identity $Im(\Dir_{\xi })^{\bot} = Ker(\Dir_{\xi }^{\ast })$ implies that $Ker(\Delta_{\xi })=\{ 0 \}$ and 
that $Im(\Delta_{\xi })=Im(\Dir_{\xi }^{\ast })=Ker(\Dir_{\xi })^{\bot }=L^{2}(S^{+} \otimes E)$. Therefore, 
$\Delta_{\xi }$ is a bounded bijective operator from $H^{2}(S^{+} \otimes E)$ to $L^{2} (S^{+} \otimes E)$. 
By the closed graph theorem, we conclude that its inverse is also bounded. 
\end{proof}

\section{Properties of the Green's operator} \label{greensect}
\begin{defn} \label{green}
Let us call the bounded linear inverse of $\Dir_{\xi}^{\ast } \Dir_{\xi }$ provided by Theorem \ref{Hilbertisom}
the \emph{Green's operator} of the Dirac-Laplace operator, and denote it by 
$$
    G_{\xi }: L^{2} (S^{+} \otimes E) \lra H^{2} (S^{+} \otimes E). 
$$
\end{defn}
In this section we list the properties of this operator that we will need in later chapters. 
\begin{lem} \label{diaglem}
$G_{\xi }$ is diagonal with respect to the decomposition $S^{+}\otimes E= (\Omega^{0} \otimes E) \oplus (\Omega^{2}\otimes E)$. 
\end{lem}
\begin{proof}
Since $G_{\xi }$ is the inverse of $\Delta_{\xi }$, it is sufficient to prove the statement for this latter 
operator. This comes from the identity 
$$
     \Delta_{\xi } = \Dir_{\xi }^{\ast } \Dir_{\xi } = (D_{\xi }^{\ast }-D_{\xi }) (D_{\xi }-D_{\xi }^{\ast })= 
     -D_{\xi }^{\ast } D_{\xi } - D_{\xi } D_{\xi }^{\ast }, 
$$ 
which is satisfied since $D_{\xi }$ is flat. 
\end{proof}
\begin{lem} \label{estgreen}
There exist $K,K'>0$ such that for $|\xi |$ sufficiently large and for any positive spinor $\psi \in H^{1}(S^{+} \otimes E)$,
the following estimates hold: 
\begin{align} 
 \left\Vert  G_{\xi } \psi   \right\Vert_{L^{2}(\C )} & \leq K |\xi |^{-2} \Vert \psi \Vert_{L^{2}(\C )} \label{estgreen1} \\
 \left\Vert  G_{\xi } \psi   \right\Vert_{H^{1}(\C )} & \leq K' |\xi |^{-1} \Vert \psi \Vert_{L^{2}(\C )} \label{estgreen2}
\end{align}     
\end{lem}
\begin{proof}
Since by definition, for any $\psi $ the positive spinor $G_{\xi } \psi $ is the solution $\varphi $ of 
$$
    \Delta_{\xi } \varphi = \psi , 
$$
the estimates (\ref{estgreen1}) and (\ref{estgreen2}) can be rewritten respectively as 
\begin{align} 
 \left\Vert  \varphi   \right\Vert_{L^{2}(\C )} & \leq K |\xi |^{-2} \Vert \Delta_{\xi } \varphi \Vert_{L^{2}(\C )} \label{estlapl1} \\
 \left\Vert  \varphi   \right\Vert_{H^{1}(\C )} & \leq K' |\xi |^{-1} \Vert \Delta_{\xi } \varphi \Vert_{L^{2}(\C )}. \label{estlapl2}
\end{align}   
Call \emph{$\xi $-energy} of $\varphi $ over all $\C $  the quantity 
\begin{equation} \label{energyglobaldef}
      E( \xi ; \varphi ) = \int_{\C } |\nabla^{+}_{\xi } \varphi |^{2} + |\Phi_{\xi } \otimes  \varphi |^{2} |\d z |^{2}.
\end{equation}
By partial integration, the Weitzenb\"ock formula (\ref{weitzen}) and Cauchy's inequality we have 
\begin{align} \label{laplest}
    E( \xi ; \varphi ) & =  \int_{\C } \langle \varphi , \Delta_{\xi } \varphi \rangle  |\d z |^{2} \\
     & \leq \Vert \varphi \Vert_{L^{2}} \Vert \Delta_{\xi } \varphi \Vert_{L^{2}} . \notag
\end{align}
Now, as we will see from (\ref{confhiggs}), on the complementary of a finite union of disks 
$\Delta(q_{k}(\xi), \varepsilon_{0}|\xi |^{-1})$ we have the point-wise lower bound 
\begin{equation} \label{higgsptwise}
     \absl \Phi_{\xi } \otimes \varphi \absr^{2} \geq c |\xi |^{2} | \varphi  |^{2} 
\end{equation}
for some $c>0$. Furthermore, we can choose $\varepsilon_{0}$ sufficiently small so that the balls 
$\Delta(q(\xi ), 2\varepsilon_{0}|\xi |^{-1})$ are disjoint and do not meet $P$ for $|\xi|$ large. Setting 
$$
     B_{\xi }:= \bigcup_{q(\xi ) \in \Sigma_{\xi }} \Delta(q(\xi ), \varepsilon_{0}|\xi |^{-1})  
$$
we then deduce the estimation 
\begin{equation} \label{higgsint}
    \int_{\C \setminus B_{\xi } } \absl \Phi_{\xi } \otimes \varphi \absr^{2} |\d z |^{2} \geq  
    c |\xi |^{2} \int_{\C \setminus B_{\xi } } \absl \varphi \absr^{2} |\d z |^{2} .
\end{equation}
Of course, extending this inequality over the disks $\Delta(q(\xi ), \varepsilon_{0}|\xi |^{-1}) $ is not possible, 
since $\Phi_{\xi }$ has a zero in $q(\xi )$. However, the integral of $|\Phi_{\xi } \otimes \varphi |^{2} + |\nabla^{+}_{\xi } \varphi|^{2}$ does
control $|\xi |^{2}$ times that of $|\varphi |^{2}$ on the whole plane; that is, we have: 
\begin{clm} \label{estgreenclm}
There exists $c >0$ such that for $| \xi |$ sufficiently large and for any spinor $\varphi $ we have
\begin{equation} \label{globalxi2est}
      E(\xi ; \varphi ) \geq c |\xi |^{2} \int_{\C } \absl \varphi \absr^{2} |\d z |^{2} 
\end{equation}
\end{clm} 
\begin{proof}
By Kato's inequality $E(\xi ; \varphi )$ can be bounded from below by 
$$
     \int_{\C } \absl \Phi_{\xi } \otimes \varphi \absr^{2} + \absl \d | \varphi | \absr^{2} |\d z |^{2}  .
$$
By (\ref{higgsint}), it only remains to show that for any $q(\xi ) \in \Sigma_{\xi }$ this integral bounds from above 
$c |\xi |^{2} \int_{\Delta(q(\xi ), \varepsilon_{0}|\xi |^{-1})} |\varphi |^{2} |\d z |^{2} $, for some $c>0$ 
(not necessarily the same as before). But since on the annulus 
$$
    \Delta(q(\xi ), 2\varepsilon_{0}|\xi |^{-1}) \setminus \Delta(q(\xi ), \varepsilon_{0}|\xi |^{-1})
$$ 
we already have the estimation (\ref{higgsptwise}), this is just a consequence of (\ref{regineq}) applied at 
the point $q(\xi )$ instead of $p_{j}$ to the function $g=| \varphi |$, with $\varepsilon = \varepsilon_{0}|\xi |^{-1}$ 
and $\delta =0$. 
\end{proof}
By the claim and (\ref{laplest}), we have 
$$
  c |\xi |^{2} \Vert  \varphi \Vert_{L^{2}(\C )}^{2} \leq  \Vert  \varphi \Vert_{L^{2}(\C )} \Vert \Delta_{\xi } \varphi \Vert_{L^{2}(\C )}, 
$$
and after dividing both sides by $\Vert  \varphi \Vert_{L^{2}(\C )}$, we get (\ref{estlapl1}). 

Plugging (\ref{estlapl1}) into (\ref{laplest}), we obtain 
\begin{equation} \label{energylaplineq}
     E(\xi ; \varphi ) \leq K | \xi |^{-2} \Vert \Delta_{\xi } \varphi \Vert_{L^{2}(\C )}^{2}. 
\end{equation} 
On the other hand, by the definitions 
\begin{align*}
    \nabla_{\xi }^{+} = & \nabla^{+}  - \frac{\xi }{2} \d z + \frac{\bar{\xi }}{2} \d \bar{z} \\
    \Phi_{\xi } = & \Phi - \frac{\xi }{2} \d z - \frac{\bar{\xi }}{2} \d \bar{z}
\end{align*}
we obtain the point-wise bounds 
\begin{align*}
   \frac{1}{2} \absl \Phi \otimes \varphi \absr^{2} - \frac{3}{2}|\xi |^{2} | \varphi |^{2} \leq & \absl \Phi_{\xi } \otimes \varphi \absr^{2} \leq 
   2\absl \Phi \otimes \varphi \absr^{2} + |\xi |^{2} | \varphi |^{2} \\
    \frac{1}{2} \absl \nabla^{+}   \varphi \absr^{2} - \frac{3}{2}|\xi |^{2} | \varphi |^{2} \leq & \absl \nabla_{\xi }^{+}  \otimes \varphi \absr^{2} \leq 
   2\absl  \nabla^{+}  \varphi \absr^{2} + |\xi |^{2} | \varphi |^{2} 
\end{align*}
and therefore 
\begin{equation} \label{energyh1ineq}
   \frac{1}{2}\Vert \varphi \Vert_{H^{1}(\C )}^{2} - (3|\xi |^{2} +1) \Vert \varphi \Vert_{L^{2}(\C )}^{2} \leq  E(\xi ; \varphi  ) \leq 
    2\Vert \varphi \Vert_{H^{1}(\C )}^{2} + (2|\xi |^{2} +1) \Vert \varphi \Vert_{L^{2}(\C )}^{2}. 
\end{equation}
Putting together this with (\ref{energylaplineq}) and (\ref{estlapl1}), we get 
\begin{align*}
    \Vert \varphi \Vert_{H^{1}(\C )}^{2} \leq & 2E(\xi ; \varphi  ) + (6|\xi |^{2} +2) \Vert \varphi \Vert_{L^{2}(\C )}^{2} \\ 
    \leq & 2E(\xi ; \varphi  ) + 7|\xi |^{2}  \Vert \varphi \Vert_{L^{2}(\C )}^{2} \\ 
    \leq & ( 2K + 7K^{2} ) | \xi |^{-2} \Vert \Delta_{\xi } \varphi \Vert_{L^{2}(\C )}^{2}, 
\end{align*}
whence (\ref{estlapl2}). 
\end{proof}

We now investigate what happens to the Green's operator when $\xi$ is close to one of the points of $\hat{P}$. 

\begin{lem}\label{estgreenlog}
There exist $K,K'>0$ such that for $|\xi -\xi_{l}|$ sufficiently small and for any positive spinor $\psi \in H^{1}(S^{+} \otimes E)$,
the following estimates hold: 
\begin{align} 
 \left\Vert  G_{\xi } \psi   \right\Vert_{L^{2}(\C )} & \leq K |\xi -\xi_{l}|^{-2} \Vert \psi \Vert_{L^{2}(\C )} \label{estgreenlog1} \\
 \left\Vert \Dir_{\xi } G_{\xi } \psi  \right\Vert_{L^{2}(\C )} & \leq K'' |\xi -\xi_{l}|^{-1} \Vert \psi \Vert_{L^{2}(\C )} \label{estgreenlog2}
\end{align}  
\end{lem}
\begin{proof}
Analogous to Lemma \ref{estgreen}. Notice that by partial integration and the Weitzenb\"ock formula
(\ref{weitzen}) one has 
\[
      \left\Vert \Dir_{\xi }\varphi \right\Vert_{L^{2}(\C )}^{2} = E(\xi;\varphi)
\]
for any positive spinor $\varphi$. Using this and setting $G_{\xi}\psi=\varphi$ the inequalities to prove can be rewritten as 
\begin{align} 
 \left\Vert  \varphi   \right\Vert_{L^{2}(\C )} & \leq K |\xi -\xi_{l}|^{-2} \Vert \Delta_{\xi } \varphi \Vert_{L^{2}(\C )} \label{estlapllog1} \\
  E(\xi;\varphi) & \leq K'' |\xi -\xi_{l}|^{-2} \Vert \Delta_{\xi } \varphi \Vert_{L^{2}(\C )}^{2}. \label{estlapllog2}
\end{align} 
The behavior (\ref{confhiggslog}) of the Higgs field shows that outside of a finite union of disks 
$\Delta(q_{k}(\xi), \varepsilon_{0}|\xi -\xi_{l}|^{-1})$ there exists $c>0$ for which we have the point-wise lower bound 
\begin{equation} \label{higgsptwiselog}
     \absl \Phi_{\xi } \otimes \varphi \absr^{2} \geq c |\xi -\xi_{l}|^{2} | \varphi  |^{2} .
\end{equation}
It follows that denoting by $B_{\xi}$ the union of all the above mentioned disks where this estimate may fail,
we have the inequality 
\begin{equation} \label{higgsintlog}
    \int_{\C \setminus B_{\xi } } \absl \Phi_{\xi } \otimes \varphi \absr^{2} |\d z |^{2} \geq  
    c |\xi -\xi_{l}|^{2} \int_{\C \setminus B_{\xi } } \absl \varphi \absr^{2} |\d z |^{2}.
\end{equation}
It is not possible to extend this inequality to the whole plane; however, we have again 
\begin{clm} 
There exists $c >0$ such that for $|\xi -\xi_{l}|$ sufficiently small and for any spinor $\varphi$ we have 
\begin{equation} 
      E(\xi ; \varphi ) \geq c |\xi -\xi_{l}|^{2} \int_{\C } \absl \varphi \absr^{2} |\d z |^{2} 
\end{equation}
\end{clm}
\begin{proof}
Similar to Claim \ref{estgreenclm}, using Kato's inequality and (\ref{regineq}) rescaled conveniently by the 
homothety $w=(\xi -\xi_{l})z$.
\end{proof}
This together with (\ref{laplest}) then shows 
\[
   c |\xi -\xi_{l}|^{2} \Vert  \varphi \Vert_{L^{2}(\C )}^{2} \leq  \Vert  \varphi \Vert_{L^{2}(\C )} \Vert \Delta_{\xi } \varphi \Vert_{L^{2}(\C )}, 
\]
which gives us (\ref{estlapllog1}). Plugging this back into (\ref{laplest}), we obtain (\ref{estlapllog2}).

\end{proof}

\section[Exponential decay for harmonic spinors]{Exponential decay results for harmonic spinors}\label{harmspinsect}
In this section we give some analytic properties of $\Delta_{\xi }$-harmonic spinors. They will 
be needed in Section \ref{flattr}, where we study the transformed flat connection. More precisely, they
will allow us to multiply any $L^{2}$ harmonic section by exponential factor so that the result remains
in $L^{2}$. They will also be of use in the computation of the parabolic weights of the transform in 
Section \ref{parweightsect}.

First we set some further notation. Fix $\xi \in \hat{\C} \setminus \hat{P}$, and let $\varphi $ be a harmonic 
negative spinor with respect to $\Dir_{\xi } \Dir_{\xi }^{\ast }$ and $p \in \C \setminus P$ any point of the plane. 
Finally, for any spinor $\psi $ (not necessarily harmonic), call \emph{$\xi $-energy of $\psi $ in the disk 
$\Delta(p,\varepsilon )$} the quantity
\begin{equation} \label{energydef}
      E(p,\varepsilon, \xi ; \psi  ) = \int_{\Delta(p,\varepsilon )} |\nabla^{+}_{\xi } \psi  |^{2} + |\Phi_{\xi } \otimes \psi |^{2} .
\end{equation}
\begin{lem} \label{basicest}
Suppose that there exists $\varepsilon_{0} > 0$, $R>0$ and $c>0$ such that the disk 
$\Delta(p,(R+1)\varepsilon_{0} )$ is disjoint from $P$, and all of the eigenvalues of $\theta_{\xi }$ in any 
point of this disk are bounded below in absolute value by $c>0$. Under these assumptions, 
we have the inequality 
\begin{equation}
     E(p,\varepsilon_{0}, \xi ; \varphi ) \leq e^{-2 c R \varepsilon_{0} } \left( 2\Vert \varphi \Vert_{H^{1}(\C)}^{2} 
     + (2|\xi |^{2}+1) \Vert \varphi \Vert_{L^{2}(\C)}^{2} \right) . 
\end{equation}
\end{lem}
\begin{proof}
Denote by $C(p,r )$ the boundary of $\Delta(p,r )$, and by $\frac{\partial }{\partial n}$ an outward-pointing
unit normal vector to it. Stokes' formula gives 
\begin{align*}
     E(p,r, \xi ; \varphi )= \int_{\Delta(p,r )} &
     \left( (\nabla^{+}_{\xi })^{\ast } \nabla^{+}_{\xi }\varphi + (\Phi_{\xi } \otimes)^{\ast } \Phi_{\xi } \otimes \varphi ,\varphi  \right) \\ + 
     & \int_{C(p,r )} \left( \left(\nabla^{+}_{\xi } \right)_{\frac{\partial }{\partial n}} \varphi ,\varphi  \right) r \d \theta .
\end{align*}
Since $\varphi $ is $\Delta_{\xi }$-harmonic, the Weitzenb\"ock formula (\ref{weitzen}) implies that the first term on 
the right-hand side vanishes. Therefore, by the tic-tac-toe inequality, we have 
$$
     E(p,r, \xi ; \varphi ) \leq \frac{1}{2} \int_{C(p,r )} \frac{1}{c} \absl \nabla^{+}  \varphi \absr^{2} 
     + c|\varphi |^{2}  r \d \theta .
$$
On the other hand, we have 
$$
    \frac{\d E(p,r, \xi ; \varphi )}{\d r} = \int_{C(p,r )} \absl \nabla^{+}_{\xi } \varphi \absr^{2} 
    + |\Phi_{\xi } \otimes \varphi|^{2} r \d \theta. 
$$
By assumption, for $r \leq (R+1)\varepsilon_{0}$ we have the estimate 
$$
     \int_{C(p,r )} |\Phi_{\xi } \otimes \varphi|^{2} r \d \theta \geq c^{2} \int_{C(p,r )} | \varphi|^{2} r \d \theta . 
$$
Putting together these estimates, we see that 
$$
     \frac{\d E(p,r, \xi ; \varphi )}{\d r} \geq 2c E(p,r, \xi ; \varphi ),
$$
whence 
$$
     \frac{\d \log E(p,r, \xi ; \varphi )}{\d r} \geq 2c.
$$
Integrating this inequality from $r=\varepsilon_{0}$ to $r= (R+1)\varepsilon_{0}$, 
we obtain
\begin{align*}
     \log E(p,\varepsilon_{0} , \xi ; \varphi ) & \leq 2c [\varepsilon_{0} - (R+1)\varepsilon_{0}] 
                                         + \log E(p,(R+1)\varepsilon_{0} , \xi ; \varphi ). 
\end{align*}
Taking exponential of both sides, we get
\begin{align*}
     E(p,\varepsilon_{0} , \xi ; \varphi ) & \leq e^{-2c R \varepsilon_{0}} 
                                         E(p,(R+1)\varepsilon_{0} , \xi ; \varphi ) \\
                                   & \leq e^{-2c R \varepsilon_{0}} E( \xi ; \varphi ), 
\end{align*}
and we conclude using (\ref{energyh1ineq}). 
\end{proof}

Next, we use the above lemma to obtain exponential decay results in terms of $\xi$  for the energy of 
harmonic spinors when $\xi$ is large, first in a fixed disk of $\C$ away from the singularities $P$, 
then near infinity in $\C$. In the first case, the statement is as follows. 
\begin{lem} \label{expdecord}
Let $p \in {\C } \setminus P$ be arbitrary, and let $\varepsilon_{0}>0$ be such that the distance between  
$p$ and $P$ is at least $3\varepsilon_{0}$. Then for $|\xi |$ sufficiently large we have the estimate
$$
     \Vert \varphi \Vert^{2}_{H^{1}(\Delta(p, \varepsilon_{0}) )} \leq e^{- \varepsilon_{0} |\xi |/3 } \Vert \varphi \Vert_{H^{1}(\C)}^{2}
$$ 
for any $\Delta_{\xi }$-harmonic spinor $\varphi $.
\end{lem} 
\begin{proof}
Since $p$ is away from $P$, in the Higgs field $\theta_{\xi } = \theta - \xi\d z/2$ the term $\theta$ is bounded on $\Delta(p,2\varepsilon_{0})$. 
Therefore, if $|\xi |$ is sufficiently large, then the eigenvalues of $\theta_{\xi }$  on this disk are bounded below in absolute
value by $|\xi | /4$. Apply Lemma \ref{basicest} with $R=1$ and $c=|\xi | /4$ to get 
\begin{align*}
   E(p,\varepsilon_{0}, \xi ; \varphi ) & \leq e^{- \varepsilon_{0} |\xi |/2 } \left( 2\Vert \varphi \Vert_{H^{1}(\C)}^{2} +
   (2|\xi |^{2} +1) \Vert \varphi \Vert_{L^{2}(\C)}^{2} \right) \\ 
    & \leq 5 e^{- \varepsilon_{0} |\xi |/2 }  |\xi |^{2} \Vert \varphi \Vert_{H^{1}(\C)}^{2} \\ 
    & \leq \frac{1}{33} e^{- \varepsilon_{0} |\xi |/3 } \Vert \varphi \Vert_{H^{1}(\C)}^{2}
\end{align*}
for $\xi $ sufficiently large. On the other hand, we have 
\begin{align} \label{h1energyineq}
    \Vert \varphi \Vert^{2}_{H^{1}(\Delta(p, \varepsilon_{0}) )} = & 
    \int_{\Delta(p, \varepsilon_{0})} |\varphi |^{2} + \absl \nabla^{+} \varphi \absr^{2} + \absl \Phi \otimes \varphi \absr^{2} \notag \\
    \leq & \int_{\Delta(p, \varepsilon_{0})} 2 |\xi |^{2} |\varphi |^{2} + \absl \nabla^{+}_{\xi } \varphi \absr^{2} + \absl \Phi_{\xi } \otimes \varphi \absr^{2}\\
    \leq & 33 \; E(p,\varepsilon_{0}, \xi ; \varphi ) \notag,
\end{align}
where the last line is a consequence of $|\Phi_{\xi } \otimes \varphi |^{2} \geq |\xi |^{2}| \varphi |^{2} /16 $ in
$\Delta(p,\varepsilon_{0})$. Putting together these two estimates, we get the lemma. 
\end{proof}
In the second case, we have the following statement.
\begin{lem} \label{expdecinf}
For any $\xi \notin \hat{P}$ there exists $R_{0}= R_{0}(\xi )>0$, $K= K(\xi )>0$ and $c=c(\xi )>0$ such that 
for any $\Delta_{\xi }$-harmonic spinor $\varphi $ and all $R>R_{0}$ the following estimate holds:
\begin{equation*}
       \Vert \varphi \Vert_{H^{1}(\C \setminus \Delta(0,2R ))}^{2} \leq K e^{-R c} \Vert \varphi \Vert_{H^{1}(\C )}^{2}.
\end{equation*}
Furthermore, if $|\xi |$ is sufficiently large, we can choose $c= |\xi |/3$ and $R_{0},K$ constants 
independent of $\xi $. 
\end{lem}
\begin{proof}
The proof is an amalgam of that of Lemmata \ref{basicest} and \ref{expdecord}. 
Define the \emph{$\xi $-energy at infinity} of a spinor by the integral 
\begin{equation} \label{energyinfdef}
      E(\infty , R , \xi ; \varphi ) = \int_{ \C \setminus \Delta(0,R )} |\nabla^{+}_{\xi } \varphi |^{2} + |\Phi_{\xi } \otimes  \varphi |^{2} .
\end{equation}
Choose $R_{0}>0$ and $c_{0}$ such that for $|z|>R_{0}$ the eigenvalues of $\theta_{\xi }(z)$ are all bigger in
absolute value then $c_{0}$. Clearly, such a choice is possible because $\xi \notin \hat{P}$. Moreover, 
for $|\xi |$ sufficiently large one can put $c_{0} = |\xi | /4$ and $R_{0}$ a constant only depending on 
the initial data $\theta $. For $r \geq R_{0}$, we have the estimate 
$$
      - E(\infty , r , \xi ; \varphi ) \geq -\frac{1}{2} \int_{C(0,r )} \frac{1}{c_{0} } \absl \nabla^{+}_{\xi }  \varphi \absr^{2} +
      c_{0} |\varphi |^{2}  r \d \theta .
$$
On the other hand, we have 
$$
    \frac{\d E(\infty ,r, \xi ; \varphi )}{\d r} = - \int_{C(0,r )} \absl \nabla^{+}_{\xi }  \varphi \absr^{2} 
    + |\Phi_{\xi } \otimes \varphi|^{2} r \d \theta. 
$$
By assumption, we have also
$$
     \int_{C(0,r )} |\Phi_{\xi } \otimes \varphi|^{2} r \d \theta \geq c_{0}^{2} \int_{C(0,r )} | \varphi|^{2} r \d \theta . 
$$
Putting together these estimates, we see that for $r \geq  R_{0}$
$$
     \frac{\d E(\infty ,r , \xi ; \varphi )}{\d r} \leq  - 2c_{0} E(\infty ,r , \xi ; \varphi ),
$$
whence 
$$
     \frac{\d \log E (\infty ,r, \xi ; \varphi )}{\d r} \leq - 2c_{0}.
$$
Integrating this inequality from $R$ to $2R$ and using (\ref{energyh1ineq}), we obtain 
\begin{align*}
     E(\infty , 2R, \xi ; \varphi ) \leq & E( \xi ; \varphi ) e^{-R c_{0}} \\
     \leq & (|\xi |^{2}+3) e^{- R c_{0} } \Vert \varphi \Vert_{H^{1}(\C)}^{2}.
\end{align*} 
On the other hand,  
\begin{align*}
     E(\infty , 2R, \xi ; \varphi ) & \geq \int_{\C \setminus \Delta(0,2R)} \absl \Phi_{\xi } \otimes \varphi \absr^{2} \\
     & \geq c_{0}^{2} \int_{\C \setminus \Delta(0,2R)} |\varphi |^{2} 
\end{align*}
implies 
$$
     K_{0} E(\infty , 2R, \xi ; \varphi ) \geq \Vert \varphi \Vert_{H^{1}(\C \setminus \Delta(0,2R ))}^{2}
$$
for some $K_{0}>0$. This gives the lemma for $\xi $ in a finite region. The case of $|\xi |$ large also 
follows noting that $K$ depends at most polynomially on $\xi $. 
\end{proof} 

Since a $\Delta_{\xi }$-harmonic spinor is subharmonic in the usual sense, the above results also imply point-wise 
exponential decay on harmonic spinors: 
\begin{lem} \label{ptwiseexpdec}
Suppose $R>R_{0}$. Then there exists $K,c>0$ such that for any $|z | > 2R+1$ and any $\Delta_{\xi }$-harmonic spinor 
$\varphi $ we have 
$$
      |\varphi (z) | \leq K e^{-Rc} \Vert \varphi \Vert_{H^{1}(\C )}^{2}.
$$
\end{lem}
\begin{proof}
Because of the condition $|z | > 2R+1$, the disk $\Delta(z,1) $ centered at $z$ of radius $1$ is contained 
in $\C \setminus \Delta(0,2R)$. On the other hand, by subharmonicity of $\varphi $ with respect to the usual 
Laplace operator, we have 
\begin{align*}
     |\varphi (z) | & \leq K_{0} \int_{\Delta(z,1)} |\varphi (w)| |\d w |^{2} \\
     &  \leq K_{1} \left( \int_{\Delta(z,1)} |\varphi (w)|^{2} |\d w |^{2} \right)^{1/2} \\
     &  \leq K_{1} \left( \int_{\C \setminus \Delta(0,2R)} |\varphi (w)|^{2} |\d w |^{2} \right)^{1/2}
\end{align*}
We conclude using Lemma \ref{expdecinf}. 
\end{proof}

\chapter{The transform of the integrable connection}
In this chapter, we define the transformed parabolic integrable connection induced by the deformation $D_{\xi }$. 
First, in Section \ref{flattr}, we define the underlying flat bundle; then in Section \ref{flatext} 
we show that its behavior at infinity verifies appropriate asymptotic conditions. This then allows us 
to apply the results of \cite{Biq-Boa} in order to define an extension into a parabolic integrable connection over 
the singularity at infinity; the same thing for other singularities follows from \cite{Sim}. 

Before starting these points, we need however
to introduce some notation. Recall first that $\hat{P}$ was defined as the set $\{ \xi_{1} , \ldots , \xi_{n'} \}$ 
of eigenvalues of the second-order term of $D$ at infinity. Let $\hat{H} \ra \hat{\C} \setminus \hat{P}$
denote the trivial Hilbert bundle with fibers $L^{2}(\C , S^{-}\otimes E)$. 
By Theorem \ref{laplker}, the transformed bundle $\hat{E}$ 
can be given as the vector bundle whose fiber over $\xi \in \C
\setminus \hat{P}$ is the kernel of the adjoint Dirac operator
$(\Dir_{\xi })^{\ast }$. By the same theorem, such an element is also $\Delta_{\xi }$-harmonic. 
Now remark that on the bundle $\hat{H}$ there exists a
hermitian metric $\langle.,.\rangle $ which is canonical once a hermitian metric
$h(.,.)$ is fixed on $E$: for any two elements 
$\hat{f}_{1},\hat{f}_{2} \in \hat{H}_{\xi }=L^{2}(\C , S^{-}\otimes E)$, 
it is defined by the $L^{2}$ inner product  
$$
     \langle \hat{f}_{1},\hat{f}_{2}\rangle= \int_{\C } h(\hat{f}_{1},\hat{f}_{2}) |\d z|^{2}. 
$$
Moreover, the trivial connection $\hat{\d}$ on the bundle $\hat{H}$ is
unitary with respect to this metric. Let $\hat{\pi}_{\xi }$ denote orthogonal
projection of $\hat{H}_{\xi }$ onto the subspace $\hat{E}_{\xi }$, and $ i $ the
inclusion $\hat{E} \inc \hat{H}$.
\begin{defn} \label{deftrhm}
We call \emph{transformed hermitian metric} the fiber metric $\hat{h}$ on $\hat{E}$ which 
is equal on the fiber $\hat{E}_{\xi }$ to the restriction of the above defined $L^{2}$ scalar product $\langle.,.\rangle$ 
to the subspace $\hat{E}_{\xi } \subset L^{2}(\C , S^{-}\otimes E)$. 
\end{defn}

\section[Construction]{Construction of the transformed flat connection}\label{flattr}
In this section we show that the transformed bundle admits an integrable
connection, which is determined only by the deformation $D_{\xi }$ and the . 
First, we describe its intrinsic construction, then we give it in terms
of an explicit formula. 

\subsection{Intrinsic definition} 
Defining a flat connection is equivalent to giving a basis of parallel
sections on a disk $B_{0}$ around each point $\xi_{0} \in \hat{\C} \setminus
\hat{P} $. Given this, in order to see that it defines
indeed a flat connection, one only needs to prove that the transition 
matrices on $B_{0} \cap B_{1}$ between two such bases (corresponding to 
points $\xi_{0}$ and $\xi_{1}$) are constant. 

So suppose 
$\xi_{0} \in \hat{\C} \setminus \hat{P}$, and let $\hat{f}_{1}( z),\ldots,\hat{f}_{\hat{r}}( z)$
be a basis of the vector space $\hat{E}_{\xi_{0}}$. On the basis of Lemma \ref{ptwiseexpdec}, 
for $\varepsilon_{0} = \varepsilon_{0} (\xi _{0}) >0 $ sufficiently small, the expressions 
\begin{align}
     \hat{f}_{j}(\xi;z) & = \hat{\pi}_{\xi}(e^{(\xi - \xi_{0})z }\hat{f}_{j}(z)) \in \hat{E}_{\xi } \label{extf} 
\end{align}
make sense for $\xi $ on the ball $B_{0}=B(\xi _{0}, \varepsilon_{0})$ of radius $\varepsilon_{0}$ centered 
at $\xi_{0}$. Therefore, (restricting $\varepsilon_{0}$ if necessary), they define an extension of the basis 
$\hat{f}_{1},\ldots,\hat{f}_{\hat{r}}$ of the vector space $\hat{E}_{\xi_{0}}$ to a trivialisation of the bundle 
$\hat{E}$ over $B_{0}$.
\begin{prop} The family of sections (\ref{extf}) for all $\xi_{0} \in \hat{\C} \setminus \hat{P}$, 
for $j \in \{ 1,\ldots,\hat{r} \}$, and for all $\xi \in B_{0}$ define a local system for  a flat connection 
$\hat{D}$ on $\hat{E} \ra \hat{\C} \setminus \hat{P}$. 
\end{prop}
\begin{defn} \label{deftrflat}
We will call $\hat{D}$ the \emph{transformed flat connection} on $\hat{\C} \setminus \hat{P}$. 
\end{defn}

\begin{proof} (Proposition)
Let $\tilde{\xi}_{0} \neq \xi_{0}$ be another point of $ \hat{\C} \setminus \hat{P}$, and 
$\hat{g}_{1}(z),\ldots,\hat{g}_{\hat{r}}(z)$ be a basis for the vector space $\hat{E}_{\tilde{\xi}_{0}}$.
According to (\ref{extf}), the local trivialisation of $\hat{E}$ near $\tilde{\xi}_{0}$ we need to 
consider is then $\hat{g}_{1}(\xi ), \ldots,\hat{g}_{\hat{r}}(\xi )$, with 
\begin{align}
     \hat{g}_{l}(\xi;z) & = \hat{\pi}_{\xi}(e^{(\xi - \tilde{\xi}_{0})z }\hat{g}_{l}(z)) \label{extg} 
\end{align}
for $\xi $ in a small disk $\tilde{B}_{0}$ around $\tilde{\xi}_{0}$. In order to show that the local 
bases (\ref{extf}) and (\ref{extg}) define indeed a local system, 
we need to show that the transition matrices $m(\xi )$ between them are 
independent of the point $\xi \in B_{0} \cap \tilde{B}_{0}$. 
We will make use of the following: 
\begin{lem} \label{trans} For any $\xi, \xi' \in B_{0}$, and any $k_{0} \in
  ker(D_{\xi_{0} }|S^{-}\otimes E)$ we have 
\[
     \hat{\pi}_{\xi' } \left( e^{(\xi' - \xi )z } \hat{\pi}_{\xi}( e^{(\xi -\xi_{0})z} k_{0}(z) )
     \right) =  \hat{\pi}_{\xi' } (  e^{(\xi' -\xi_{0})z} k_{0}(z) ). 
\]
\end{lem}
\begin{proof} (Lemma) Set $k_{\xi }(z)=e^{(\xi -\xi_{0})z} k_{0}(z)$;
we need to prove that 
\[
  \hat{\pi}_{\xi' }[  e^{(\xi' - \xi )z} \hat{\pi}_{\xi}( k_{\xi }(z) ) ] = \hat{\pi}_{\xi' }(  e^{(\xi' - \xi )z}  k_{\xi }(z)  ), 
\]
or equivalently that 
\[
     \hat{\pi}_{\xi' }[ e^{(\xi' - \xi )z} ( \Id - \hat{\pi}_{\xi} ) ( k_{\xi } )  ] = 0, 
\]
which is still equivalent to 
\begin{equation}   \label{orthrel}
      e^{(\xi' - \xi )z} ( \Id - \hat{\pi}_{\xi} ) ( k_{\xi } ) \bot \hat{E}_{\xi '}. 
\end{equation}
Since $\hat{\pi}_{\xi}$ is orthogonal projection to $\hat{E}_{\xi }$, we have  
\begin{equation} \label{eorth}
     ( \Id - \hat{\pi}_{\xi} ) ( k_{\xi } ) \in \hat{E}_{\xi }^{\bot }.   
\end{equation}
Moreover, observe that for $\xi_{0}$ and $\xi$ fixed, the relation 
\begin{equation} \label{imgauge}
    e^{(\xi - \xi_{0})z}.D_{\xi_{0} } = D_{\xi_{0} } -  (\xi  - \xi_{0}  )\d z \land  
   = D_{\xi },  
\end{equation}
holds, and so 
\begin{equation} \label{impluse}
      k_{\xi } = e^{(\xi -\xi_{0})z} k_{0} \in e^{(\xi -\xi_{0})z} ker(D_{\xi_{0} }) \subset
      ker(D_{\xi }) = Im(D_{\xi }^{\ast })^{\bot }= Im(D_{\xi }) \oplus \hat{E}_{\xi }.
\end{equation}
From (\ref{eorth}) and (\ref{impluse}) it follows that $( \Id - \hat{\pi}_{\xi} ) k_{\xi } \in Im(D_{\xi })$. 
Now using (\ref{imgauge}) for $(\xi'- \xi ) $ instead of $(\xi-\xi_{0})$,
we deduce that $e^{(\xi' - \xi )z}( \Id - \hat{\pi}_{\xi} ) k_{\xi } \in Im( D_{\xi' } )$, 
whence (\ref{orthrel}). This finishes the proof of the lemma. 
\end{proof}

Let us now come back to the study of the transition matrix: let $\xi , \xi' \in B_{0} \cap \tilde{B}_{0} $, 
and suppose we have 
\begin{equation} \label{tm}
    \hat{f}_{j}(\xi ) = \sum_{l=1}^{\hat{r}} m_{jl} \hat{g}_{l}(\xi ), 
\end{equation}
where $(m_{jl})$ is the transition matrix between the two bases at the point $\xi $. 
Lemma \ref{trans} means that for $|\xi - \xi '|$ sufficiently small, we have 
\begin{align}
         \hat{f}_{j}(\xi' ) & = \hat{\pi}_{\xi '}( e^{(\xi' - \xi )z} \hat{f}_{j}(\xi )  ) \label{bpch} \\
         \hat{g}_{l}(\xi' ) & = \hat{\pi}_{\xi '}( e^{(\xi' - \xi )z} \hat{g}_{l}(\xi )  ). \label{bpchg} 
\end{align}
Now plugging (\ref{tm}) into (\ref{bpch}), then using (\ref{bpchg}) we obtain 
\begin{align*}
    \hat{f}_{j}(\xi' ) & = \hat{\pi}_{\xi '} \left( e^{(\xi' - \xi )z}  \sum_{l=1}^{\hat{r}} m_{jl}
      \hat{g}_{l}(\xi ) \right) \\
    & = \sum_{l=1}^{\hat{r}} m_{jl} \hat{\pi}_{\xi '}( e^{(\xi' - \xi )z}\hat{g}_{l}(\xi )) \\ 
    & = \sum_{l=1}^{\hat{r}} m_{jl}  \hat{g}_{l}(\xi' ),    
\end{align*}
so the transition matrix at the point $\xi '$ is the same as the one at $\xi $, 
whence we obtain the Proposition. 
\end{proof}


\subsection{Explicit description} 
We now give an explicit formula for the flat connection constructed above. 
In the sequel we follow \cite{Jardim}. First define a unitary connection on $\hat{E}$ with respect to 
the transformed hermitian metric by 
\begin{equation} \label{trunit}
      \hat{\nabla } = \hat{\pi}_{\xi} \circ \hat{\d} \circ i . 
\end{equation}
The fact that this connection is indeed $\hat{h}$-unitary can be seen as follows: 
let $f,g \in \Gamma(\hat{E})$ be local sections around $\xi_{0}$, then from
orthogonality of $\hat{\pi}_{\xi}$ to $\hat{E}$ with respect to the norm $\langle .,. \rangle $
we have in $\xi_{0}$
\begin{align*}
    \hat{\d} (\hat{h}(\hat{f},\hat{g})) & =  \hat{\d} \langle \hat{f},\hat{g} \rangle = \langle \hat{\d } \hat{f},\hat{g} \rangle + 
    \langle \hat{f},\hat{\d }\hat{g} \rangle  \\ & = \langle \hat{\nabla  } \hat{f},\hat{g} \rangle + \langle  \hat{f},\hat{\nabla } \hat{g} \rangle 
    = \hat{h}(\hat{\nabla  } \hat{f},\hat{g}) + \hat{h}( \hat{f},\hat{\nabla } \hat{g}), 
\end{align*}
where $\hat{\d}$ stands for exterior differentiation of functions along the coordinate $\xi$ as 
well as for the trivial connection with respect to $\xi$ on the trivial Hilbert bundle $\hat{H}$. 
Finally, we define an endomorphism-valued $(1,0)$-form 
(a candidate to be a transformed Higgs field) by mapping a 
$\Delta_{\xi }$-harmonic section $\hat{f}(\xi;z)$ to 
\begin{equation} \label{trhiggs}
     \hat{\theta }_{\xi }(\hat{f}(\xi;z))= -\frac{1}{2}\hat{\pi}_{\xi}(z\hat{f}(\xi;z)) \d \xi, 
\end{equation}
where $\d \xi $ stands for the standard generator of the holomorphic 
$(1,0)$-forms on $\hat{\C}$. This field will indeed be holomorphic provided
that the original metric $h$ is harmonic (see Section \ref{sectharm}). 
\begin{prop} \label{explflattr}
The connection $\hat{\nabla }+ 2\hat{\theta }$ is equal to the transformed flat connection $\hat{D}$ defined above. 
\end{prop}
\begin{proof}
We need to show that for all $\xi_{0}$ and all $f(z) \in \hat{E}_{\xi_{0}}$, the
local $\hat{D}$-parallel section in $\xi \in B_{0}$ given by 
\begin{equation} \label{parallel}
     \hat{f}(\xi;z)= \hat{\pi}_{\xi}(e^{(\xi - \xi_{0})z }\hat{f}(z))
\end{equation}
is parallel in $B_{0}$ with respect to $\hat{\nabla }+ 2\hat{\theta }$. First, let us check it in $\xi_{0}$:
\begin{align*}
      ((\hat{\nabla }+ 2\hat{\theta })\hat{f})(\xi_{0})  = 
     \hat{\pi}_{\xi_{0}}[(\hat{\d}\hat{f})(\xi_{0} )- z\hat{f}(\xi_{0} )\d \xi ] .
\end{align*}
We observe that by (\ref{parallel}) we have 
\begin{align*}
    (\hat{\d} \hat{f})(\xi_{0} )= (\hat{\d} \hat{\pi}_{\xi})_{\xi_{0}} \hat{f}(\xi_{0} ) + 
    \hat{\pi}_{\xi_{0} }( z \hat{f}(\xi_{0}) \d \xi  ),  
\end{align*}
hence 
\[
       ((\hat{\nabla }+ 2\hat{\theta })\hat{f})(\xi_{0})  = \hat{\pi}_{\xi_{0}} [ (\hat{\d} \hat{\pi}_{\xi})_{\xi_{0}} \hat{f}(\xi_{0} )] .
\]
Now $\hat{\pi}_{\xi} \circ \hat{\pi}_{\xi} = \hat{\pi}_{\xi}$ implies 
\[
      \hat{\d} \hat{\pi}_{\xi} \circ \hat{\pi}_{\xi} + \hat{\pi}_{\xi} \circ \hat{\d} \hat{\pi}_{\xi} = \hat{\d} \hat{\pi}_{\xi}, 
\]
therefore 
\[
       \hat{\pi}_{\xi_{0}} [ (\hat{\d} \hat{\pi}_{\xi})_{\xi_{0}} \hat{f}(\xi_{0} )] =  (\hat{\d} \hat{\pi}_{\xi})_{\xi_{0}}
       \circ  ( \Id - \hat{\pi}_{\xi_{0}} ) \hat{f}(\xi_{0} ) =0, 
\]
since $\hat{\pi}_{\xi_{0}}$ is the projection to $ \hat{E}_{\xi_{0}} $ and 
$\hat{f}(\xi_{0}) \in \hat{E}_{\xi_{0}}$. 

Next, fix an arbitrary $\xi \in B_{0}$. Then, as we have just shown,
the local section defined for $|\xi ' - \xi |$ sufficiently small by 
$$
     \hat{f}'(\xi ') = \hat{\pi}_{\xi '} (e^{(\xi' - \xi )z}\hat{f}(\xi;z ))
$$ 
is parallel in $\xi $ (compare with (\ref{parallel}), setting $\xi_{0}= \xi , \xi = \xi '$). 
But Lemma \ref{trans} tells us that the local sections $\hat{f}'$ and $\hat{f}$ coincide in a 
neighborhood of $\xi$; in particular $\hat{f}$ is parallel in $\xi $. 
\end{proof}

The following is now immediate: 

\begin{prop} \label{unittr} The unitary part of the transformed flat connection $\hat{D}$ is 
$$   
     \hat{D}^{+}= \hat{\nabla } + \hat{\theta} - \hat{\theta}^{\ast } = 
     \hat{\pi }_{\xi } \circ (\hat{\d } - \frac{1}{2}z \d \xi \land + \frac{1}{2} \bar{z} \d \bar{\xi} \land ).  
$$ 
\end{prop}
\begin{defn} We will call the above unitary connection $\hat{D}^{+}$ the \emph{transformed 
unitary connection}. The covariant derivative associated to it will be denoted $\hat{\nabla }^{+} $. 
\end{defn}
\begin{rk} The fact that the formula for the transformed unitary connection involves 
extra multiplication terms by $z$ and $\bar{z}$ compared to the usual formulae 
of other Nahm transforms is an artifact: as we will see in the next chapter, the transform admits 
an interpretation from the point of view of Higgs bundles, in which the formula for 
the transformed unitary connection agrees with the usual one. 
\end{rk}

\section[Extension over the singularities]{Extension over the singularities}\label{flatext}
At this point, it should be pointed out that a priori we have no guarantee that the constructed flat  
connection is indeed of the form given in \cite{Biq-Boa} (and therefore extends nicely over the singularities); 
that is, in an orthonormal basis with respect to its harmonic metric it is not necessarily the model
(\ref{polard}) up to a perturbation described in (\ref{pertj}) and (\ref{pertinf}). However, there is a theorem of O. Biquard 
and M. Jardim which allows us to show that this is the case. Namely,  Theorem 0.1 of \cite{Biq-Jar} states the following: 
\begin{thm} \label{biqjar}
Let $\tilde{A}$ be an $SU(2)$-instanton on $\R^{4}$, invariant with respect to 
the additive subgroup $\Z \frac{\partial}{\partial x_{3}} \oplus \Z \frac{\partial}{\partial x_{4}}$, and suppose that its curvature 
$F_{\tilde{A}}$ has quadratic decay at infinity (that is, $\absl F_{\tilde{A}} \absr = O(r^{-2})$, where
$r^{2}=x_{1}^{2}+x_{2}^{2}$). Then there exists a gauge near infinity in which $\tilde{A}$ is asymptotic 
to the following model: 
\begin{align*}
      \tilde{A}_{0} = \d + i \Big{(} \lambda_{1} \d x_{3} + & \lambda_{2} \d x_{4} + (\mu_{1} \cos \theta -\mu_{2} \sin \theta ) \frac{\d x_{3}}{r} \\
              & + (\mu_{1} \sin \theta + \mu_{2} \cos \theta ) \frac{\d x_{4}}{r} + \alpha \d \theta 
              \Big{)},
\end{align*}
where $z=r e^{i\theta} $ are coordinates for the $(x_{1},x_{2})$-plane. Moreover, the difference $a$ between 
$\tilde{A}$ and this model satisfies 
$$
     |a|=O(r^{-1-\delta }), \; \absl \nabla_{\tilde{A}_{0}} a \absr = O(r^{-2-\delta }). 
$$
\end{thm}
In order to be able to apply this result to our case, consider the Euclidean space $(\R^{4})^{\ast }$ spanned by 
orthonormal vectors $\frac{\partial}{\partial x^{\ast }_{j}}$ for $j=1,2,3,4$, and identify the subspace spanned by 
$\frac{\partial}{\partial x^{\ast }_{1}}$ and  $\frac{\partial}{\partial x^{\ast }_{2}}$ with the line $\hat{\C}$ with complex 
coordinate $\xi$ underlying $\hat{\D}$. By the results of \cite{Hit}, $\hat{D}$ then induces 
an instanton $\tilde{A}$ on $(\R^{4})^{\ast }$ with singularities, invariant with respect to the subspace 
$\R \frac{\partial}{\partial x^{\ast }_{3}} \oplus \R \frac{\partial}{\partial x^{\ast }_{4}}$. 
In particular, $\tilde{A}$ is invariant with respect to $\Z \frac{\partial}{\partial x^{\ast }_{3}} \oplus \Z \frac{\partial}{\partial x^{\ast }_{4}}$,
so Theorem \ref{biqjar} can be applied to it, provided that its curvature has quadratic decay. 
In order to have an explicit description of $\tilde{A}$ and its curvature, remember that $\hat{D}$ decomposes as 
$$
     \hat{D} = \hat{\nabla}^{+} + \hat{\theta} + \hat{\theta}^{\ast}, 
$$
where $\hat{\nabla}^{+}$ is the transformed unitary connection, $\hat{\theta}$ the field defined in (\ref{trhiggs}) and 
$\hat{\theta}^{\ast}$ its adjoint with respect to the harmonic metric of $\hat{D}$. Now as we will see in Section 
\ref{sectharm}, this harmonic metric is in fact the transformed hermitian metric $\hat{h}$ given in Definition
\ref{deftrhm}. The unitary part of $\hat{D}$ decomposes further into its $(1,0)$-  and $(0,1)$-part:  
$$
    \hat{\nabla}^{+} = (\hat{\nabla}^{+})^{1,0} + (\hat{\nabla}^{+})^{0,1}.
$$
Finally, we write $\hat{\vartheta}$ for the endomorphism-part of $\hat{\theta}$:
$$
     \hat{\theta} = \hat{\vartheta} \d \xi .
$$
The instanton over $(\R^{4})^{\ast}$ corresponding to $\hat{D}$ is then given by the formula 
$$
    \tilde{A} = \hat{\nabla}^{+} + \Re \hat{\vartheta} \d x^{\ast}_{3} + \Im \hat{\vartheta} \d x^{\ast}_{4},
$$
where we recall that 
$$
    \frac{\partial}{\partial \xi} = \frac{1}{2} \left( \frac{\partial}{\partial x^{\ast}_{1}} -  \frac{\partial}{\partial x^{\ast}_{2}} \right)
$$
is the natural complex coordinate of $\hat{\C}$, and the connection $\hat{\nabla}^{+}$ on $(\R^{4})^{\ast}$ acts as
$\hat{\nabla}^{+}$ along $\hat{\C}$ and as the trivial connection along $\R \frac{\partial}{\partial x^{\ast }_{3}} \oplus \R \frac{\partial}{\partial x^{\ast }_{4}}$. 
Furthermore, as it can be seen from the results in Section 1 of \cite{Hit}, we then have the formula 
\begin{align}
      F_{\tilde{A}} = - & [ \hat{\vartheta} , \hat{\vartheta}^{\ast }] (\d x^{\ast}_{1} \land \d x^{\ast}_{2} + \d x^{\ast}_{3} \land \d x^{\ast}_{4}) \notag \\
       & + (\hat{\nabla}^{+})_{x^{\ast}_{1}} \Re \hat{\vartheta} (\d x^{\ast}_{1} \land \d x^{\ast}_{3} - \d x^{\ast}_{2} \land \d x^{\ast}_{4})  \\
       & + (\hat{\nabla}^{+})_{x^{\ast}_{1}} \Im \hat{\vartheta} (\d x^{\ast}_{1} \land \d x^{\ast}_{4} + \d x^{\ast}_{2} \land \d x^{\ast}_{3}) \notag,   
\end{align}
where we have written $(\hat{\nabla}^{+})_{x^{\ast}}$ to denote the action of the unitary connection in the 
$\frac{\partial}{\partial x^{\ast}}$-direction. Hence, before we can apply Theorem \ref{biqjar} we need to check the following: 
\begin{thm} \label{comdec} There exists a constant $K>0$ such that the commutator $ [\hat{\vartheta},\hat{\vartheta}^{\ast}]$
is bounded by $K |\xi |^{-2}$ as $\xi \ra \infty$. The same estimation holds for $\hat{\nabla }^{+} \hat{\vartheta}$. 
\end{thm}
\begin{proof}
We start with the case of the commutator. 
Let $\hat{f}(\xi;z) \in \hat{E}_{\xi } = Ker (\Dir_{\xi })^{\ast } $ be arbitrary; we wish to show the estimate 
$$
     \absl [\hat{\vartheta }, \hat{\vartheta }^{\ast } ] \hat{f}(\xi) \absr_{\hat{h}} \leq K |\xi |^{-2} | \hat{f}(\xi)|_{\hat{h}},
$$
with $K$ independent of $\hat{f}$ and of $\xi $. Recall the well-known formula from Hodge theory: 
\begin{equation} \label{projform}
     \hat{\pi }_{\xi } = \Id - \Dir_{\xi } G_{\xi } \Dir_{\xi }^{\ast }. 
\end{equation}
Using this, we obtain 
\begin{align}
     [\hat{\vartheta }, \hat{\vartheta }^{\ast } ] \hat{f}(\xi) & = 
     -\frac{1}{2}\hat{\pi }_{\xi }(z \hat{\pi }_{\xi }(\bar{z} \hat{f}(\xi)) - 
     \bar{z} \hat{\pi }_{\xi }(z \hat{f}(\xi))) \notag \\  
     & = \frac{1}{2} \hat{\pi }_{\xi } ( z \Dir_{\xi } G_{\xi } \Dir_{\xi }^{\ast } (\bar{z} \hat{f}(\xi)) 
     - \bar{z} \Dir_{\xi } G_{\xi } \Dir_{\xi }^{\ast } (z \hat{f}(\xi))  ). \label{higgscomm}
\end{align}
Since $D_{\xi }$ is a connection, the following commutation relations hold:
\begin{align*}
           & [D_{\xi }, z]  = \d z \land  & & [D_{\xi }, \bar{z} ]  = \d \bar{z} \land \\
           & [D_{\xi }^{\ast }, z]  = \frac{\partial }{\partial \bar{z}} \contr & &
           [D_{\xi }^{\ast }, \bar{z}] =\frac{\partial }{\partial {z}},
\end{align*}
where $\contr$ stands for contraction of a differential form by a vector field. 
It follows immediately 
\begin{align}
       [\Dir_{\xi }, z]  & = -[\Dir_{\xi }^{\ast }, z]  = \d z \land - \frac{\partial }{\partial \bar{z}} \contr = \d z \cdot \label{commrelz} \\
       [\Dir_{\xi }, \bar{z}] & = -[\Dir_{\xi }^{\ast }, \bar{z}]  =  \d \bar{z} \land -  \frac{\partial }{\partial {z}} \contr = \d
       \bar{z} \cdot \label{commrelzbar}
\end{align}
where the Clifford multiplication $\cdot $ is defined by these formulae. 
Plugging these in the expression (\ref{higgscomm}), using $\Dir_{\xi }^{\ast }\hat{f}(\xi;z)=0$ and 
$\hat{\pi }_{\xi } |_{ Im \Dir_{\xi }^{\ast }} =0$ together with the definition of $\hat{h}$, we get 
\begin{align} \label{commineq}
       \absl [\hat{\vartheta }, \hat{\vartheta }^{\ast } ] \hat{f}(\xi) \absr_{\hat{h}} & = \frac{1}{2}
       \left\Vert \hat{\pi }_{\xi }  \left(  \d z \cdot G_{\xi } \d \bar{z}\cdot  \hat{f}(\xi) -  
        \d \bar{z} \cdot  G_{\xi } \d z \cdot \hat{f}(\xi) \right) \right\Vert_{L^{2}(\C )} \notag \\ 
        & \leq  \frac{1}{2} \left\Vert  G_{\xi } \d \bar{z} \cdot  \hat{f}(\xi) \right\Vert_{L^{2}(\C )} + 
         \frac{1}{2} \left\Vert  G_{\xi }  \d z \cdot  \hat{f}(\xi) \right\Vert_{L^{2}(\C )}, 
\end{align}
since the norm of the orthogonal projection of a vector to a subspace is at most the norm of the vector and 
the action of Clifford multiplication by $\d z$ and $\d \bar{z}$ is point-wise bounded. 
We conclude by the first statement of Lemma \ref{estgreen}. 

Next, let us come to $\hat{\nabla }^{+} \hat{\vartheta}$. Similarly to the above, using (\ref{projform}) and the
commutation formulae (\ref{commrelz})-(\ref{commrelzbar}) we obtain 
\begin{align*}
      \left( \hat{\nabla }^{+} \hat{\vartheta} \right) \hat{f}(\xi) = & 
      \left( \hat{D}^{+} \circ \hat{\vartheta } - \hat{\vartheta } \circ \hat{D}^{+} \right) \hat{f}(\xi) \\
      = & \hat{\pi }_{\xi } \left( \hat{\d} - \frac{z}{2} \d \xi + \frac{\bar{z}}{2} \d \bar{\xi } \right) 
      \hat{\pi }_{\xi } \left( - \frac{z}{2} \right) \hat{f}(\xi) \\ 
      & - \hat{\pi }_{\xi } \left( - \frac{z}{2} \right) 
      \hat{\pi }_{\xi } \left( \hat{\d} - \frac{z}{2} \d \xi + \frac{\bar{z}}{2} \d \bar{\xi } \right)  \hat{f}(\xi)\\
      = & \hat{\pi }_{\xi } \Big{[} \left( \hat{\d} - \frac{z}{2} \d \xi + \frac{\bar{z}}{2} \d \bar{\xi } \right) 
       \Dir_{\xi } G_{\xi } \Dir_{\xi }^{\ast } \left( \frac{z}{2} \hat{f}(\xi) \right) \\
       & - \frac{z}{2} \Dir_{\xi } G_{\xi } \Dir_{\xi }^{\ast } 
       \left( \hat{\d} - \frac{z}{2} \d \xi + \frac{\bar{z}}{2} \d \bar{\xi } \right) \hat{f}(\xi) \Big{]} \\
      = & \hat{\pi }_{\xi } \Big{[} \left( \frac{1}{2} \d \xi \land \d z - \frac{1}{2} \d \bar{\xi } \land \d \bar{z } \right) \cdot 
      G_{\xi } \frac{\d z}{2} \cdot \hat{f}(\xi) \\
      & \hspace{.8cm} - \frac{\d z}{2} \cdot G_{\xi } \left( \frac{1}{2} \d \xi \land \d z - 
      \frac{1}{2} \d \bar{\xi } \land \d \bar{z } \right) \cdot \hat{f}(\xi)
      \Big{]} \\
       & +
      \hat{\pi }_{\xi } \left[  \hat{\d} \Dir_{\xi } G_{\xi } \frac{\d z}{2} \cdot \hat{f}(\xi) - \frac{\d z}{2} \cdot G_{\xi }
        \Dir_{\xi }^{\ast } \hat{\d} \hat{f}(\xi) \right]
\end{align*}
(here $\d z$ and $\d \bar{z}$ act on the spinors by Clifford multiplication, whereas $\d \xi $ and 
$\d\bar{\xi }$ by wedge product). Noticing that $| \d \xi |= | \d \bar{\xi }| =2$, the first term in the last
expression can be treated exactly as in (\ref{commineq}). For the second term, one only needs to remark 
that the commutation relations 
\begin{align*}
     \left[ \hat{\d} ,  D_{\xi } \right] = & \left[ \hat{\d} ,  D - \frac{\xi }{2} \d z + 
     \frac{\bar{\xi }}{2} \d \bar{z}  \right] \\
      = & -  \frac{\d \xi \land \d z \land }{2} + \frac{\d \bar{\xi } \land \d \bar{z} \land}{2} 
\end{align*}
and 
\begin{align*}
     \left[ \hat{\d} ,  D_{\xi }^{\ast } \right] = & -  \frac{\d \xi \land }{2} \frac{\partial }{\partial \bar{z}} \contr 
     + \frac{\d \bar{\xi } \land }{2} \frac{\partial }{\partial {z}} \contr 
\end{align*}
show that 
\begin{align*}
      \left[ \hat{\d} ,  \Dir_{\xi } \right] = - \left[ \hat{\d} ,  \Dir_{\xi }^{\ast } \right] = 
      - \frac{1}{2} \d \xi \land \d z \cdot + \frac{1}{2} \d \bar{\xi } \land \d \bar{z } \cdot 
\end{align*}
holds. Therefore we can proceed again as in (\ref{commineq}). 
\end{proof}

On the basis of Theorem \ref{biqjar}, the behavior of the transformed flat connection at infinity satisfies 
the hypothesis considered in \cite{Biq-Boa}. Namely, in a suitable gauge its difference from a model 
with second-order pole is in the weighted Sobolev-space $L^{1,2}_{-2+\delta}(\Omega^{1}\otimes E)$ considered in Section 2 
of that article. Indeed, passing to a coordinate $w=z^{-1}$, $|w|=\rho$ in which the double pole is in $0$, the norm of 
the perturbation is $O(\rho^{1+\delta})$, whereas that of its derivative is also $O(\rho^{1+\delta})$ (because the norm of $1$-forms
near infinity is $|\d z|=|\d w|/|w|=1$), and we conclude since $\rho^{1+\delta}/ \rho^{2} \in L^{2}_{\delta -2}$.  
It follows from the results of its Sections 7 and 8 that the analytic flat connection 
$\hat{D}$ defined outside infinity extends to an algebraic integrable connection with a parabolic structure 
on the singular fiber at infinity. On the other hand, such an extension over logarithmic singularities 
(that is, singularities in which the eigenvalues of $\hat{D}$ or equivalently those of $\hat{\vartheta}$ have at most 
first-order poles) is ensured by Theorem 2 of \cite{Sim}. Therefore, by Theorem \ref{eigenlog} the flat 
connection $\hat{D}$ on $\hat{\C} \setminus \hat{P}$ can be extended into a meromorphic integrable connection 
on $\CPt$ with parabolic structures at the singularities. 
\begin{defn} \label{transfext}
The \emph{transformed meromorphic integrable connection} is the meromorphic integrable connection with
parabolic structure in the singularities induced by the above extension procedures, subject to local changes 
of holomorphic trivialisations near the singularities to take all weights between $0$ and $1$. 
We will continue to denote it by $(\hat{E},\hat{D})$. The underlying extension will be 
called \emph{transformed extension} of the transformed bundle. 
\end{defn}
\begin{rk}
We will see in Section \ref{parweightsect} that the parabolic structures are adapted to the harmonic metric; 
namely, the weight $0\leq \hat{\alpha}_{k} <1$ of a subspace $F_{k}\hat{E}|_{p}$ of a singular fiber corresponds in 
local coordinate $z$ vanishing at the puncture to a decay bounded above by $|z|^{2\hat{\alpha}_{k}}$ of the norm 
of a parallel section extending an element of $F_{k}\hat{E}|_{p}$, as measured by the harmonic metric.
However, in Sections \ref{logext} and \ref{infext} we will construct a different 
extension over the punctures -- more suited to analytical study --, where the behavior of the norm of 
parallel sections near the singular points will no longer be bounded. We then pass back to the transformed 
extension in Corollary \ref{pardegcor}, where we remark that it is the one that establishes a "good" correspondence. 
\end{rk}

\chapter[Higgs bundle interpretation]{Interpretation from the point of view of Higgs bundles} \label{algint}
Let $(E,D,h)$ be a Hermitian bundle with integrable connection. 
Throughout this chapter, we suppose that the original metric $h$ is harmonic. This metric then defines a Higgs
bundle $(\E, \theta )$ starting from the integrable connection, via the procedure described in Section \ref{higgsintro}. 
We first prove that the transformed metric $\hat{h}$ is then harmonic for $\hat{D}$. 
Next, we give an interpretation of the transformed Higgs bundle of $(\E, \theta )$ in terms of the hypercohomology
of a sheaf map over $\CP$. These results 
will then be used to define the \emph{induced extension} $^{i}\hat{\E}$ of the 
transformed bundle over the punctures $\hat{P} \cup \{ \infty \}$, and to compute the topology and 
the singularity parameters of this extension of the transformed Higgs bundle. 
This will enable us to eventually compute the topology and 
the singularity parameters of the transformed Higgs bundle with respect to its transformed extension 
given in Definition \ref{transfext}. 

\section[Integrable connection and Higgs bundle]{The link with the transformed integrable connection }
Recall that we have defined the deformation of the Higgs bundle by the
formula (\ref{higgsdef}), and we write $D''_{\xi }$ for the $D''$-operator of
this deformation. Explicitly, we have 
\[
    D''_{\xi } = D''+ \theta_{\xi },
\]
where $\theta_{\xi }= \theta - \xi/2  \d z$. Moreover, as we have noticed in Section \ref{transf}, nonabelian Hodge theory
identifies the deformation of the Higgs bundle structure (\ref{higgsdef}) and that of the integrable
connection via the unitary gauge transformation 
\begin{equation*}
       g(z,\xi ) =e^{ [\bar{\xi} \bar{z}- \xi z]/2}. 
\end{equation*}
In other words, writing $g_{\xi }=g(.,\xi )$ for the gauge transformation restricted to
the fiber $\hat{H}_{\xi }$, we have 
\begin{equation} \label{gaugetr}
     g_{\xi }.D_{\xi }=D_{\xi }^{H}  = D - \frac{\xi}{2} \d z \land - \frac{\bar{\xi}}{2} \d \bar{z} \land.
\end{equation}
Since the gauge transformation $g_{\xi }$ is unitary, in addition to (\ref{gaugetr}) we have as well 
\begin{equation} \label{gaugetradj}
     g_{\xi }.D_{\xi }^{\ast }=(D_{\xi }^{H})^{\ast }.
\end{equation}
\begin{defn} The operator $\Dir_{\xi }^{H}= D_{\xi }^{H} - (D_{\xi }^{H})^{\ast }$ will be 
referred to as the \emph{Higgs Dirac operator}. In the same way, we let $\Dir_{\xi }''$ stand 
for the Dirac operator $D_{\xi }'' - (D_{\xi }'')^{\ast }$. 
The \emph{transformed smooth bundle underlying the Higgs bundle} 
is the bundle $\hat{V}$ over $\hat{\C} \setminus \hat{P}$ whose fiber over $\xi$ 
is the first $L^{2}$-cohomology space $L^{2}H^{1}(\EC_{\xi }^{H})$ of the operator $D_{\xi }^{H}$. 
\end{defn}
\begin{prop}\label{trvbprop}
This way we define a smooth vector bundle $\hat{V}$. Furthermore, there exists
a canonical bundle isomorphism between the smooth bundle $\hat{E}$ underlying the 
transformed integrable connection and the smooth bundle $\hat{V}$ underlying the transformed 
Higgs bundle. 
\end{prop}
\begin{proof}
Theorem \ref{coker} tells us that the transformed bundle underlying the integrable connection 
is the bundle of first $L^{2}$-cohomologies of $D_{\xi}^{int}$. 
For any $\xi$, the gauge transformation $g_{\xi }$ of $E$ induces a
natural isomorphism between the $L^{2}$-cohomology spaces of the complexes 
(\ref{ellcompl})  and 
\begin{equation}\label{unittrans}
        \Omega^{0}\otimes E \xrightarrow{\text{$g_{\xi }.D^{int}_{\xi}$}}  \Omega^1 \otimes E
        \xrightarrow{\text{$g_{\xi }.D^{int}_{\xi}$}} \Omega^2 \otimes E.
\end{equation} 
which is just $\EC_{\xi }^{H}$. In Theorem \ref{Fredholm} we have shown that the $0$-th and $2$-nd 
cohomology of $\EC_{\xi }$ vanishes for all $\xi\in \hat{\C} \setminus \hat{P}$, whereas Corollary \ref{cokerbundle}
implies that the cohomology spaces $L^{2}H^{1}(\EC_{\xi })$ define a smooth vector bundle over 
$\hat{\C} \setminus \hat{P}$. 
This then implies the same thing for $\EC_{\xi }^{H}$, whence the bundle isomorphism
between the bundles over $\hat{\C} \setminus \hat{P}$ in question. 
\end{proof}
Theorem \ref{laplker} has the following interpretation: 
\begin{thm} \label{laplkerhiggs} The first $L^{2}$-cohomology $\hat{V}_{\xi}=L^{2}H^{1}(\EC_{\xi}^{H})$ of the operator 
$D_{\xi }^{H}$ is canonically isomorphic to the kernel of the adjoint Dirac operator 
\begin{equation} \label{diradjhiggs}
     (\Dir_{\xi}^{H})^{\ast }: L^{2}(S^{-} \otimes E) \lra L^{2}(S^{+} \otimes E)
\end{equation}
on its domain, or alternatively to the kernel of the Laplace operator 
\begin{equation} \label{lapldirachiggs}
     \Delta_{\xi }^{H} = \Dir_{\xi }^{H} (\Dir_{\xi}^{H})^{\ast }  :
     L^{2}(S^{-} \otimes E) \lra L^{2}(S^{-} \otimes E)
\end{equation}
on its domain.  
\end{thm}
\begin{proof} 
Apply the gauge transformation $g$ to Theorem \ref{laplker} and notice that 
(\ref{gaugetr}) and (\ref{gaugetradj}) imply 
\begin{equation} \label{gaugetrdir}
    g_{\xi }.\Dir_{\xi }^{\ast } =(\Dir_{\xi }^{H} )^{\ast }
\end{equation}
and 
\begin{equation} \label{gaugetrlapl}
    g_{\xi }.\Delta_{\xi }=\Delta_{\xi }^{H}; 
\end{equation}
and in particular that 
\begin{equation} \label{gaugetrcoker}
     g_{\xi } (Ker (\Dir_{\xi }^{\ast })) = Ker ((\Dir_{\xi }^{H})^{\ast })
\end{equation}
and 
\begin{equation} \label{gaugetrker}
     g_{\xi } (Ker (\Delta_{\xi })) = Ker (\Delta_{\xi }^{H}).
\end{equation}
\end{proof}

This result enables us to put similar definitions as in the integrable deformation case. 
\begin{defn}
The hermitian bundle metric on $\hat{V}$ given by $L^{2}$ scalar product of the $(\Dir^{H}_{\xi})^{\ast}$-harmonic 
representative will be called the \emph{transformed hermitian metric}, and will be denoted by $\hat{h}$. 
Also, $\hat{\pi }_{\xi }^{H}$  will stand for $\hat{h}$-orthogonal projection of $L^{2}(S^{-} \otimes E)$ onto $\hat{V}$.
\end{defn}
\begin{rk}
Starting from a Higgs bundle with any Hermitian metric (not necessary harmonic), we can define in the same 
way its transform on the transformed bundle $\hat{V}$.
\end{rk}
Next, we recollect the above considerations in terms of the transformed bundles. 
\begin{prop} \label{gaugeprop} The family of gauge transformations $g$ induce a Hermitian bundle isomorphism between $\hat{E}$ 
and $\hat{V}$. 
Furthermore, the fiber $\hat{V}_{\xi }$ can be identified with  
the first $L^{2}$-cohomology of the single complex associated to 
the following double complex, denoted by $\mathcal{D}_{\xi }$:
\[
\xymatrix{
         \Omega^{0,1}\otimes E \ar[r]^{\theta_{\xi } \land} & \Omega^{2}\otimes E\\
        \Omega^{0}\otimes E \ar[r]^{\theta_{\xi } \land} \ar[u]^{\dbar^{\E}} & \Omega^{1,0}\otimes E. \ar[u]_{\dbar^{\E}} 
}
\]
\end{prop}
\begin{rk} Notice that commutativity of this diagram follows from the hypothesis
$\dbar^{\E}\theta =0$, which is just the definition of the harmonicity of $h$. 
\end{rk}
\begin{proof}
By (\ref{gaugetrker}),the $D_{\xi}$-harmonic representative of a class is mapped by $g$ into a 
$D_{\xi}^{H}$-harmonic class. Since the transformed metric from both points of view is 
induced by $L^{2}$-norm of the harmonic representatives, and $g$ is unitary, this gives 
the first statement. 
For the second, remark that by Theorem \ref{simpson}, the Laplace operator 
$\Delta_{\xi }^{H} $ is equal (up to a factor of 2)
to the Laplace operator $\Delta_{\xi }''= \Dir''_{\xi } (\Dir''_{\xi})^{\ast } $, therefore their
kernels coincide. This then identifies $\hat{V}$ with the first $L^{2}$-cohomology of the 
complex  
\begin{equation}\label{higgsellcompl}
     \Omega^{0}\otimes E \xrightarrow{\text{$D_{\xi}''$}}  \Omega^1 \otimes E
     \xrightarrow{\text{$D_{\xi}''$}} \Omega^2 \otimes E. 
\end{equation} 
Finally, recall that the formula 
\[ 
   D''_{\xi } = \dbar^{\E} + \theta_{\xi }
\]
gives the decomposition of $D''_{\xi } $ into its $(0,1)$- and $(1,0)$-part respectively. 
This means that the complex (\ref{higgsellcompl}) is the single complex associated to the double
complex  $\mathcal{D}_{\xi }$. However, it is not necessarily true that the domain of $D''_{\xi }$ is the sum 
of the domain of $\dbar^{\E}$ and that of $\theta_{\xi}$, it could in principle be larger. 
Still, the two $L^{2}$-cohomologies are the same. 
Indeed, suppose $f=f^{1,0} \d z + f^{0,1} \d \bar{z}\in L^{2}(\Omega^{1}\otimes E)$ is in the kernel of $D''_{\xi}$, that is 
\begin{equation}\label{dbarplustheta}
      \dbar^{\E}f^{1,0}\d z + \theta_{\xi }\land f^{0,1} \d \bar{z}=0.
\end{equation}
We wish to represent the $D''_{\xi}$-cohomology class of $f$ by a class 
$\tilde{f}^{1,0} \d z +\tilde{f}^{0,1}\d \bar{z}$ such that $\dbar \tilde{f}^{1,0}\in L^{2}$ and 
$\theta_{\xi}\tilde{f}^{0,1}\in L^{2}$. Away from logarithmic singularities, 
one can simply choose $f$ itself, for there locally ${f}^{0,1}\in L^{2}$
implies $\theta_{\xi}{f}^{0,1}\in L^{2}$ and by (\ref{dbarplustheta}) then $\dbar^{\E}f^{1,0}\in L^{2}$ as well. 
Thus we only need to modify $f$ in a neighborhood of the logarithmic  punctures. 
By Claim \ref{crsol} near any such puncture we can find $g\in L^{2}(E)$ such that $\theta_{\xi}g\in L^{2}(\Omega^{1,0}\otimes E)$ and 
$$
     f^{0,1} \d \bar{z} + \dbar^{\E}g=0.
$$
Using $\dbar^{\E}\theta_{\xi}=0$, the last two identities then also imply 
$$
     \dbar^{\E} (f^{1,0}\d z + \theta_{\xi }g)=0.
$$
Put $\tilde{f}^{1,0}\d z=f^{1,0}\d z + \theta_{\xi }g$;
as both $f^{1,0}$ and $\theta_{\xi }g$ are supposed to be in $L^{2}$, so is $\tilde{f}^{1,0}\d z$. 
This then shows that $f$ is cohomologuos in the $L^{2}$ complex of 
(\ref{higgsellcompl}) to a class locally represented by a section $\tilde{f}^{1,0} \d z$, 
where $\tilde{f}^{1,0}\in L^{2}$ and $\dbar^{\E}\tilde{f}^{1,0}\in L^{2}$. In different terms 
$\tilde{f}^{1,0}\d z\in\Dm (\dbar^{\E})$, and this shows that the first $L^{2}$-cohomology of 
(\ref{higgsellcompl}) is indeed equal to that of $\mathcal{D}_{\xi }$. 
\end{proof}

Next, let us investigate what the transformed integrable connection $\hat{D} $ and
its unitary part $\hat{D}^{+}$ become under this gauge transformation. 
Notice that since the gauge transformation $g$ is
unitary, the orthogonal projection $\hat{\pi}$ onto $\hat{E}$ is transformed into the
orthogonal projection $\hat{\pi}^{H}$ onto $\hat{V}$, with respect to
the same $L^{2}$-metric on the fibers; in different terms $g_{\xi }.\hat{\pi}_{\xi}=\hat{\pi}_{\xi}^{H}$. 
The image of the transformed integrable connection $\hat{D}$ under the
gauge transformation $g$ in the point $\xi $ is given by 
\begin{align} \label{htrint}
     \hat{D}^{H}& = g.\hat{D}\notag \\
     & = g.(\hat{\pi}_{\xi}\circ (\hat{\d}-z \d \xi \land ))\\
     & = \hat{\pi}_{\xi}^{H} \left( \hat{\d} - \frac{1}{2}(z\d \xi \land + \bar{z} \d\bar{\xi } \land ) \right),\notag 
\end{align}
(see (\ref{trunit}), (\ref{trhiggs}) and Proposition \ref{explflattr}), 
and that of the candidate Higgs field is the endomorphism  
\begin{align} \label{htrhiggs}
    \hat{\theta }^{H} & = g.\hat{\theta }\notag \\
    &= g.(\hat{\pi}_{\xi}\circ (-z/2 \d \xi \land ) ) \\
    &= - \frac{1}{2} \hat{\pi}_{\xi}^{H} (z \d \xi \land). \notag
\end{align}
Therefore, if we decompose the transformed flat connection in the point of view of Higgs bundles into its 
unitary and self-adjoint part, we obtain 
\begin{align} \label{hdecomp}
     (\hat{D}^{H})^{+}=  \hat{\pi}_{\xi}^{H} ( \hat{\d}) &&
     (\hat{D}^{H})^{sa}= \hat{\theta }^{H}+(\hat{\theta }^{H})^{\ast}
\end{align}
(these formulae can also be deduced directly from Proposition \ref{unittr}). This then gives the 
desired interpretation of the transformed unitary connection $\hat{D}^{+}$ in this point of view. 
\begin{defn} \label{defholtrafo}
We let $\dbar^{\hat{\E}}$ stand for the $(0,1)$-part of $(\hat{D}^{H})^{+}$. 
Moreover, we call the holomorphic bundle $\hat{V}$ with partial connection $\dbar^{\hat{\E}}$ 
the \emph{transformed holomorphic bundle} and we denote it by $\hat{\E}$. 
\end{defn}

\section{Harmonicity of the transformed metric} \label{sectharm}
In this section we prove the following result:
\begin{thm} 
If the original metric $h$ is harmonic, then the same thing is true for $\hat{h}$. 
\end{thm}
\begin{proof}
First remark that by (\ref{hdecomp}), the formula for $\dbar^{\hat{\E}}$ is $\hat{\pi}_{\xi}^{H}(\hat{\d}^{0,1})$. 
Also, the $(1,0)$-part of $(\hat{D}^{H})^{sa}$ is just $\hat{\theta }^{H}$.
By definition, harmonicity of the transformed metric $\hat{h}$ resumes then in the
equation 
\begin{equation} \label{hharm}
    \dbar^{\hat{\E}} \hat{\theta }^{H}=0.
\end{equation}
By Proposition \ref{gaugeprop} we have $\hat{V}_{\xi } = L^{2}H^{1}(D_{\xi }'')$,
with  $D''_{\xi } = D''- \xi/2 \d z $. From this formula it is clear
that $D_{\xi }''$ depends holomorphically on $\xi $, so we are in the situation
described in part 3.1.3 of \cite{Don-Kronh} of chain complexes 
$$
     \Omega^{0}\otimes E \xrightarrow{\text{$D_{\xi}''$}}  \Omega^1 \otimes E
     \xrightarrow{\text{$D_{\xi}''$}} \Omega^2 \otimes E 
$$
varying holomorphically with $\xi $. There it is shown that if the 
first cohomology spaces $\hat{V}_{\xi }$ of these complexes are all finite 
dimensional, of the same dimension, then the bundle $\hat{V}$ constructed out of them
over the parameter space of $\xi$ carries a natural holomorphic structure. 
Explicitly, this is given by by saying that a section $f \in \Gamma(\hat{V})$ in a
neighborhood of $\xi_{0}$ is holomorphic if and only if it admits a lift 
$\tilde{f} \in \Gamma(Ker(D_{\xi}''|_{\Omega^{1}}))$ which is itself holomorphic with respect to
the holomorphic structure induced by the $(0,1)$-part $\hat{\d }^{0,1}$ of
the trivial connection $\hat{\d }$ on the Hilbert bundle $\hat{H}$. 
This holomorphic structure is the same as the
one defined by the operator $ \dbar^{\hat{\E}}$, since both are induced 
by $\hat{\d}^{0,1}$ and $ \hat{\pi}^{H } $. The section $\hat{\theta }^{H} \in End( \hat{V}) \otimes \Omega_{\hat{\C}}^{1,0}$ 
is then holomorphic for this holomorphic structure if
and only if it maps each holomorphic section $f$ into a holomorphic
section. In particular, this is the case if it admits a lift 
\[ \xymatrix{ Ker(D_{\xi}''|_{\Omega^{1}}) \ar[r]^(0.4){\Theta} & Ker(D_{\xi}''|_{\Omega^{1}}) \otimes
  \Omega_{\hat{\C}}^{1,0} \\
  \hat{V}_{\xi } \ar[u] \ar[r]^(0.4){\hat{\theta }^{H}} & \hat{V}_{\xi } \otimes
  \Omega_{\hat{\C}}^{1,0},  \ar[u] 
   }
\]   
such that 
\begin{enumerate}
\item{$\Theta$ passes to the quotient $Ker(D_{\xi}''|_{\Omega^{1}}) \ra
    Ker(D_{\xi}''|_{\Omega^{1}})/Im(D_{\xi}''|_{\Omega^{0}})=\hat{V}_{\xi }$, the quotient
    being $\hat{\theta }^{H}$, and  } \label{hol1}
\item{$\Theta$ is holomorphic with respect to the holomorphic structure induced
  by $\hat{\d}^{0,1} $.} \label{hol2}
\end{enumerate}
Recall from Section \ref{l2coh} that $Ker(D_{\xi}''|_{\Omega^{1}})$ is a closed Hilbert subspace of
$\hat{H}_{\xi }$; call $\pi_{Ker(D_{\xi}'')}$ orthogonal projection of $\hat{H}_{\xi }$ to it. We now claim that the map 
\begin{align*}
     \Theta  : Ker(D_{\xi}''|_{\Omega^{1}}) & \lra Ker(D_{\xi}''|_{\Omega^{1}}) \otimes
     \Omega_{\hat{\C}}^{1,0} \\ 
     \tilde{f}_{\xi } & \mapsto -\frac{1}{2} \pi_{Ker(D_{\xi}'')}(z\tilde{f}_{\xi }(z)) \d \xi 
\end{align*}
verifies the hypotheses needed. 

For (\ref{hol1}), we need to show $\Theta ( Im ( D_{\xi}''|_{\Omega^{0}} ) ) \subseteq Im (D_{\xi}''|_{\Omega^{0}} ) $. 
Let $g_{\xi }$ be a local section of the trivial Hilbert bundle $ L^{2}(E) \ra \hat{\C} $. Then we have 
\begin{align*}
      \Theta  (D_{\xi}''g) &= -\frac{1}{2}\pi_{Ker(D_{\xi}'')}(z D_{\xi}''g_{\xi }) \d \xi \\
      & = -\frac{1}{2} \pi_{Ker(D_{\xi}'')}( D_{\xi}''(zg_{\xi }(z))) \d \xi \\
      & = -\frac{1}{2} D_{\xi}''(zg_{\xi }(z)) \d \xi, 
\end{align*}
because the operator $D_{\xi}''=\dbar^{\E}+\theta_{\xi }$ commutes with multiplication by
$z$, and $Im( D_{\xi}''|_{\Omega^{0}} ) \subseteq Ker ( D_{\xi}''|_{\Omega^{1}} )$. This shows
that $Im( D_{\xi}''|_{\Omega^{0}} )$ is invariant by $\Theta  $; the quotient is
clearly $\hat{\theta }^{H}$.

Next come to (\ref{hol2}): we remark that the formula defining $\Theta $ only
depends on $\xi$ via the projection $\pi_{Ker(D_{\xi}'')}$. But since the
operator $D_{\xi}''$ depends holomorphically in $\xi $, so do the subspaces 
$Ker(D_{\xi}'')$, and since the metric is independent of $\xi $, the same thing
is true for the projections $\pi_{Ker(D_{\xi}'')}$. This shows that $\Theta $, and so
$\hat{\theta }^{H}$ is holomorphic in $\xi $. 
\end{proof}

\section{Identification with hypercohomology}\label{ident}
In this section we will often use basic properties of hypercohomology; for an introduction to this topic, 
we refer to \cite{Griff-Harr} and \cite{Dem}.

Before we start, we need to define the functional spaces 
\begin{align*}
     \tilde{L}_{\xi }^{2}(E) & = \Dm(D_{\xi }''|_{ \Omega^{0}\otimes E }) \\
     & = \{ u \in  L^{2}( E ) \ : \  \theta_{\xi } \land  u, \dbar^{\E} u \in L^{2} \}\\
     \tilde{L}_{\xi }^{2}(\Omega^{0,1} \otimes E) & = \Dm(D_{\xi }''|_{ \Omega^{0,1}\otimes E }) \\ 
     & = \{ v \d \bar{z} \in L^{2}( \Omega^{0,1}\otimes E ) \ : \ \theta_{\xi } \land  v \d \bar{z} \in L^{2} \} \\
     \tilde{L}^{2}(\Omega^{1,0} \otimes E) & = \Dm(D_{\xi }''|_{ \Omega^{1,0}\otimes E }) \\
     & = \{ u \d z \in  L^{2}( \Omega^{1,0}\otimes E ) \ : \ \dbar^{\E} (u \d z) \in L^{2} \},  
\end{align*}
for the Euclidean metric $|\d z|^{2}$ on $\C$ and the hermitian metric $h$ on the fibers, 
adapted to the parabolic structure with weights $\{\alpha_{1}^{j},\ldots,\alpha_{r}^{j} \}$. 
Notice that we may drop the index $\xi $ of these spaces, since they all coincide: indeed, in a
logarithmic singularity the deformation $ \xi \d z $ is bounded, and at infinity the
condition $\xi \notin \hat{P} $ implies that no eigenvalues of $\theta_{\xi }$ vanish,
and this gives equivalence of the corresponding norms exactly as in Lemma
\ref{lemmaone}. We identify these functional spaces to the sheaves of their local
sections. 
In what follows, we are going to define sheaves $\E$ and $\F$ of sections of 
$\Omega^{0} \otimes E $ and $\Omega^{1,0} \otimes E $ respectively on $\CP$ with the property that 
the $L^{2}$-cohomology $L^{2}H^{\bullet }(D_{\xi }'')$ of (\ref{higgsellcompl}) 
identifies to the hypercohomology  $\H^{\bullet }(\E \xrightarrow{\theta_{\xi } \land } \F)$ of
the sheaf map $\E \xrightarrow{\theta_{\xi } \land } \F$. This latter is then
explicitly given in terms of a sky-scraper sheaf over the zero set $\Sigma_{\xi }$
of $det(\theta_{\xi })$ by a simple use of the spectral sequence of the double
complex.

\subsection{Definition and resolution of the sheaves} \label{sheaves}

Recall that the parabolic structure on $\E$ with adapted Hermitian fiber metric means that the holomorphic
bundle $\E$ on $\C \setminus P$ has a natural extension to all $\CP$:
the holomorphic sections at a singular point are the holomorphic sections
outside the singularity which are bounded with respect to the metric. 
By an abuse of language, for $U \subset \CP $ an open set let $\E|_{U}$ be the set of
holomorphic sections of the bundle $\E$ in U. In other words, we denote
by $\E$ the sheaf of local holomorphic sections of $\E$ (extended over the punctures as above).

Next, let us define $\F$: for an open set $U \subset \CP$ containing no singular point,
let $\F|_{U}$ be the set of $\dbar^{\E}$-holomorphic sections of
$\Omega^{1,0} \otimes E$. If $U$ contains exactly one singular point $p_{j} \in P $ (and does
not contain the infinity), then let 
$\F|_{U}$ be the set of $\dbar^{\E}$-meromorphic sections $\sigma \d z$ of
$\Omega^{1,0} \otimes E$ such that $\sigma$ be $\dbar^{\E}$-meromorphic in $U$ with only one simple
pole at $p_{j}$, and such that its residue in this point be contained in the subspace
$Im(Res(\theta,p_{j}))$. Finally, if $U$ contains the infinity (but no other singular
points), then let $\F|_{U}$ be 
the set of all $\dbar^{\E}$-meromorphic sections $\sigma \d z$ of $\Omega^{1,0} \otimes E$ with a
double pole at infinity, and no other poles in $U$. Notice that since in
the coordinate $w=1/z$ of $\CP$ the section $\d z$ has a double pole at infinity, 
this amounts to say that $\sigma $ is a $\dbar^{\E}$-holomorphic section of $E$ in
$U$. Writing $\sigma = \sum_{k} f_{k}^{\infty } \sigma_{k}^{\infty }$ in the holomorphic basis
(\ref{holtrivinf}) at infinity, it is still the same thing to say that $f_{k}^{\infty }$ be a
holomorphic function in $U$ for all $k$ (in particular bounded at infinity). 
It is easy to check that this way we defined a sheaf.  

We introduce some further notation: set $\tilde{r} = \sqrt{1+|z|^{2}} $ on $\C$; 
then for $a \in \{ 0,1 \}$ we denote by 
$\tilde{r} \tilde{L}^{2}(\Omega^{a,0}\otimes E)$ the space of sections $u$ of $\Omega^{a,0}\otimes E$ such
that $\tilde{r}^{-1}u \in \tilde{L}^{2}$. This way we only loosen the condition on
the behavior of $u$ at infinity with respect to $\tilde{L}^{2}$, namely that $r^{-1}u$ be in
$\tilde{L}^{2}$ in a neighborhood of infinity. It is immediate that there exist an inclusion 
of vector spaces 
\begin{equation}\label{inclvsp}
       \tilde{L}^{2}(\Omega^{a,0}\otimes E) \hookrightarrow \tilde{r} \tilde{L}^{2}(\Omega^{a,0}\otimes E).
\end{equation}
\begin{lem} \label{lemrese}
The sequence 
\begin{eqnarray} \label{resole}  
              \E \inc \tilde{r} \tilde{L}^{2}(E) \xrightarrow{\dbar^{\E}} \tilde{L}^{2}(\Omega^{0,1} \otimes E)
\end{eqnarray}
is a resolution of $\E$. 
\end{lem}
\begin{proof}
It is known that away from the singularities, the sequence of usual $L^{2}$-sections with respect to
Euclidean metric gives a resolution of the sheaf of holomorphic sections. 
Therefore, we only need to show that  (\ref{resole}) is a resolution at the
singularities. 

Consider first $p_{j} \in P$. We first prove that (\ref{resole}) is 
locally exact in $\tilde{r} \tilde{L}^{2}(E)$.  
Let $E$ be trivialized in $ \Delta(p_{j},\varepsilon ) $ by the local
sections $\{\sigma_{k}^{j} \}$ given in (\ref{holtriv}). As we have seen in 
(\ref{modelhiggs}), in this trivialisation up to a perturbation term 
$\theta = diag (\lambda_{k}^{j}) {\d z }/{z}$, with
$\lambda^{j}_{k}=(\mu^{j}_{k}-\beta^{j}_{k})/2$, and the parabolic weights are given by
$\alpha^{j}_{k}=\Re (\mu^{j}_{k}) -[\Re (\mu^{j}_{k})]$. By definition,
any holomorphic section $\sigma $ of $E^{j}$ can be given as a sum $\sum_{k}\phi^{j}_{k}
\sigma_{k}^{j}$, where $\phi^{j}_{k}$ are holomorphic functions defined in 
$\Delta(p_{j},\varepsilon )$, in particular bounded by a constant $K$. 
This implies that $\sigma \in L^{2}(E)$, so that 
$\sigma \in \tilde{L}^{2}(E)$ if and only if $\theta \land \sigma \in L^{2}$. 
Recall that $L^{2}$ is defined with respect
to the parabolic structure $\{ \alpha_{k}^{j} \}$, and that the perturbation
term in $\theta$ behaves as $O(r^{-1+\delta })$ with $\delta >0$, where $r=|z-p_{j}|$. This implies that 
\begin{align*}
     \int_{\Delta(p_{j},\varepsilon )} |\theta \sigma |^{2} & \leq K' \int \sum_{k=1}^{r_{j}}
     |r^{-1+\delta}\sigma_{k}^{j}|^{2 } 
     + K' \int  \sum_{k=r_{j}+1}^{r} |r^{-1}\sigma_{k}^{j}|^{2} \\
     & \leq K''  \int \sum_{k=1}^{r_{j}} |r^{-1+\delta}|^{2 } +  K'' \int  \sum_{k=r_{j}+1}^{r}
     |r^{-1+\alpha_{k}^{j}}|^{2}. 
\end{align*}
By Hypothesis \ref{main}, $\alpha_{k}^{j} > 0$ for all $j \in \{ r_{j}+1,\ldots,r \}$. 
It then follows that this last expression is finite, which proves that 
any holomorphic section of $E$ is in $\tilde{L}^{2}$. On the other hand, 
if a section $\sigma = \sum_{k}\phi_{k}^{j}\sigma_{k}^{j}$ of $E$ is meromorphic in
$p_{j}$, then there is at least one $k \in \{ 1,\ldots,r\}$ such that $\phi_{k}^{j}$
has a pole in $p_{j}$. Suppose $k \in \{ 1,\ldots,r_{j}\}$: then
$|\phi_{k}^{j}\sigma_{k}^{j}|\sim 1/r$, and $\sigma $ is clearly not in $L^{2}$. Suppose
now $k \in \{ r_{j}+1,\ldots,r \}$: then again by Hypothesis \ref{main} we have 
$\lambda_{k}^{j} \neq 0$, and therefore $|\theta \land \phi_{k}^{j}\sigma_{k}^{j}| \sim r^{-2+\delta }$, and
so $\theta \land \sigma \notin L^{2}$.  Hence, the sections of $\tilde{L}^{2}(\Delta(p_{j},\varepsilon ),E)$ in the kernel
of $\dbar^{\E}$ are exactly the local holomorphic sections of $E$, in other
words the local sections of $\E$. This shows local exactness in $\tilde{L}^{2}(E)$. 

The next thing we show is that in $\Delta(p_{j},\varepsilon )$ the complex
(\ref{resole}) is exact at $\tilde{L}^{2}(\Omega^{0,1} \otimes E)$: 
let $v \d \bar{z} \in \tilde{L}^{2}(\Delta(p_{j},\varepsilon ), \Omega^{0,1} \otimes E)$ be
an arbitrary section; for $\varepsilon >0$ sufficiently small we wish to find 
$u \in \tilde{L}^{2}(\Delta(p_{j},\varepsilon ),E)$ such that 
\begin{equation} \label{dbareq}
            \dbar^{\E}u=v \d\bar{z}
\end{equation}   
We can suppose without restricting generality that $v= f \sigma_{k}^{j}$, with 
$f $ a function defined in $\Delta(p_{j},\varepsilon )$. Since $\sigma_{k}^{j}$ is a
holomorphic section of $E$, solving (\ref{dbareq}) boils down to solving
the usual Cauchy-Riemann equation on the disk 
\begin{equation} \label{dbareq2}
           \frac{\partial g}{\partial \bar{z} } = f
\end{equation} 
with $u= g \sigma_{k}^{j} \in \tilde{L}^{2}(\Delta(p_{j},\varepsilon ), E)$. Exactness
near a singularity at a finite point is given by the following claim: 
\begin{clm} \label{crsol}
For $f \in L^{2}$ the  equation (\ref{dbareq2}) has a
solution $g$ such that $g r^{-1+\delta } \in L^{2}$ for any $\delta >0$. 
For $f$ such that $fr^{\alpha } \in L^{2}$ with $0<\alpha <1$, (\ref{dbareq2}) has a solution $g$ 
such that $g r^{-1+\alpha  } \in L^{2}$. 
\end{clm}
\begin{proof} 
The first statement is established combining the usual resolution of the
Cauchy-Riemann equation for $f \in L^{2}$ by an $L^{2,1}$-function $g$ and the
estimation (\ref{regineq}). 

The second one is a direct consequence of Proposition I.3 of \cite{Biq91}. One
might also prove it by direct estimations on the solution given by the
Cauchy kernel, as in \cite{At-Pat-Sin}.  
\end{proof}

Now let us come back to exactness at a singularity in a finite point: 
for the regular case $k\in \{ 1,\ldots,r_{j} \}$ we have $f\in L^{2}$ and  
$|\theta\land g \sigma_{k}^{j}| \leq |g|r^{-1+\delta }$, so we can apply directly the first
statement of the claim; 
for the singular case $k\in \{ r_{j}+1,\ldots,r \}$ by definition $|\theta\land f
\sigma_{k}^{j} \d \bar{z}| \sim  |f|r^{-1+\alpha }$ is in $L^{2}$ with $\alpha >0$ by
Hypothesis \ref{main}, therefore we can apply the second statement of
the claim. Remark that in this case even a stronger condition then the
assumption $fr^{\alpha } \in L^{2}$ of the claim holds. However, we will need the
claim in its full generality to show exactness at infinity.  

We now come to exactness at infinity. Recall that $\xi \notin \hat{P}$ implies
$\theta_{\xi } $ is an isomorphism
$L^{2}(\Omega^{0,b}) \ra L^{2}(\Omega^{1,b})$ for $b \in \{ 0,1 \}$. Therefore, the sections at
infinity of the sheaves $\tilde{L}^{2}(\Omega^{0,b})$ and $L^{2}(\Omega^{0,b})$
coincide. First, we consider exactness in
$\tilde{r} \tilde{L}^{2}(E)= \tilde{r} {L}^{2}(E)$: by the definition of $\E$, its local sections are the 
holomorphic linear combinations $\sigma = \sum_{k} \phi_{k}^{\infty }\sigma_{k}^{\infty } $. First we
check that these sections verify $r^{-1}\sigma \in L^{2}$: since $|\phi_{k}^{\infty }| \leq K$
and $|\sigma_{k}^{\infty }| \sim r^{-\alpha_{k}^{\infty}}$ with $\alpha_{k}^{\infty} >0$ by Hypothesis \ref{main}, 
we see that $r^{-1}\phi_{k}^{\infty }\sigma_{k}^{\infty } \in L^{2}$. On the other hand, if we have a
section $\sigma = \sum_{k} \phi_{k}^{\infty }\sigma_{k}^{\infty } $ in the kernel of
$\dbar^{\E}$, then for all $k$ the function $\phi_{k}^{\infty }$ is either
holomorphic or meromorphic; but if $r^{-1}\sigma \in L^{2}$, then it implies that 
$\phi_{k}^{\infty }$ is holomorphic for all $k$. This proves exactness in the first term. 

Next we come to the term $L^{2}(\Omega^{0,1} \otimes E )$: for a section $v \d \bar{z} \in
L^{2}(\C \setminus \Delta(R), \Omega^{0,1}\otimes E)$ we search $u \in rL^{2}(\C \setminus
\Delta(R), E)$  such that $ \dbar^{\E}u=v$. Suppose $v=f \sigma_{k}^{\infty }$ and 
$u=g \sigma_{k}^{\infty }$ again. In the coordinate $w=1/z=\rho e^{-\theta }$ on $\Delta(0,1/R)$
we find (for simplicity we took $R=1$ and wrote $\Delta= \Delta(0,1/R)$ ): 
\begin{align*}
     \int_{\Delta} | f |^{2} \rho^{2\alpha -4} |\d w |^{2} & = \int_{\C \setminus \Delta} | f |^{2} 
            r^{-2\alpha } |\d z |^{2} < \infty \\
     \int_{\Delta} | g |^{2} \rho^{2\alpha -2} |\d w |^{2} & = \int_{\C \setminus \Delta} | g |^{2} 
            r^{-2-2\alpha } |\d z |^{2} < \infty . 
\end{align*}
On the other hand, the Cauchy-Riemann equation 
$$
    \frac{\partial g}{\partial \bar{z} }=f
$$
transforms into 
$$
    \frac{\partial g}{\partial \bar{w} }= -\frac{f}{\bar{w}^{2}}. 
$$
and we conclude applying Claim \ref{crsol} to $-f/\bar{w}^{2}$. 

\end{proof}

We can also show the counterpart of Lemma \ref{lemrese} for $\F$: 
\begin{lem} \label{lemresf} The complex
\begin{equation} \label{resolf}
    \F \inc \tilde{r} \tilde{L}^{2}(\Omega^{1,0} \otimes E) \xrightarrow{\dbar^{\E}} L^{2}(\Omega^{1,1} \otimes E)
\end{equation}
is a resolution of $\F$.
\end{lem}
\begin{proof}
Away from the singularities this is also given by classical elliptic
theory, therefore we focus our attention on a neighborhood of a
singular point. 

Let us first treat the case of a singularity at a finite point $p_{j} \in P$. A
local section of $\F$ is then by definition a section
$\sigma = \sum_{k}\phi_{k}^{j}\sigma_{k}^{j} \d z$ such that $\phi_{k}^{j}$ is holomorphic for 
$k \in \{ 1,\ldots r_{j} \}$ and has a pole of order at most one in $p_{j}$ for 
$k \in\{r_{j}+1,\ldots r\}$. 
From the form of the parabolic structure, it follows that
$|\phi_{k}^{j}\sigma_{k}^{j}| \sim O(1)$ for $k \in \{ 1,\ldots r_{j} \}$ and 
$|\phi_{k}^{j}\sigma_{k}^{j}| \sim O(r^{-1+\alpha_{k}^{j}})$ for  $k \in\{r_{j}+1,\ldots r\}$. 
By Hypothesis \ref{main} we have $\alpha_{k}^{j}>0$, thus $\sigma \in L^{2}(\Omega^{1,0}\otimes
E)$. On the other hand, if a section $\sigma = \sum_{k}\phi_{k}^{j}\sigma_{k}^{j} \d z$ 
of $\Omega^{1,0}\otimes E$ satisfies $\dbar^{\E}\sigma=0 $, but $ \sigma \notin L^{2}(\Omega^{1,0}\otimes E)$ then
either $\phi_{k}^{j}$  has a pole for some $k \in \{ 1,\ldots r_{j} \}$ or $\phi_{k}^{j}$
has an at least double pole for some $k \in\{r_{j}+1,\ldots r\}$, and therefore $\sigma $ is not a
local section of $\F$. This shows exactness in the first term. 

Consider now exactness at the second term in $\Delta(p_{j},\varepsilon )$: here 
we need to solve  (\ref{dbareq2}), for $f \in L^{2}$ with the solution $g$ in $L^{2}$
in the regular case; and for $f$ such that $fr^{\alpha } \in L^{2}$ with the solution $g$ such
that $gr^{\alpha } \in L^{2}$ in the singular case. Both follow from Claim \ref{crsol}.

There now remains to show exactness at infinity: this is done similarly to
the case of $\E$. 
\end{proof}

\subsection{Hypercohomology and $L^{2}$-cohomology}

We can use the results of the last section in order to deduce the following: 
\begin{prop}\label{algl2}
The first $L^{2}$-cohomology $\hat{V}_{\xi } = L^{2}H^{1 }(D_{\xi }'')$ of 
(\ref{higgsellcompl}) is isomorphic to the hypercohomology  $\H^{1 }(\E \xrightarrow{\theta_{\xi } \land }\F)$. 
\end{prop}
\begin{proof} By Lemmas \ref{lemrese} and \ref{lemresf}, $\theta_{\xi }$ defines a
morphism of resolutions 
\begin{equation} \xymatrix{ \tilde{L}^{2}(\Omega^{0,1}\otimes E) \ar[r]^{\theta_{\xi } \land }    &  {L}^{2}(\Omega^{1,1}\otimes E) \\
              \tilde{r} \tilde{L}^{2}(E) \ar[r]^(.4){\theta_{\xi } \land} \ar[u]^{\bar{\partial}^{\E} } &
              \tilde{r} \tilde{L}^{2}(\Omega^{1,0}\otimes E)  \ar[u]_{\bar{\partial}^{\E} } \\
              \mathcal{E} \ar[r]^{\theta_{\xi } \land} \ar@{^{(}->}[u] & \mathcal{F} \ar@{^{(}->}[u] .
} \label{mr}
\end{equation}
Therefore, by general theory, the hypercohomology of the sheaf map $\E
\xrightarrow{\theta_{\xi } \land }\F$ identifies to the cohomology of the single
complex formed by the double complex $\mathcal{D}^{r}_{\xi }$:  
\begin{eqnarray} \label{drcomplex} 
              \xymatrix{ \tilde{L}^{2}(\Omega^{0,1}\otimes E) \ar[r]^{\theta_{\xi } \land}    &  {L}^{2}(\Omega^{1,1}\otimes E) \\
              \tilde{r} \tilde{L}^{2}(E) \ar[r]^(.4){\theta_{\xi } \land} \ar[u]^{\bar{\partial}^{\E} } &
              \tilde{r} \tilde{L}^{2}(\Omega^{1,0}\otimes E)  \ar[u]_{\bar{\partial}^{\E} }. 
            }
\end{eqnarray}
We show that the first cohomology of the single complex of this double
complex is isomorphic to the first cohomology of the single complex associated
to the double complex $\mathcal{D}_{\xi }$: 
\begin{eqnarray} \label{dcomplex} 
              \xymatrix{ \tilde{L}^{2}(\Omega^{0,1}\otimes E) \ar[r]^{\theta_{\xi } \land}    &  {L}^{2}(\Omega^{1,1}\otimes E) \\
              \tilde{L}^{2}(E) \ar[r]^(.4){\theta_{\xi } \land} \ar[u]^{\bar{\partial}^{\E} } &
              \tilde{L}^{2}(\Omega^{1,0}\otimes E)  \ar[u]_{\bar{\partial}^{\E} }. 
            }
\end{eqnarray}
We define a map 
$$
     \iota : H^{1}(\mathcal{D}_{\xi }) \lra  H^{1}(\mathcal{D}^{r}_{\xi })
$$
as follows: represent a cohomology class of $H^{1}(\mathcal{D}_{\xi })$ by a
couple 
$$
    (\kappa\d \bar{z},\nu \d z )\in \tilde{L}^{2}(\Omega^{0,1}\otimes E) \oplus \tilde{L}^{2}(\Omega^{1,0}\otimes E),
$$ 
and use the inclusion (\ref{inclvsp}) to map it into the cohomology class represented by the same couple 
$(\kappa,\nu )$ in $H^{1}(\mathcal{D}^{r}_{\xi })$. This is well defined, since if 
$(\kappa \d \bar{z}+\dbar^{\E}\lambda ,\nu \d z + \theta_{\xi }\lambda )$ 
is a couple in $H^{1}(\mathcal{D}_{\xi })$ representing the same class
as $(\kappa \d \bar{z},\nu \d z)$, for $\lambda \in \tilde{L}^{2}(E)$, then in particular 
$\lambda \in\tilde{r} \tilde{L}^{2}(E)$, and so the two couples are cohomologuous in
$H^{1}(\mathcal{D}^{r}_{\xi })$ as well. This also shows that $\iota $ is injective.

We only need to prove surjectivity: suppose we have a couple 
$(\kappa\d \bar{z},\nu \d z)\in \tilde{L}^{2}(\Omega^{0,1}\otimes E) \oplus \tilde{r} \tilde{L}^{2}(\Omega^{1,0}\otimes E)$
representing a class in $H^{1}(\mathcal{D}^{r}_{\xi })$. It is clearly sufficient
to prove that this class can be represented by a couple vanishing in a
neighborhood of infinity. Since $\theta_{\xi }$ is an isomorphism at infinity, we
can put (restricting to a smaller neighborhood of infinity if necessary) 
$\lambda = \theta_{\xi }^{-1}(\nu \d z)$. This is then a section in $\tilde{r} \tilde{L}^{2}(E)$,
and the couple $(\kappa\d \bar{z} - \dbar^{\E} \lambda ,\nu \d z - \theta_{\xi }\lambda) $ is
cohomologuous to $(\kappa\d \bar{z},\nu \d z)$ in $H^{1}(\mathcal{D}^{r}_{\xi })$. By
definition, the $(1,0)$-term of this couple vanishes at infinity. The same
thing is true for the $(0,1)$-part, because $\theta_{\xi } (\kappa\d \bar{z} -
\dbar^{\E} \lambda) = -\dbar^{\E}(\nu \d z - \theta_{\xi }\lambda) =0$ near infinity and $\theta_{\xi }$ is an
isomorphism there. This finishes the proof of the proposition, for the
$L^{2}$-cohomology of (\ref{higgsellcompl}) is by definition the cohomology
of the single complex associated to $\mathcal{D}_{\xi }$.  
\end{proof}

\subsection{The spectral curve} \label{specsect}
In the explicit identification of the hypercohomology, the following notions will be of much importance. 
Recall that (up to wedge product by $\d z$) $\theta_{\xi }$ is a meromorphic section of $End(E)$ over $\CP$. 
\begin{defn} \label{specsetdef}
For $\xi \in \hat{\C} \setminus \hat{P}$, the set of zeros of $det{(\theta_{\xi })}$ is called the \emph{spectral set}
corresponding to $\xi$. We denote it by $\Sigma_{\xi }$. 
\end{defn}
\begin{lem} \label{efdiv}
For each $\xi \in \hat{\C} \setminus \hat{P}$, the spectral set  is an effective divisor of $\CP$, 
in other words a finite set of points with multiplicities in $\N$. 
\end{lem}
\begin{proof}
The section $\det{(\theta_{\xi })}$ of $End(V)$ is holomorphic with respect to $\dbar^{\E}$. We only 
need to check it does not vanish identically for any $\xi$. Suppose there exists $\xi$ such that 
$$
     det{(\theta_{\xi}(q))} = 0
$$
for all $q \in \C \setminus P$. In different terms, $\theta$ has a constant eigenvalue over $\C\setminus P$; 
in particular, the residue of this eigenvalue at infinity is $0$. This contradicts $\lambda^{\infty}_{k}\neq 0$ 
for all $k \in \{ 1,\ldots,n \}$ (see (\ref{mainii}) of Hypothesis \ref{main}).
\end{proof}
A basic property is the following. 
\begin{clm} \label{merospec}
The points of $\Sigma_{\xi }$ define a multi-valued meromorphic function of $\xi \in \hat{\C}$. 
\end{clm}
\begin{proof}
By assumption, $det(\theta_{\xi }(z))$ depends holomorphically on 
$\xi \in \hat{\C}$ and meromorphically on $z$. We conclude using the implicit function theorem, 
namely that the solutions of a meromorphic equation depending holomorphically on a variable are 
meromorphic in this variable. 
\end{proof}
\begin{defn} 
The graph of the multi-valued meromorphic function 
\begin{align*}
     \hat{\C} \setminus \hat{P} & \lra \CP \\
     \xi & \mapsto \Sigma_{\xi }
\end{align*}
is called the \emph{spectral curve} of the Higgs bundle. It is denoted by $\Sigma $. 
\end{defn}
This object was first studied by N. Hitchin in \cite{Hit2}. 
By Claim \ref{merospec} the spectral curve is an analytic subvariety
$$ 
     \Sigma \xrightarrow{\j} (\hat{\C } \setminus \hat{P}) \times  \CP , 
$$
of (complex) dimension one. (Here $\j$ stands for inclusion.)
Moreover, by construction it is naturally a branched cover of $\hat{\C }$ via projection to the first factor. 

Here is an important property. 
\begin{prop} \label{red}
The spectral curve $\Sigma$ is reduced; in other words, $det(\theta_{\xi })$ vanishes only up to the first order except
for a finite set of points of $\Sigma$.
\end{prop}
\begin{proof} Suppose $\Sigma$ has infinitely many points $(q,\xi)$ where $det(\theta_{\xi })$ vanishes up to order higher
than one. Since $\Sigma$ has a natural extension into a compact curve in $\CP \times \CPt$ (see Section \ref{ext}), 
this means that for any $\xi$ some zero $q(\xi) \in\Sigma_{\xi}$ of $\theta_{\xi}$ has multiplicity higher than one; in different
terms, some irreducible component of $\Sigma$ has multiplicity higher than one. In particular, as $\xi \ra \infty$, at least two 
of the $q_{k}(\xi)$ must have the same Laurent expansions. This is impossible by (\ref{asymbehinf}) 
and the assumption $\lambda^{j}_{k}\neq \lambda^{j}_{k'}$ for $k \neq k'$ made in (\ref{maini}) of Hypothesis \ref{main}. 
\end{proof}

\subsection{Explicit computation of the hypercohomology}
Let us now compute the hypercohomology of 
\begin{equation} \label{sheafmap} 
        \E \xrightarrow{\theta_{\xi } \land }\F
\end{equation}
Consider arbitrary algebraic resolutions of the sheaves $\E$ and $\F$ such
that $\theta_{\xi } \land  $ induce a morphism of resolutions 
\begin{equation} \xymatrix{ K^{0,1} \ar[r]^{\theta_{\xi } \land }    &  K^{1,1} \\
              K^{0,0} \ar[r]^{\theta_{\xi } \land} \ar[u]^{\delta  } &
             K^{1,0} \ar[u]_{\delta  } \\
              \mathcal{E} \ar[r]^{\theta_{\xi } \land} \ar@{^{(}->}[u] & \mathcal{F.} \ar@{^{(}->}[u] 
} \label{algmr}
\end{equation}
For example, one might take resolutions by \v{C}ech cochains. 
By definition, the first filtration $K_{p}$ of the single complex associated to 
(\ref{algmr}) is given by 
\begin{align*}
     K_{0} & = (K^{0,1} \oplus  K^{0,0} ) \oplus ( K^{1,1} \oplus  K^{1,0} ) \\
     K_{1} & =  K^{1,1} \oplus  K^{1,0} . 
\end{align*}
The first page of the spectral sequence corresponding to this filtration is
given by 
\begin{equation} \xymatrix{ (\Cs^{0})^{[1]}(\CP)    &  (\Cs^{1})^{[1]}(\CP)  \\
              (\Cs^{0})^{[0]} (\CP) \ar[u]^{\delta  } &
            (\Cs^{1})^{[0]}(\CP) \ar[u]_{\delta  }
} \label{ss1}
\end{equation}
where $\Cs^{j}$ is the $j$-th cohomology sheaf of the map (\ref{sheafmap}),
and the vertical sequences come from resolutions 
\begin{align*}
       \Cs^{0} & \inc  (\Cs^{0})^{[0]} \xra{\delta } (\Cs^{0})^{[1]} \\
       \Cs^{1} & \inc  (\Cs^{1})^{[0]} \xra{\delta } (\Cs^{1})^{[1]}
\end{align*}
by taking global sections. Let us now describe explicitly the cohomology sheaves. 
Recall from definition \ref{specsetdef} that $q \in \Sigma_{\xi}$ are exactly the points where 
the map $\theta_{\xi }(q): E(q) \ra E(q)$ is not surjective. 
After all this preparation, we  have the following characterization: 
\begin{lem} 
The cohomology sheaf  $\Cs^{0}$  of order $0$ of the sheaf map 
(\ref{sheafmap}) is $0$. If $\det{(\theta_{\xi })}$ has a zero of order $1$
in all points of $q \in \Sigma_{\xi }$, then the first cohomology sheaf $\Cs^{1}$ is the
sky-scraper sheaf $\Sky_{\xi }$ whose stalk over a point $q \in \Sigma_{\xi }$ is the
finite-dimensional subspace $coKer(\theta_{\xi}(q)) \subset E(q)$, and all other stalks are $0$. 
\end{lem}
\begin{rk}
The cokernel of $\theta_{\xi}(q)$ is naturally identified with the orthogonal of the image with respect 
to the fiber metric, or, which is the same thing, with the kernel of $\theta_{\xi}^{\ast}(q)$. 
This allows us to think of $coKer(\theta_{\xi}(q))$ as a subspace of $E(q)$.
\end{rk}
\begin{proof}
Let us start with $\Cs^{0}$: suppose we have a section $\phi \in \E|_{U}$ on 
an open set $U \subset \CP$ such that $\theta_{\xi } \phi =0$. Since on the open subset 
$U \setminus \Sigma_{\xi }$  the map $\theta_{\xi }: E(q) \ra E(q) $ is an isomorphism,
we deduce that $\phi =0$ on this set. But a holomorphic section vanishing on
an open set vanishes everywhere, thus $\phi =0$ on all of $U$. This gives the
first statement of the lemma. 

We now come to $\Cs^{1}$: let $U \subset \CP$ be an open subset. If $U \cap  \Sigma_{\xi }
= \emptyset  $ then $\theta_{\xi }$ is an invertible holomorphic endomorphism of $\E$ on $U$, therefore 
$\Cs^{1}|_{U} =0$. Suppose now $U$ contains exactly one point $q \in \Sigma_{\xi}$. 
Then, for any section $\phi \in \E|_{U}$ the vector $(\theta_{\xi } \phi)(q) $ lies by
definition in the image of $\theta_{\xi }(q)$, which is just the orthogonal of
$coKer( \theta_{\xi }(q))$. Therefore, this latter is contained in
$\Cs^{1}|_{U}$. On the other hand, the condition that $\theta_{\xi }$ has a zero
of order $1$ in $q$ means that any section $\psi \in \E|_{U}$ such that 
$\psi(q) \bot coKer(\theta_{\xi }(q))$ is in $Im(\theta_{\xi })$. This proves the second statement. 
\end{proof}
\begin{rk} By Proposition \ref{red}, the condition of $\det{(\theta_{\xi })}$ having a first-order zero in
  all points of $ \Sigma_{\xi }$ is generic in $\xi$: it is verified for all $\xi$ except 
  for twice the eigenvalues of $\theta(q)$ for the finite number of points 
  $q$ of $\Sigma$ of multiplicity higher than one. For the discrete set of $\xi$
  where there exists a $q \in \Sigma_{\xi}$ with a multiple zero, one introduces
  the flag 
$$
    E(q) = F_{0}E(q) \supseteq coKer(\theta_{\xi }(q))=F_{1}E(q)\supseteq   \ldots
    \supset F_{r_{q}}E(q)=\{ 0\} , 
$$ 
the subscript of $F$ being the order of zero of $\theta_{\xi }^{\ast}(q)$ along the given subspace, 
and proves that the cohomology sheaf $\Cs^{1}|_{U}$ over an open set containing $q$ as the only 
element of $\Sigma_{\xi }$ is in this case equal to the jet space
$$
     \bigoplus_{m=1}^{r_{q}-1}F_{m}E(q).  
$$ 
The assumptions that for fixed $j \in \{ 1,\ldots,n \}$ all the $\lambda^{j}_{k} $ be different for $k \in \{ r_{j}+1,\ldots,r \} $ 
and for fixed $l \in \{ 1,\ldots,n' \}$ all the  $\lambda^{\infty }_{k}$  be different for $k \in \{ 1+a_{l}, \ldots a_{l+1} \}$ 
(see (\ref{maini}) and  (\ref{mainii}), Hypothesis \ref{main}), mean that in the punctures 
of $\CPt$ the limit states have first-order zeros. 
\end{rk}

Now since a resolution of the sky-scraper sheaf $\Sky_{\xi }$ is given by 
$$
     \Sky_{\xi } \inc \Sky_{\xi } \ra 0, 
$$
the first page of the hypercohomology spectral sequence (\ref{ss1}) becomes
$$
     \xymatrix{ 0    &  0  \\
              0 \ar[u]^{\delta  } &
             \bigoplus_{q \in \Sigma_{\xi }} coKer(\theta_{\xi }(q)). \ar[u]_{\delta  }
} \label{ss1b}
$$
All this implies the following: 
\begin{prop} \label{hss} The hypercohomology spectral sequence corresponding to the
  first filtration collapses in its first page, and we have a natural
  isomorphism 
$$
    \H^{1 }(\E \xrightarrow{\theta_{\xi } \land }\F) \simeq  \bigoplus_{q \in \Sigma_{\xi }} coKer(\theta_{\xi }(q)).
$$
\end{prop}
\begin{proof}
This is a consequence of the standard fact that a spectral sequence
collapses as soon as non-zero elements only appear in one of its rows. 
Furthermore, an explicit isomorphism can be given as follows: fix a radially invariant 
bump-function $\chi$ on the unit disk $\Delta \subset \C$, equal to $0$ on the boundary
of $\Delta$ and to $1$ in $0$, and such that $\d \chi $ is supported
on the annulus $ 1/3<r<2/3 $. For any complex number $a \neq 0$ set $\chi_{a }(z)=
\chi(z/a)$. Now choose $\varepsilon_{0}>0$ so that the distance in $\C $ between any two distinct points
of the finite set $P \cup \Sigma_{\xi }$ is at least $3\varepsilon_{0}$. 
For any $(v_{q})_{q \in \Sigma_{\xi }} \in \oplus  coKer(\theta_{\xi }(q))$ consider the section
$v_{\varepsilon_{0}} =\sum_{q \in \Sigma_{\xi } } v_{q} \chi_{\varepsilon }(z-q)$. 
Because $\d \chi_{\varepsilon_{0}}$ is
supported on the annulus $\varepsilon_{0}/3 <r< 2\varepsilon_{0}/3 $, the section 
$\dbar^{\E}(v_{\varepsilon_{0}} \d z)\in \Omega^{1,1}\otimes E$ is supported outside a neighborhood
of $\Sigma_{\xi }$. Since this latter is the zero set of $\det{(\theta_{\xi })}$, it then
follows that there exists a section $t_{\varepsilon_{0}} \d \bar{z} \in \Omega^{0,1}\otimes E$ such that 
$\theta_{\xi }\land (t_{\varepsilon_{0}}\d \bar{z}) + \dbar^{\E}(v_{\varepsilon_{0}}\d z)=0$, and
$t_{\varepsilon_{0}}$ is supported on the support of $\dbar^{\E} v_{\varepsilon_{0}}$, that is
outside a neighborhood of $\Sigma_{\xi }$ and of infinity. The couple 
$(v_{\varepsilon_{0}}\d z,t_{\varepsilon_{0}}\d \bar{z})$  therefore defines a cocycle in the
single complex associated to $ \mathcal{D}_{\xi }$, and using Proposition \ref{algl2} we can define a map 
\begin{align} \label{psi}
    \Psi_{\xi } : \bigoplus_{q \in \Sigma_{\xi }} coKer(\theta_{\xi }(q)) & \lra H^{1}(\mathcal{D}_{\xi }) 
     = \H^{1 }(\E \xrightarrow{\theta_{\xi } \land }\F) \notag \\
     (v_{q})_{q \in \Sigma_{\xi }} & \mapsto [(v_{\varepsilon_{0}}\d z,t_{\varepsilon_{0}}\d \bar{z})], 
\end{align}
where $[(v_{\varepsilon_{0}}\d z,t_{\varepsilon_{0}}\d \bar{z})]$ stands for 
the cohomology class in $H^{1}(\mathcal{D}_{\xi })$ of this couple. 

We need to show that this map does not depend on $\varepsilon_{0} >0$ chosen,
provided that it is sufficiently small as explained above. Consider 
therefore the section $v_{\varepsilon_{1}}$ for $\varepsilon_{1} < \varepsilon$. 
Since in the union of the disks of radius $ \varepsilon_{1}/3 $ around the
elements of $\Sigma_{\xi }$ we have $v_{\varepsilon_{1}} = v_{\varepsilon_{0}}$, and
$\theta_{\xi }$ is invertible outside this set, there exists a section $u \in \Gamma(E)$
such that $\theta_{\xi }u + v_{\varepsilon_{1}} \d z=v_{\varepsilon_{0}}\d z$. Then, as in the
proof of Proposition \ref{algl2}, the couple $(v_{\varepsilon_{0}}\d z,t_{\varepsilon_{0}}\d \bar{z})$
is equal to $(v_{\varepsilon_{1}}\d z+ \theta_{\xi }u ,t_{\varepsilon_{1}}\d \bar{z}+\dbar^{\E}u)$, and
the two couples define the same cohomology class in $H^{1}(\mathcal{D})$. 
This then allows us to fix $\varepsilon_{0}>0$ sufficiently small once and for all.

In a similar way, one can prove that $\Psi_{\xi }$ is independent of the actual cut-off function $\chi $ as well. 

Finally, the inverse of $\Psi_{\xi }$ can be obtained as follows: let the cohomology class 
$\eta \in  H^{1}(\mathcal{D}_{\xi })$ be represented by a $1$-form $\eta^{1,0} \d z + \eta^{0,1} \d \bar{z}$, where 
$\eta^{1,0}$ and $\eta^{0,1}$ are sections of $E$. Then we have 
\begin{equation} \label{psiinv}
     \Psi_{\xi }^{-1} \eta = (eval_{q } \eta^{1,0} )_{ q \in \Sigma_{\xi } }, 
\end{equation} 
where $eval_{q } \eta^{1,0}$ stands for evaluation of the section $\eta^{1,0}$ in the point $q$. 
\end{proof}
\begin{rk}
Notice that the formula (\ref{psiinv}) is independent of the $1$-form representative of $\eta$;
in particular, the $(1,0)$-part of the harmonic representative of a cohomology class 
$\Psi_{\xi }(v_{q})_{q \in \Sigma_{\xi }}$ vanishes in the $q \in \Sigma_{\xi }$ where $v_{q}=0$. 
\end{rk}

\section[Extension over the singularities]{Extension of the Higgs bundle over the singularities} \label{ext}
The interpretation of the holomorphic bundle underlying the transformed
Higgs bundle in terms of hypercohomology established in the 
previous section allows us to extend it over the singular points $\hat{P} \cup
\{ \infty \}$ in the parameter space ${\CPt}$. At each puncture, we need to do two 
things: first, define the fiber of the transformed vector bundle over it. 
This then extends the holomorphic structure induced by $\dbar^{\hat{\E}}$ over the puncture 
in a natural way: a holomorphic section through the singular point will be a 
continuous section in a neighborhood of it, that is holomorphic in the punctured neighborhood.  
(Continuity is defined at the same time as the exceptional fiber.) 
The second thing to do then is to give an explicit basis of holomorphic sections with respect 
to this extended holomorphic structure. It is important to note that the extensions $^{i}\hat{\E}$ 
we define here are \emph{not} the transformed extensions given in Definition \ref{transfext}, 
but rather ones induced by the original Higgs bundle, and for which computations are more comfortable. 
This is why we will call $^{i}\hat{\E}$ the \emph{induced extension}. 
We study the link between these two extensions in Section \ref{secttop}. 

\subsection{Extension to logarithmic singularities} \label{logext}
First, we consider the case of points of the set $\hat{P}$. 
We shall now describe the extension $^{i}\hat{\E}$ over such a point. 
Notice first that as the deformation $\theta_{\xi }$ has a well-defined extension over these 
points, its hypercohomology spaces are also well-defined there. In particular, in view of 
Proposition \ref{algl2}, we may extend the transformed vector bundle $\hat{V}$ by putting 
\begin{equation*}
        \hat{V}_{\xi_{l}}=\H^{1}(\E \xrightarrow{\theta_{\xi_{l} } \land }\F)
\end{equation*} 
This is the definition of the fiber over such a point. 

In order to give explicit representatives of holomorphic sections, 
let us examine what happens to the fiber $\hat{V}_{\xi }$
when $\xi $ approaches one of the points of $\hat{P}=\{ \xi_{1}, \ldots, \xi_{n'} \}$, say $\xi_{l}$. 
First, let us find the spectral points. 
\begin{clm}\label{nobranching}
As $\xi \ra \xi_{l}$, exactly $m_{l}=a_{l+1}-a_{l}$ branches of the meromorphic functions $q_{k}\in \Sigma_{\xi}$
converge to infinity, while all others remain in a bounded region of $\C$. Moreover, 
labelling the spectral points converging to infinity by $q_{1+a_{l}}(\xi),\ldots,q_{a_{l+1}}(\xi)$, 
they admit the asymptotic behavior 
\begin{equation} \label{zerobeh}
     {q}_{k}(\xi )=\frac{2\lambda_{k }^{\infty }}{(\xi-\xi_{l} )} + O(|\xi-\xi_{l}|^{-\delta}), 
\end{equation}
where $\delta>0$ can be chosen arbitrarily small. In particular, the branches converging to $\infty\in\CP$ of the 
spectral curve are not ramified over the point $\xi_{l}$. 
\end{clm} 
\begin{proof}
As it can be seen from (\ref{modelhiggsinf}), exactly $m_{l}$
of the eigenvalues of the leading order term near infinity of the Higgs  
field $\theta_{\xi }$ converges to $0$. Recall from Definition \ref{specsetdef} that $\Sigma_{\xi }$ is 
the vanishing set of $det(\theta_{\xi })$. This implies that (counted with multiplicities) exactly $m_{l}$ of the
points $q(\xi ) \in \Sigma_{\xi }$ converge to infinity; label these by $1+a_{l},\ldots,a_{l+1}$. All the other 
spectral points remain therefore bounded. By assumption (see (\ref{modelhiggsinf})) 
in a holomorphic trivialisation of the bundle $\E$ in a neighborhood of $\infty \in \CP$, ignoring the factor 
$\d z$ the field $\theta_{\xi }$ is of the form
$$
      \frac{1}{2} (A - \xi \Id) + \frac{C}{z} + O(z^{-2}), 
$$
where $O(z^{-2})$ stands for holomorphic terms independent of $\xi$. 
Suppose first that the field is exactly equal to the polar part in this formula, 
in other words the $O(z^{-2})$ term is equal to $0$. Then the solutions 
$\tilde{q}_{1}(\xi),\ldots,\tilde{q}_{r}(\xi)$ are clearly given by 
\begin{equation*} 
     \tilde{q}_{k}(\xi )=\frac{2\lambda_{k }^{\infty }}{(\xi-\xi_{l} )}.
\end{equation*}
In general, since $det(\theta_{\xi })$ is holomorphic in $z$, we can apply Rouch\'e's theorem to compare 
the position of the zeros of $det(\theta_{\xi })$ with those of the polar part studied above. 
This yields that the solutions $q_{k}(\xi ) \in \C$ of $det(\theta_{\xi })(q(\xi ))=0$ near infinity are close to 
$\tilde{q}_{k}(\xi )$; more precisely for any $\delta>0$, there exists $K>0$ such that for all $|\xi-\xi_{l}|$ 
sufficiently small we have 
$$
     \left\vert q_{k}(\xi )-\tilde{q}_{k}(\xi ) \right\vert < K |\xi-\xi_{l}|^{-\delta}.
$$
Remark here that as $\xi \ra \xi_{l}$ the behavior of $|\xi-\xi_{l}|^{-\delta}$ is small compared to 
$|\tilde{q}_{k}(\xi)|=c|\xi-\xi_{l}|^{-1}$. In other words, we have the expansion (\ref{zerobeh})
so that $q_{k}(\xi )$ converges indeed to infinity asymptotically proportionally to $(\xi-\xi_{l} )^{-1}$ for 
$a_{l} < k \leq  a_{l+1}$, while all other holomorphic families of zeros of $\det(\theta_{\xi })$ remain bounded. 

The condition that the $\lambda_{1+a_{l}}, \ldots, \lambda_{a_{l+1}}$ are all distinct (see (\ref{mainii}), Hypothesis
\ref{main}) now implies that there is no splitting of the solutions at infinity, 
that is to say locally near $\xi=\xi_{l}$ any $q_{k}(\xi)$ with $a_{l} < k \leq  a_{l+1}$ itself forms a meromorphic 
function without branching. Indeed, the occurrence of a branching at infinity implies that the Puiseux series 
of the corresponding solutions agree, which is not the case here because of the asymptotic behaviors 
(\ref{zerobeh}) with different leading coefficients.  
\end{proof}

Now, recall that for fixed $\xi \in \hat{\C} \setminus \hat{P}$, in the 
explicit description of $\hat{V}_{\xi }$ given in the proof of Proposition
\ref{hss}, we considered the zeros $q_{k}(\xi )$ for  $ k=1,\ldots,r$ of $det(\theta_{\xi })(q)$, 
and for each $q_{k}(\xi )$ an element $v_{k}(\xi )$ of the subspace 
$coKer(\theta_{\xi })_{q_{k}(\xi )} \subset E_{q_{k}(\xi )}$. Then we extended each
$v_{k}(\xi)$ holomorphically into a neighborhood of $q_{k}(\xi )$, and multiplied the section
we obtained by a bump-function equal to $1$ in a small disk around $q_{k}(\xi)$
and to $0$ on the boundary of a slightly larger disk. 
This section of $\F$ constituted the $(1,0)$-part of the element
in $\H^{1}(\E \xrightarrow{\theta_{\xi }} \F) \simeq \hat{V}_{\xi }$, and we chose the $(0,1)$-part in
such a way that the couple be in $Ker(D_{\xi }'')$. 
In what follows, we wish to do the same thing, but for all $\xi $ in a
neighborhood of $\xi_{l}$ at the same time. 

Let us consider one meromorphic family of zeros $q_{k}(\xi )$ with $a_{l} < k \leq  a_{l+1}$. 
We have just seen that $q_{k}(\xi )$ converges to $\infty $ as $\xi \ra \xi_{l}$;
therefore, we need to take a holomorphic section of $\E$ at infinity, extending an element of the
cokernel of $\theta_{\xi_{l}} $. One can check from formula (\ref{modelhiggsinf})
that this cokernel is equal to the vector subspace of the fiber $\F_{\infty }= E_{\infty } \otimes \d z$ generated by 
$\{ \sigma_{m}^{\infty } (\infty ) \d z \}_{m=1+a_{l}}^{a_{l+1}}$, where $\{ \sigma_{m}^{\infty } \}_{m=1}^{r}$ 
is the holomorphic trivialisation of $\E$ at infinity considered in (\ref{holtrivinf}).
Furthermore, since the metric $h$ is mutually bounded with the diagonal model 
$$
     diag(|z|^{-2\alpha^{\infty}_{k}}), 
$$
the orthogonal of the image of $\theta_{\xi}$ in $E(q_{k}(\xi))$ converges to $\sigma_{k}^{\infty}(\infty)$ as $\xi \ra \xi_{l}$. 
Let $\varsigma_{k}(z )$ be a holomorphic extension of $\sigma_{k}^{\infty }(\infty) $ to a neighborhood of infinity such that for any 
$\xi \in \hat{\C}$ sufficiently close to $\xi_{l}$, the vector $\varsigma_{k}(q_{k}(\xi ) )\d z$ be in the cokernel of 
$\theta_{\xi }(q_{k}(\xi ))$. 
Such an extension exists because $\theta_{\xi }$ varies holomorphically with $\xi $ and by Claim \ref{nobranching} 
$q_{k}(\xi )$ is a genuine (single-valued) meromorphic function of $\xi$. 
A holomorphic section $\hat{\sigma}_{k}$ of $\hat{\E}$ around $\xi_{l}$ is then given by the section constructed as follows: 
for $\xi $ sufficiently close to $\xi_{l}$ such that $\varsigma_{k} $ is defined in $q_{k}(\xi ) $, set 
\begin{equation} \label{vlog}
    v_{k } (z, \xi ) = \chi_{\varepsilon_{0} (\xi -\xi_{l})^{-1}} (z- q_{k}(\xi )) \varsigma_{k}(z ) ,
\end{equation}
where we recall from the proof of Proposition \ref{hss} that $\chi_{\varepsilon_{0} (\xi -\xi_{l})^{-1}}$ is 
a bump-function on a disk centered at $0$ and of diameter $\varepsilon_{0} |\xi -\xi_{l}|^{-1}$ with 
$\varepsilon_{0}$ sufficiently small only depending on the parameters of the initial connection, 
fixed once and for all. 
(The importance of this choice will become clear in Theorem \ref{parweightthmlog}.) 
Also, let $t_{k } (z,\xi ) \d \bar{z} \in \Gamma (\C , E \otimes \Omega^{0,1})$ be the unique solution of the equation 
\begin{equation} \label{tlog}
     \dbar^{\E}  v_{k }(z,\xi ) \d z = - \theta_{\xi } t_{k }(z,\xi ) \d \bar{z}. 
\end{equation}
Then consider the cohomology class $\hat{\sigma}^{l}_{k}(\xi)$ in 
$\H^{1 }(\E \xrightarrow{\theta_{\xi } \land }\F) \simeq \hat{V}_{\xi } $ of the couple 
$(v_{k } (z, \xi ) \d z, t_{k } (z,\xi ) \d \bar{z})$ defined as above. Since the choice of $\varsigma_{k}$ is independent
of $\xi $ and moreover $\theta_{\xi }$ and  $ q_{k}(\xi )$ depend holomorphically on $\xi $, it follows that 
$\hat{\sigma}^{l}_{k}$ is $\dbar^{\hat{\E}}$-holomorphic in $\xi $ outside of $\xi_{l}$. 
\begin{defn}
Let the extension $^{i}\hat{\E}$ of $\hat{\E}$ to $\xi_{l}$ be defined by the holomorphic trivialisation 
given by the sections $\hat{\sigma}^{l}_{k}$ for all choice of $k \in \{ 1+a_{l}, \ldots, a_{l+1} \}$ and for some 
holomorphic extension $\varsigma_{k}$ of $\sigma^{\infty }_{k}(\infty )$ such that for any $\xi \in \hat{\C}$ sufficiently close to 
$\xi_{l}$, we have $\varsigma_{k}(q_{k}(\xi ) )\d z \in coKer(\theta_{\xi }(q_{k}(\xi ) ) )$.
\end{defn}

\subsection{Extension to infinity} \label{infext} 
In order to define the fiber over infinity, we first rephrase what we have done until now to obtain the
holomorphic bundle $\hat{\E}=(\hat{V},\dbar^{\hat{\E}})$ underlying the 
transformed Higgs bundle: we considered the sheaves $\E$ and $\F$ over $\CP$, we pulled them back to 
$\CP \times \hat{\C}$ by the projection map $\pi_{1}$ on the first factor, and formed the sheaf map 
$$
    \pi_{1}^{\ast }\E \xra{\theta_{\bullet } }  \pi_{1}^{\ast }\F 
$$
equal to $\theta_{\xi }$ on the fiber $\CP \times \{ \xi \}$. We then defined the vector
bundle 
$$
      \hat{V}_{\bullet  } = \H^{1}(\pi_{1}^{\ast }\E \xra{\theta_{\bullet  }  } \pi_{1}^{\ast }\F),
$$
over $\hat{\C} \setminus \hat{P}$ and we let $\dbar^{\hat{\E}}$ be the 
partial connection induced by $\hat{\d}^{0,1}$. In what follows, we keep on
writing $\E$ and $\F $ for their pull-back to the product, whenever this
does not cause confusion. Notice that $\theta_{\bullet }$ is
holomorphic in both coordinates. We wish to extend the hypercohomology of
this sheaf map over infinity; we will be done if we can extend the map
$\theta_{\bullet}$  over infinity in a holomorphic manner. Indeed, the hypercohomology 
of a holomorphic family of sheaf morphisms is a holomorphic vector bundle 
over the base space of the deformations, in our case $\CPt$. Notice that by definition 
$\theta_{\xi }= \theta - \xi/2 \d z \land $, so it becomes singular as we let $\xi $ converge to infinity. However, we can
slightly change the sheaf $\F $ in such a way that there exist a natural
extension of $\theta_{\bullet }$. Again, we follow \cite{Jardim}.  

Consider the projections $\pi_{j}$ to the $j$-th coordinate in the product
manifold $\CP \times \widehat{\mathbf{CP}}^{1}$, and set 
$\tilde{\F } = \pi_{2}^{\ast }\mathcal{O}_{\widehat{\mathbf{CP}}^1}(1) \otimes \F  $. 
Recall that $\mathcal{O}_{\widehat{\mathbf{CP}}^1}(1)$ admits two global
holomorphic sections $s_{0}$ and $s_{\infty }$, characterized by the
fact that if $\hat{U}_{0}$ and $\hat{U}_{\infty }$ are the standard neighborhoods of 
$0\in \widehat{\mathbf{CP}}^1$ and $\infty \in  \widehat{\mathbf{CP}}^1$ with coordinates $\xi $ and $\zeta=\xi^{-1} $
vanishing in $0$ and $\infty $ respectively, then we have 
\begin{align}
              s_{0}(\xi )= \xi && s_{\infty }(\xi )=1 && \mbox{in } \hat{U}_{0} \label{o1zero} \\
              s_{0}(\zeta )= 1 && s_{\infty }(\zeta )=\zeta && \mbox{in } \hat{U}_{\infty }. \label{o1inf}
\end{align}
Notice that here $\xi $ is the standard coordinate of $\C$ we used to define
$\theta_{\xi }$. Therefore for $\eta \in \widehat{\mathbf{CP}}^1$ we put 
\begin{align}
       \tilde{\theta}_{\eta } & : \E \lra \tilde{\F } \\
       \tilde{\theta}_{\eta } & =  s_{\infty }(\eta ) \otimes \theta  -\frac{1}{2} s_{0}(\eta ) \otimes \d z \land , \label{globdefhiggs}
\end{align} 
We remark that by (\ref{o1zero}), on $\hat{U}_{0}=\C $ we have 
$\tilde{\theta}_{\xi }= \theta - \xi/2 \d z \land =\theta_{\xi }$, so $\tilde{\theta}_{\bullet }$ is indeed an extension of the 
deformation $\theta_{\bullet }$ to infinity. 
Therefore, in what follows we keep on writing $\theta $ for $\tilde{\theta}$ whenever this does
not cause any confusion.
In the same manner, we see that 
$$
 {\theta }_{\infty }= - \frac{1}{2} s_{0}(\xi )_{\xi = \infty } \otimes \d z \land : \E \lra \F \otimes \mathcal{O}_{\CPt }(1)_{\xi = \infty }. 
$$
From the definition of the sheaves $\E$ and $\F$ one can see 
that the cohomology sheaves of
this map are $\Cs^{0}(\d z \land )=0$ and $\Cs^{1}(\d z \land )= \Sky_{\infty }$, 
the sky-scraper sheaf supported in points of $P$ and
having stalk equal to $ s_{0}(\xi )_{\xi = \infty } \otimes coKer(Res(\theta,p)) $ in $p \in P$. 
Therefore, as in Proposition \ref{hss}, we obtain that the first
hypercohomology space of this map equals $  s_{0}(\xi )_{\xi = \infty } \otimes \left( \oplus_{p\in P} coKer(Res(\theta,p)) \right)$, 
and all its other hypercohomology spaces vanish. The extension of the vector bundle 
$\hat{V}$ to infinity is then given by setting $\hat{V}_{\eta  }
= \H^{1}(\E \xra{{\theta }_{\eta }} \tilde{\F} )$ for all $\eta \in \CPt \setminus \hat{P}$. 
In particular, any local section at $\zeta = 0$ of $\hat{E}$ is a family of sections of 
the sheaf $\tilde{\F}$, and therefore can be written 
\begin{equation} \label{holsectinf}
     s_{0}(\zeta ) \otimes \psi(z, \zeta ),
\end{equation}
where $\psi(z, \zeta )$ are sections of $\F $ depending on the parameter $\zeta $. 
\begin{defn} 
The extension $^{i}\hat{\E}$ of the holomorphic structure of $\hat{\E}$ to infinity is the extension  
whose holomorphic sections at infinity can be written as in (\ref{holsectinf}), with $\psi(z, \zeta )$  
holomorphic in $\zeta $. 
\end{defn}

We come to the explicit description of a holomorphic section of $^{i}\hat{\E}$ at $\xi = \infty $ 
with respect to this extension. We make a similar construction as in the case of logarithmic singularities: 
first, we make a basic remark. 
\begin{clm} \label{convzeros}
As $\xi \ra \infty $, all zeros of $det(\theta_{\xi })$ converge to one of the points of $P$. Moreover, supposing 
$q(\xi)\ra p_{j}$, we have the asymptotic behavior 
\begin{equation} \label{asymbehinf}
     q(\xi) 
     = p_{j}+2 \frac{\lambda_{k}^{j}}{\xi} + O(\xi^{-2+\delta}),  
\end{equation}
where $\lambda_{k}^{j}$ is a non-vanishing eigenvalue of the residue of $\theta$ at $p_{j}$ and $\delta>0$ can be chosen 
arbitrarily small. In particular, the spectral curve is not branched over the point $\xi=\infty$. 
\end{clm}
\begin{proof}
Let us consider the deformation of the Higgs field in terms of the coordinate $\zeta=\xi^{-1 }$ in $\hat{U}_{\infty }$. 
As we see from (\ref{o1inf}) and (\ref{globdefhiggs}), it is given by 
$$
     {\theta }_{\zeta }=\zeta \theta - \frac{1}{2} \d z \land . 
$$
Notice that as $\zeta \ra 0$, the first term on the right-hand side in a fixed
point $z \in \CP \setminus P$ becomes insignificant, and $\theta_{\zeta }(z)$ converges to
$-1/2\d z \land $. Therefore, for $|\zeta |$ 
sufficiently small, all zeros of $det(\theta_{\xi })$ are in a neighborhood of $P$. 
In order to determine the asymptotic of this convergence,
remember that in a holomorphic trivialisation of $E$ in some neighborhood of $p_{j}$ the Higgs field is equal to the
model (\ref{modelhiggs}) up to terms in $O(z-p_{j})$. As in the case $\xi \ra \xi_{l}$, the solutions are close 
to those of the diagonal model $det(diag(\theta_{\zeta}(\tilde{q})))=0$ (see Claim \ref{nobranching}). 
This equation is 
$$
     \Pi_{k=1}^{r} \left(\frac{\zeta \lambda_{k}^{j}}{\tilde{q}-p_{j}} - \frac{1}{2} \right)  =0. 
$$ 
The solutions $\tilde{q}_{k}^{j}(\zeta )$ are clearly given by  
$$
    \tilde{q}_{k}^{j}(\zeta )  = p_{j}+2 \zeta \lambda_{k}^{j} = p_{j}+2 \frac{ \lambda_{k}^{j}}{\xi }.
$$
Here the upper index of the solution stands for the point $p_{j}\in P$ it converges to, and 
the lower index $k \in \{r_{j}+1,\ldots,r\}$ is determined by the extension of the cokernel of $\theta_{\zeta}$ 
at the point. An application of Rouch\'e's theorem gives again the claim. 

Finally, $\Sigma$ is not ramified at $\xi=\infty$ because this would imply that at least two of the $q_{k}(\xi)$ 
admit the same Puiseux expansion, which is impossible because of (\ref{asymbehinf}) and (\ref{maini}) of 
Hypothesis \ref{main}. 
\end{proof}
Furthermore, by Claim \ref{merospec} the points of $\Sigma_{\xi }$ define a 
multi-valued meromorphic function in the variable $\xi $ near infinity. 
Let $q_{k}^{j}(\xi ) \in \Sigma_{\xi }$ be such a holomorphically varying zero of 
$det(\theta_{\xi } )$, and suppose it converges to $p_{j} \in P$ as $\xi \ra \infty $. 
We can let the index $k$ to vary from $r_{j}+1$ to $r$. 
Consider the diagram 
$$
\xymatrix@!=16pt{
             & \hspace{3pt} \Sigma_{} \ar@{^{(}->}[d]^(.4){\j} & \\
             & \CP \times \CPt \ar[dl]_(.55){\pi_{1}} \ar[dr]^{\pi_{2}} & \\
             \CP & & \CPt    }
$$
where $\j$ is inclusion and the two other arrows are canonical projections. In order to define a local 
holomorphic section of the transformed bundle, we need to choose elements of $coKer(\theta_{\xi }(q_{k}^{j}(\xi )))$ 
for all $\xi $, such that they depend holomorphically with $\xi $. It is clear that this is equivalent 
to choose a local holomorphic section $\psi$ of $\j^{\ast } \pi_{1}^{\ast } \F$ over the branch $(q_{k}^{j}(\xi  ),\xi)$ 
near the point $(p_{j}, \infty )$ such that for all $\xi $, we have $\psi(q_{k}^{j}(\xi ), \xi) \in coKer(\theta_{\xi }(q_{k}^{j}))$.
Since any local section of $\F $ near $p_{j}$ multiplied by $(z-p_{j})$ is a local section of 
the sheaf $\E \otimes \d z$, the section $(q_{k}^{j}(\xi )-p_{j}) \psi $  of $\j^{\ast } \pi_{1}^{\ast } \F$ 
near $(p_{j},\infty )$ is in fact a local holomorphic section of $\j^{\ast } \pi_{1}^{\ast } (\E \otimes \d z)$ on the branch 
$( q_{k}^{j}(\xi ),\xi )$ of the spectral curve $\Sigma \subset \CP \times \CPt$. Furthermore, because of Claim \ref{convzeros}, 
$(q_{k}^{j}(\xi ), \xi ) \mapsto q_{k}^{j}(\xi )$ is a simple cover near $p_{j}$ without branching. 
In particular, for all $q$ sufficiently close to $p_{j}$ there exists a unique $\xi (q)$ 
such that $q=q_{k}^{j}(\xi (q))$. Therefore, $(q_{k}^{j}(\xi )-p_{j}) \psi(q_{k}^{j}(\xi ), \xi)$ is the lift from $\CP$ of a section 
$\varsigma_{k}^{j} (z) \d z$ of $\E \otimes \d z$ in a neighborhood of $p_{j} $, such that for all $q$ we have 
\begin{equation} \label{varsigmacoker}
   \varsigma_{k}^{j} (q) \d z \in coKer(\theta_{\xi (q)}(q)). 
\end{equation}
In particular, $\varsigma_{k}^{j} (p_{j}) \d z \in coKer(\theta_{\infty } (p_{j}))=E_{sing} \otimes \d z$, as it can easily 
be checked using formula (\ref{globdefhiggs}). Conversely, we may consider any section $\varsigma_{k}^{j}(z)$
satisfying (\ref{varsigmacoker}), lift $\varsigma_{k}^{j}(z) \d z$ to a section of $\j^{\ast } \pi_{1}^{\ast } (\E \otimes \d z)$, 
and divide the result by $q-p_{j}$ to obtain $\psi$. 
Fix now for all $k=\{ r_{j}+1,\ldots,r \}$ a section $\varsigma_{k}^{j} $ satisfying (\ref{varsigmacoker}). 
All that we have said above motivates the definition:
\begin{equation} \label{vinf}
      v_{k}^{j}(z, \xi )= \chi_{\varepsilon_{0} \xi^{-1}}(z-q_{k}^{j}(\xi  )) \frac{ \varsigma_{k}^{j} (z) }{z-p_{j}} \otimes s_{0} (\xi ),
\end{equation}
where we recall again from the proof of Proposition \ref{hss} that $\chi_{\varepsilon_{0}\xi^{-1 }}$ is a 
bump-function over the disk of radius $\varepsilon_{0}/|\xi |$. Remark that evaluation of 
$v_{k}^{j}(z, \xi ) \d z$ in $z=q_{k}(\xi )$ is by definition in
the cokernel of $\theta_{\xi }$. Also, as in the case of logarithmic singularities, for all $\xi $ close to 
infinity, let $t_{k}^{j}(z,\xi )$ be the unique section of $E $ satisfying the equation (\ref{tlog}) 
for all $z$, in other words such that $D_{\xi }''(v_{k}^{j}(z,\xi ) \d z, t_{k}^{j}(z,\xi ) \d \bar{z} )=0$. 
A holomorphic trivialisation of $^{i}\hat{\E}$ at infinity is then given by the $D_{\xi }''$-harmonic representatives 
$\hat{\sigma}_{k}^{\infty}(\xi )$ of the couples $(v_{k}^{j} (z,\xi ) \d z, t_{k}^{j}(z, \xi ) \d \bar{z} )$ for all  
$k=\{ r_{j}+1,\ldots,r \}$ and all $j=\{1,\ldots,n\}$.

\section{Singularities of the transformed Higgs field} \label{sing}
In this part, we describe the eigenvalues of the singular parts of the transformed Higgs
field $\hat{\theta}^{H}$ at the singularities. This establishes points (\ref{iv}), (\ref{vi}) and (\ref{vii}) 
of Theorem \ref{mainthmhiggs}. 

\subsection{The case of a logarithmic singularity} 
Recall from (\ref{htrhiggs}) that the transformed Higgs field is defined as multiplication by the coordinate 
$-z/2$ of a harmonic spinor, followed by projection onto harmonic forms. 
\begin{lem} \label{treigen}
The set of eigenvalues of the transformed Higgs field $\hat{\theta }^{H}$ on the fiber $\hat{E}_{\xi }^{H}$ 
(with multiplicities) is equal to $-\Sigma_{\xi }/2$ (with multiplicities), where $\Sigma_{\xi }$ is the set of zeros 
of $det(\theta_{\xi })$. 
\end{lem}
\begin{proof}
Let a cohomology class in the space $ \hat{E}_{\xi }^{H}= H^{1}(\mathcal{D}_{\xi })$ (see \ref{dcomplex})
be represented by $1$-forms $(v(\xi ) \d z, t(\xi ) \d \bar{z}) \in (\Omega^{1,0} \oplus \Omega^{1,0}) \otimes E$. 
Since this spinor is not harmonic, first of all we need a technical result:  
\begin{clm} Let $(v(\xi )\d z, t(\xi )\d \bar{z}) \in (\Omega^{1,0} \oplus \Omega^{1,0}) \otimes E $ be annihilated by $D_{\xi }''$. 
Then we have 
$$
     \hat{\pi}^{H}_{\xi }(z  \hat{\pi}^{H}_{\xi }(v(\xi )\d z, t(\xi )\d \bar{z}) ) 
     = \hat{\pi}^{H}_{\xi }(z (v(\xi )\d z, t(\xi )\d \bar{z}) ).
$$
In words, the action of the Higgs field can be computed on any representative section in 
$Ker(D_{\xi }'')$. 
\end{clm}
\begin{proof} This is straightforward: we need to show 
$$
      \hat{\pi}_{\xi }^{H}(z  (\mbox{Id} -\hat{\pi}^{H}_{\xi })(v(\xi )\d z, t(\xi )\d \bar{z}) ) = 0, 
$$
which is equivalent to 
$$
      z \Dir_{\xi }'' G_{\xi } (\Dir_{\xi }'')^{\ast } (v(\xi )\d z, t(\xi )\d \bar{z}) \bot \hat{E}_{\xi }^{H}. 
$$
Now the only thing to remark is that if  $(v(\xi )\d z, t(\xi )\d \bar{z}) \in Ker(D_{\xi }'')$, then this implies that
$$
     (\Dir_{\xi }'')^{\ast } (v(\xi )\d z, t(\xi )\d \bar{z}) = (D_{\xi }'')^{\ast } (v(\xi )\d z,t(\xi )\d \bar{z}) \in \Omega^{0}\otimes E,
$$ 
and by diagonality of $G_{\xi }$ with respect to the decomposition $S^{+}\otimes E  = (\Omega^{0}\otimes E) \oplus (\Omega^{2}\otimes E)$ 
(see Lemma \ref{diaglem}), also 
$$
    G_{\xi } (\Dir_{\xi }'')^{\ast } (v(\xi )\d z, t (\xi ) \d \bar{z})\in \Omega^{0} \otimes E. 
$$
Therefore we have
$$
  \Dir_{\xi }'' G_{\xi } (\Dir_{\xi }'')^{\ast } (v(\xi )\d z, t(\xi )\d \bar{z}) 
  = D_{\xi }'' G_{\xi } (D_{\xi }'')^{\ast } (v(\xi )\d z,t(\xi )\d \bar{z}), 
$$
and we conclude using the commutation relation 
$$
        [z, D_{\xi }'']=0
$$
combined with $Im(D_{\xi }'') \bot \hat{E}_{\xi }^{H}$. 
\end{proof}

The proof of the lemma is now immediate: via the map (\ref{psiinv}), 
\[
     \Psi_{\xi }^{-1} (z ( v(\xi )\d z, t(\xi )\d \bar{z} )) = ( q \cdot  eval_{q } v(\xi ) )_{ q \in \Sigma_{\xi } }
\]
multiplication by $z$ goes over to multiplication by $q$ in the point $q \in \Sigma_{\xi }$, and via (\ref{psi}) 
this is then re-transformed into multiplication by the constant $q$ on the component of $v(\xi )$ 
localized near $q$. 
\end{proof}

\begin{thm} \label{eigenlog}
The eigenvalues of the transformed Higgs field $\hat{\theta }^{H}$ have first-order poles in the points of
$\hat{P}$. Furthermore, the non-vanishing eigenvalues of its residue in the puncture $\xi_{l}$ are 
equal to $\{ -\lambda^{\infty }_{1+a_{l}}, \ldots , -\lambda^{\infty }_{a_{l+1}} \}$, where $\{ \lambda^{\infty }_{1+a_{l}}, \ldots, \lambda^{\infty }_{a_{l+1}} \}$ 
are the eigenvalues  of the residue of the original Higgs field $\theta $ at infinity, 
restricted to the eigenspace of $A$ corresponding to the eigenvalue $\xi_{l}$.
\end{thm}
\begin{proof}
As we have seen in (\ref{zerobeh}), the point $q_{k}(\xi ) \in \Sigma_{\xi }$ converges to infinity at
the first order with $2\lambda_{k}^{\infty }(\xi - \xi_{l})^{-1}$  as $\xi \ra \xi_{l }$, where $k \in \{ 1+a_{l},\ldots,a_{l+1} \}$ 
is an index such that the eigenvalue $\lambda_{k}^{\infty }$ of the residue term of $\theta $ at infinity appears in the 
eigenspace of the second order term $A$ corresponding to the eigenvalue $\xi_{l }$. By Lemma \ref{treigen}, 
the transformed Higgs field has a logarithmic singularity at $\xi_{l}$, and 
the corresponding residue is $-\lambda_{k}^{\infty }$. 
\end{proof}

\subsection{The case of infinity} 

We wish to show the following. 
\begin{thm} \label{eigeninf}
The transformed Higgs field has a second order singularity at infinity. The set of eigenvalues of 
its leading order term is $\{ -p_{1}/2, \ldots , -p_{n}/2 \}$, where $\{ p_{1}, \ldots , p_{n} \}=P$ is the 
set of punctures of the original Higgs bundle. The multiplicity of the eigenvalue $-p_{j }/2$ 
is equal to $r-r_{j} = rk (Res(\theta , p_{j}))$. The set of eigenvalues of the residue of the transformed 
Higgs field restricted to the eigenspace of the second-order term corresponding to the 
eigenvalue $-p_{j }/2$ is $\{ -\lambda^{j}_{k} \}_{k \in \{ r_{j}+1, \ldots , r \}}$. 
\end{thm}
\begin{proof}
In Claim \ref{convzeros} we have proved that as $\zeta \ra 0$, all zeros of $det(\theta_{\zeta } )$ must converge to
one of the points of $P$. Furthermore, the expansion of a spectral point $q_{k}$ converging to 
$p_{j}$ is (\ref{asymbehinf}). 
By Lemma \ref{treigen}, on the corresponding components $\hat{\theta }^{H}$ is just multiplication 
by $- \Sigma_{\xi} \d \xi /2$. Hence, we see that the eigenvalues of the leading-order term of the transformed 
Higgs field are equal to $\{ -p_{j}/2 \}_{ j=1,\ldots,n}$, while those of its first-order term are 
$\{-\lambda_{k}^{j} \}_{j=1,\ldots,n;k=r_{j}+1,\ldots,r}$. 
\end{proof}

\section{Parabolic weights} \label{parweightsect}
Here we compute the parabolic weights of the transformed Higgs bundle with respect to the induced 
extension. 

\subsection{The case of infinity} 
\begin{thm} \label{parweightthminf}
The parabolic weight of the extension $^{i}\hat{\E}$ of the transformed Higgs bundle at infinity described in Subsection 
\ref{infext}, restricted to the eigenspace of 
$\hat{\theta }$ corresponding to the eigenvalue $-p_{j}/2$ of its second order term and the eigenvalue 
$- \lambda^{j}_{k}$ of its residue is equal to $-1+\alpha^{j}_{k}$, where $\alpha^{j}_{k}$ is the parabolic weight on 
the $\lambda^{j}_{k}$-eigenspace of the residue of the original Higgs bundle at $p_{j}$. 
\end{thm}
\begin{proof} 
We prove the statement in two steps. In the first one, we show that it is true supposing the 
original Higgs bundle only has one logarithmic point of a precise form. 
In the second one, we show how the case with only one 
logarithmic point and the exponential decay results of Section \ref{harmspinsect} imply the general case. 
\paragraph{Step 1.} Let us first suppose that the set of logarithmic singularities is reduced to 
a single point $p_{1}$, that we may take to be $0$ without restricting generality. Furthermore, 
we suppose that $E$ is a holomorphically trivial bundle over $\C$ and that in a global holomorphic 
trivialisation $\{ \sigma_{k} \}$ the Higgs field is equal to
$$
     \theta = diag\left( \frac{\lambda_{k}}{z } \right)_{k=1,\ldots,r} \d z
$$
and the metric is just
\begin{equation} \label{modeleucl}
      h(\sigma_{k}, \sigma_{k}) =  |z|^{2\alpha_{k}}.
\end{equation}
This defines a parabolic Higgs bundle with weights $\alpha_{k}$ at $0$ and $-\alpha_{k}$ at infinity, the 
field having deformation 
\begin{equation} \label{higgsfield}
      \theta_{\xi } = diag\left( \frac{\lambda_{k}}{z } - \frac{\xi}{2}  \right)_{k=1,\ldots,r} \d z 
\end{equation} 
and the $D''$-operator 
\begin{equation}
      D''_{\xi } = \dbar +  diag\left( \lambda_{k} \frac{\d z}{z}  - \frac{\xi}{2} \d z  \right)_{k=1,\ldots,r} .
\end{equation}
Recall from Subsection \ref{infext} that a representative $(v_{\xi } \d z, t_{\xi } \d \bar{z})$ 
of any spinor $ \psi_{\xi }$ is supported in the finite collection of disks 
$ \cup_{q(\xi) \in \Sigma_{\xi } } \Delta ( q(\xi) , \varepsilon_{0}|\xi |^{-1})$. By Claim \ref{convzeros}, the points $q(\xi )$ are 
given by
\begin{equation} \label{vanish}
      q_{k} (\xi ) = \frac{2\lambda_{k}}{\xi }. 
\end{equation}
Define a family of homotheties indexed by $\xi \in \hat{\C } \setminus \hat{P}$ 
\begin{align} \label{homoth}
     h_{\xi }: \C & \lra \C \notag \\
              w & \mapsto z = \frac{w}{\xi };
\end{align}
in such a way that 
\begin{align}
      h_{\xi }^{-1}(0) & =  0 \notag \\
      h_{\xi }^{-1}( q_{k}(\xi ) ) & = 2\lambda_{k} \hspace{1cm} \mbox{for} \hspace{5mm} 
      k=r_{1},\ldots,r. \label{confzero}
\end{align}
Therefore, this corresponds to a family of coordinate changes $z \leftrightarrow w$ in the plane, such that the position of 
the zeros of the Higgs field $\theta_{\xi }$ after applying $h_{\xi }^{-1}$ is constant (the $2\lambda_{k}$ for $k=r_{1},\ldots,r$), 
as well as that of the poles ($0$ and $\infty $). Moreover, $\d z = \xi^{-1} \d w$ implies 
\begin{align} \label{confhiggs}
     h_{\xi }^{\ast } \theta_{\xi } =  diag\left[ \lambda_{k} \frac{\d w}{w}  -  
     \frac{1}{2} \d w  \right]_{k=1,\ldots,r} , 
\end{align}
and so 
\begin{align} \label{confdsecond}
     h_{\xi }^{\ast }D''_{\xi } & = \dbar +  diag\left[ \lambda_{k} \frac{\d w}{w}  -  
     \frac{1}{2} \d w \right]_{k=1,\ldots,r} , 
\end{align}
where $\dbar $ stands this time for the Dolbeault operator with respect to the $w$-coordinate. 
The crucial observation is that this operator is independent of $\xi$. 
On the other hand, remark that the Euclidean metric on the base space and the fiber metric (\ref{modeleucl})
behave under these coordinate changes as 
\begin{align}
     (h_{\xi })_{\ast } |\d w |^{2} & = |\xi |^{2} |\d z |^{2} \label{confeucl} \\
     |\sigma_{k}(z)|^{2} & =  |\xi |^{-2\alpha_{k} } |w|^{2\alpha_{k}} . \label{confherm}
\end{align}
In other words, if we denote by $h^{(w)}$ the model hermitian metric on $h_{\xi }^{\ast } E$ 
equal in the basis $h_{\xi }^{\ast } \sigma_{k}$ to 
$$
      h^{(w)} = diag(|w|^{2\alpha_{k}}),
$$
then the homotheties $h_{\xi}$ induce a family of tautological isomorphisms of Hermitian fiber bundles 
\begin{align} \label{confisom}
     ( h_{\xi }^{\ast } E,  h^{(w)}) & \lra  (E,  h) \\
     (h_{\xi }^{\ast } \sigma_{k})(w ) & \mapsto |\xi |^{\alpha_{k} } \sigma_{k}(z) \notag .
\end{align}
We deduce from (\ref{confeucl}) that in the basis $h_{\xi }^{\ast } \sigma_{k}$ the pull-back $h_{\xi }^{\ast } \Delta_{\xi }$ of the 
Laplacian of the Dirac operator $\Dir_{\xi }''$ has the form 
\begin{equation}  \label{conflapl}
     |\xi |^{2} \left[\Delta + diag \absl \frac{\lambda_{k} }{w}-\frac{1}{2}  \absr^{2}_{k=1,\ldots,r} \right]  , 
\end{equation}
where $\Delta$ stands for the usual Laplace operator on functions with respect to the 
metric $|\d w |^{2}$. The operator $\Delta^{(w)}$ between brackets in this formula 
is a bounded operator from the weighted Sobolev space $H^{2}(S^{+} \otimes E, |\d w|^{2})$ to $L^{2}(S^{+} \otimes E, |\d w|^{2})$. 
The weight at $0$ is determined by the condition that for a section $u \in H^{2}$ we have 
$u/|w|^{2} \in L^{2}$, and this gives therefore exactly the maximal domain of $\Delta^{(w)}$ 
(see Theorem \ref{Hilbertisom}). We infer that the pull-back $h_{\xi }^{\ast } G_{\xi }$ of the Green's operator 
of $\Delta_{\xi }$ is 
\begin{equation}  \label{confgreen}
     |\xi |^{-2} G^{(w)}, 
\end{equation}
where $G^{(w)}$ is the inverse of $\Delta^{(w)}$. It also follows from Theorem \ref{Hilbertisom} that $G^{(w)}$ 
is a bounded linear operator from $L^{2}(S^{+} \otimes E, |\d w|^{2})$ to $H^{2}(S^{+} \otimes E, |\d w|^{2})$. 
Because $\Delta^{(w)}$ is diagonal in the basis $\sigma_{k}$, the same is true for $G^{(w)}$.
Remark that the pull-backs $h_{\xi }^{\ast }\hat{\pi}_{\xi}$ of the orthogonal projections onto $\Delta_{\xi}$-harmonic spinors
are all equal to the orthogonal projection $\hat{\pi}^{(w)}$ onto $\Delta^{(w)}$-harmonic spinors: indeed, the 
conformal factor $|\xi |^{2}$ in (\ref{conflapl}) changes neither the space of harmonic spinors nor the 
orthogonal projection operator onto them. In particular, since $\Delta^{(w)},G^{(w)}$ and $h$ are diagonal in the
basis $\sigma_{k}$, the same thing is true for all $\hat{\pi}_{\xi}$. 

Now notice that by the definition of the $\dbar^{\hat{\E}}$-holomorphic extension to infinity of the
transformed bundle given in (\ref{vinf}) and via the identification (\ref{confisom}), the sections 
$|\xi|^{\alpha_{k}} h_{\xi}^{\ast } (v_{k}(z,\xi) \d z)$ (modulo the value of the section $s_{0}$ of $\mathcal{O}_{\CPt }(1)$) 
coincide: indeed, 
$$
   |\xi|^{\alpha_{k}} \chi_{\varepsilon_{0}/\xi}(z-q_{k}(\xi )) \sigma_{k}(z) \frac{\d z}{z} = 
   \chi_{\varepsilon_{0}} \left(w - 2\lambda_{k}  \right) (h_{\xi }^{\ast }\sigma_{k})(w) \frac{\d w}{w}.
$$ 
It then follows from formula (\ref{confdsecond}) together with the definition
(\ref{tlog}) that the coefficient of $s_{0}$ in $|\xi|^{\alpha_{k}} h_{\xi}^{\ast } t_{k}(z,\xi) \d \bar{z}$ 
is also independent of $\xi$. From the fact that the projections $\hat{\pi}_{\xi}$ are also constant, we deduce that
the coefficient of $s_{0}$ in the pull-back 
\begin{equation}\label{confspinor}
      (h_{\xi}^{\ast }\hat{\sigma}_{k})(w,\xi) =   |\xi|^{\alpha_{k}} \hat{\sigma}^{\infty}_{k}(z,\xi)
\end{equation}
of the spinors $|\xi|^{\alpha_{k}} \hat{\sigma}^{\infty}_{k}(z,\xi)$ representing $ |\xi|^{\alpha_{k}}(v_{k}(z,\xi) \d z, t_{k}(z,\xi) \d \bar{z})$ 
does not depend on $\xi$. 
Therefore, denoting by $f_{k}(z,\xi)$ the coefficient of $s_{0}$ in $\hat{\sigma}^{\infty}_{k}(z,\xi)$ and by 
$(h_{\xi}^{\ast }f_{k})(w,\xi)$ the coefficient of $s_{0}$ in $(h_{\xi}^{\ast }\hat{\sigma}^{\infty}_{k})(w,\xi)$, we see 
by invariance of the $L^{2}$-norm of $1$-forms by conformal coordinate change that 
$$
     \int_{\C} \vert  f_{k}(z,\xi) \vert_{h,|d z|^{2}}^{2} |\d z|^{2} = 
     |\xi|^{-2\alpha_{k}}\int_{\C}\vert (h_{\xi}^{\ast }f_{k})(w,\xi) \vert_{h^{(w)},|d w|^{2}}^{2} |\d w|^{2}, 
$$
for all $\xi$, with the integral on the right-hand side a constant independent of $\xi$. 
On the other hand, recall from (\ref{o1zero}) that on the affine chart $\hat{U}_{0}$ of $\CP$ we have 
$s_{0}(\xi )= \xi $. Observe also that the transformed Hermitian metric $\hat{h}$ is defined in the chart
$\hat{U}_{0}$, and that for any harmonic spinor $f$ we have 
$$
     \hat{h}(\xi f,\xi f)=|\xi|^{2} \hat{h}(f,f) = |\zeta|^{-2} \hat{h}(f,f) 
$$
with $\zeta=\xi^{-1}$ the local coordinate centered at $0$ of the singularity at infinity. 
This means that the effect on the parabolic weights of multiplying by $s_{0}$ is adding $-1$. 
On the other hand, the $-\lambda_{k}$-eigenspace of the residue of the transformed Higgs bundle at infinity is spanned by
$\hat{\sigma}^{\infty}_{k}$. From all that has been said above, we deduce 
\begin{equation}\label{asymbehharmspin}
     \hat{h}(\hat{\sigma}^{\infty}_{k},\hat{\sigma}^{\infty}_{k})=M|\zeta|^{-2+2\alpha_{k}},
\end{equation}
where $M$ is independent of $\xi$; in different terms, that the parabolic weight of the transformed 
Higgs bundle at infinity on the $-\lambda_{k}$-eigenspace of the residue is equal to $-1+\alpha_{k}$. 

\paragraph{Step 2.}
Starting from now, we drop the assumption that the set of logarithmic singularities is reduced to a point. 
In this part, we patch together solutions to local problems provided by Step 1, and use the results of 
Section \ref{harmspinsect} to estimate the defect of these patched sections to be solutions of the global
problem. We find that the interaction between solutions to local problems near different punctures is 
small as $|\xi|$ gets large. 

Let $( \dbar^{\E}, \theta )$ be a Higgs bundle with some logarithmic singularities $P= \{ p_{1}, \ldots,p_{n} \}$. 
In a holomorphic trivialisation $\{ \sigma^{j}_{k} \}_{k=1}^{r} $ near each one of these points, up to terms 
in $O(1) \d z$, the Higgs field has the form  
\begin{equation}\label{dpjagain}
     \theta^{j} = \frac{A^{j}}{z-p_{j}} \d z , 
\end{equation}
where the ${A}^{j}$ are some diagonal matrices as in (\ref{dpj}). 
The deformation of these local models is 
$$
    \theta^{j}_{\xi } = \left[\frac{A^{j}}{z-p_{j}} - \frac{\xi}{2}\right] \d z, 
$$
and similarly the deformation of the local $D''$-operators $(D'')^{j}$ is 
\begin{align*}
    (D''_{\xi })^{j} & = \dbar^{\E } + \theta_{\xi } \\
    & = \dbar^{\E } + \left[ \frac{A^{j}}{z-p_{j}} - \frac{\xi}{2} \right] \d z, 
\end{align*}
Finally, that of the Dirac operator $\Dir^{j}= (D'')^{j} - ((D'')^{j})^{\ast }$ is  
$$
      \Dir^{j}_{\xi } = (D''_{\xi })^{j} - ((D''_{\xi })^{j})^{\ast },
$$
adjoint being taken relative to the harmonic metric corresponding to $(D'')^{j}$. 
Now for all $j$ we can consider the extension of $\theta^{j}$ to a trivial bundle $E^{j}$ over the whole plane 
by keeping the same formula (\ref{dpjagain}) for it, endowed with the model metric
$$
        h^{j}=diag(|z-p_{j}|^{2\alpha^{j}_{k}})_{k=1}^{r}. 
$$
It is clear that this extension only has one regular
singularity (in $p_{j}$) and an irregular one at infinity, so all the results of Step 1 hold for them. 
In particular, for representatives 
$$
     (v_{k}^{j}(z,\xi) \d z, t_{k}^{j}(z,\xi) \d \bar{z})
$$ 
as described in Subsection \ref{infext} we have a harmonic representative 
$$ 
   \hat{\sigma}^{\infty}_{k}(z,\xi) \in Ker(\Dir^{j}_{\xi })^{\ast } \subset H^{1}(\C, S^{-} \otimes E^{j})
$$ 
with 
$$
     \int_{\C} \absl \hat{\sigma}^{\infty}_{k}(z,\xi) \absr_{h^{j},|dz|^{2}}^{2} |\d z|^{2} = |\xi|^{2-2\alpha^{j}_{k}}.
$$
This growth is measured with respect to the diagonal model metric $h^{j}$; 
however, since the spinor $\hat{\sigma}^{\infty}_{k}$ is exponentially concentrated near $p_{j}$ and here $h^{j}$ 
is mutually bounded with the harmonic metric $h$ of $(\E,\theta)$, this implies 
\begin{equation}\label{diagharm}
    c|\xi|^{2-2\alpha^{j}_{k}} \leq \int_{\C} \absl \hat{\sigma}^{\infty}_{k}(z,\xi) \absr_{h,|dz|^{2}}^{2} |\d z|^{2}\leq C|\xi|^{2-2\alpha^{j}_{k}}
\end{equation}
for some $0<c<C$. Let $\chi^{j} $ be a cut-off function supported 
in a disk $\Delta(p_{j}, 3 \varepsilon_{0})$, equal to $1$ on $\Delta(p_{j}, 2\varepsilon_{0})$, such that 
$| \nabla  \chi^{j} | \leq K$. Then for $\varepsilon_{0}>0$ fixed sufficiently small, the global section of $S^{-} \otimes E$ defined by 
$$
    \hat{\sigma}(z,\xi) =  \chi^{j}(z) \hat{\sigma}^{\infty}_{k}(z,\xi)
$$
has a meaning, for the holomorphic trivialisation $\{ \sigma^{j}_{k} \} $ is defined in 
$\Delta(p_{j}, 3\varepsilon_{0})$ provided $\varepsilon_{0}$ is sufficiently small. 
Now notice that if $q(\xi ) \ra p_{j}$ as $\xi \ra \infty$ and more precisely 
$$
    q(\xi ) = p_{j} + \frac{2\lambda^{j}_{k}}{\xi } + O(|\xi |^{-2}),  
$$
in other words on the component of the transformed bundle with eigenvalue of the second-order part of 
$\hat{\theta}$ at infinity equal to $-p_{j}/2$ and eigenvalue of the residue of $\hat{\theta}$ at infinity equal to $-\lambda^{j}_{k}$, 
the holomorphic extension $\varsigma^{j}_{k}$ of the cokernel has as parabolic weight the $\alpha^{j}_{k}$ corresponding to 
the eigenspace of the eigenvalue $\lambda^{j}_{k}$ of the residue of $\theta$. 
Recall that the harmonic metric on the transformed side is just $L^{2}$-metric 
of the $\Delta_{\xi}$-harmonic representative with respect to the harmonic metric $h$ of the original Higgs bundle. 
The statement of the theorem will therefore follow once we prove that the harmonic representative of 
$\hat{\sigma}(z,\xi)$ satisfies the inequality 
\begin{equation} \label{asymbehharmspin2}
    c|\xi|^{2-2\alpha^{j}_{k}} \leq \int_{\C} \absl \hat{\pi}_{\xi}\hat{\sigma}(z,\xi) \absr_{h,|dz|^{2}}^{2} |\d z|^{2}\leq C|\xi|^{2-2\alpha^{j}_{k}}.
\end{equation}
for some $0<c<C$. Our first aim is to prove the following. 
\begin{lem} \label{lemstep4}
There exists $\delta>0$ and $K>0$ such that for $|\xi|$ sufficiently large the inequality 
$$
   \left\Vert \Dir_{\xi }^{\ast } \hat{\sigma}(\xi) \right\Vert_{L^{2}(\C )}^{2} 
   \leq K | \xi |^{2-2 \delta} \left\Vert  \hat{\sigma}(\xi) \right\Vert_{L^{2}(\C )}^{2}
$$
holds.
\end{lem}
\begin{proof}
Covering the annulus centered at $p_{j}$ of radii $2\varepsilon_{0}$ and $2R_{0}$ by a finite number of disks
of radius $\varepsilon_{0}$, we deduce from Lemmas \ref{expdecord} and \ref{expdecinf} that the 
$\Dir^{j}_{\xi }$-harmonic spinor $\hat{\sigma}^{\infty}_{k}(z,\xi)$ is concentrated in $H^{1}$-norm, up to a factor 
decreasing exponentially with $|\xi |$, in the disk $\Delta(p_{j},2\varepsilon_{0})$. 
In particular, it is concentrated up to an exponentially decreasing factor in the same disk in $L^{2}$-norm as well. 
Denoting by $\cdot $ Clifford multiplication, we have the estimation 
\begin{align*}
    \int_{\C } \absl \Dir_{\xi }^{\ast } (\chi^{j}(z)  \hat{\sigma}^{\infty}_{k}(z,\xi)) \absr^{2} | \d z |^{2} \leq & 
    \int_{\C } \absl \chi^{j}(z) \Dir_{\xi }^{\ast } \hat{\sigma}^{\infty}_{k}(z,\xi) \absr^{2} | \d z |^{2} \\
     & + \int_{\C } \absl ( \nabla  \chi^{j} )(z) \cdot  \hat{\sigma}^{\infty}_{k}(z,\xi) \absr^{2} | \d z |^{2} \\ 
    \leq & \int_{\Delta(p_{j}, 3 \varepsilon_{0}) } \absl \Dir_{\xi }^{\ast } \hat{\sigma}^{\infty}_{k}(z,\xi) \absr^{2} | \d z |^{2} \\
    & + K \int_{\Delta(p_{j}, 3 \varepsilon_{0}) \setminus \Delta(p_{j}, 2 \varepsilon_{0})} \absl \hat{\sigma}^{\infty}_{k}(z,\xi) \absr^{2} |\d z|^{2}. 
\end{align*}
Again, by Lemma \ref{expdecord} the second integral on the right-hand side is bounded by an exponentially decreasing
multiple of $\| \hat{\sigma}^{\infty}_{k}(z,\xi) \|^{2}_{L^{2}(\C)} $ as $|\xi| \ra \infty$. Therefore, we only need to treat 
$$
     \left\Vert \Dir_{\xi }^{\ast } \hat{\sigma}^{\infty}_{k}(z,\xi) \right\Vert_{L^{2}( \Delta(p_{j}, 3 \varepsilon_{0}) )}^{2}. 
$$
Remark that by hypothesis, 
$$
     (\Dir^{j}_{\xi })^{\ast } \hat{\sigma}^{\infty}_{k}(z,\xi) = 0, 
$$
so we have 
$$
     \Dir_{\xi }^{\ast } \hat{\sigma}^{\infty}_{k}(z,\xi) = \left[ \Dir_{\xi }^{\ast } -  (\Dir^{j}_{\xi })^{\ast } \right] \hat{\sigma}^{\infty}_{k}(z,\xi). 
$$
This is then bounded by 
$$
    \hat{\sigma}^{\infty}_{k}(z,\xi) O(|z-p_{j}|^{-1+\delta }),
$$
where $O(|z-p_{j}|^{-1+\delta })$ stands for a term bounded from above by a constant (independent of $\xi$) 
times $|z-p_{j}|^{-1+\delta }$, because $\Dir_{\xi}^{\ast }$ and $(\Dir^{j}_{\xi })^{\ast }$ are Dirac operators having 
the same local model at the puncture and their difference is clearly independent of $\xi$. 
In order to study this quantity, we make use of the coordinate $w= \xi(z -p_{j}) $ analogously to that 
introduced in (\ref{homoth}). Under this coordinate change, the disk $\Delta(p_{j},3 \varepsilon_{0})$ goes 
into the (varying) disk $\Delta(0, 3 \varepsilon_{0}|\xi|)$. Hence, we need to prove
\begin{align*}
       \int_{\Delta(0 , 3 \varepsilon_{0} |\xi |)} & |w |^{-2+2\delta } |\xi|^{2-2\delta } 
     \absl (h_{\xi}^{\ast}\hat{\sigma}^{\infty}_{k})(w,\xi)\absr^{2}_{|d z |^{2},h} |\xi |^{-2} |\d w|^{2} \\
     & \leq K |\xi |^{2-2\delta} \int_{\C } \absl (h_{\xi}^{\ast}\hat{\sigma}^{\infty}_{k})(w,\xi)\absr^{2}_{|d z |^{2},h} |\xi|^{-2} |\d w|^{2}
\end{align*}
Recall from (\ref{confspinor}) that in the coordinate $w$ the spinors
$|\xi|^{-\alpha^{j}_{k}}h_{\xi}^{\ast}\hat{\sigma}^{\infty}_{k}$ are independent of $\xi$. 
Therefore this boils down to 
\begin{align} \label{perturbest}
       \int_{\Delta(0 , 3 \varepsilon_{0} |\xi |)} & |w |^{-2+2\delta } 
       \absl (h_{\xi}^{\ast}\hat{\sigma}^{\infty}_{k})(w)\absr^{2}_{|d z|^{2}} |\d w |^{2}\notag \\
         & \leq  K \int_{\C } \absl (h_{\xi}^{\ast}\hat{\sigma}^{\infty}_{k})(w)\absr^{2}_{|d z |^{2}} |\d w |^{2}
\end{align}
for a suitable constant $K>0$. 
Because 
$$
    (h_{\xi}^{\ast}\hat{\sigma}^{\infty}_{k})(w) \in H^{1}(\C),
$$ 
in particular we have 
$$
    (h_{\xi}^{\ast}\hat{\sigma}^{\infty}_{k})(w) \in L^{2}(\C),
$$ 
and also
$$
   \frac{1}{w} (h_{\xi}^{\ast}\hat{\sigma}^{\infty}_{k})(w) \in L^{2}_{loc}.
$$
near the origin. This implies $|w|^{-1+\delta}(h_{\xi}^{\ast}\hat{\sigma}^{\infty}_{k})(w) \in L^{2}(\C)$. Therefore, 
$$
   K = 2\frac{\left\Vert |w|^{-1+\delta}(h_{\xi}^{\ast}\hat{\sigma}^{\infty}_{k})(w) \right\Vert_{L^{2}(\C)}^{2}}
   {\left\Vert (h_{\xi}^{\ast}\hat{\sigma}^{\infty}_{k})(w) \right\Vert_{L^{2}(\C)}^{2}}
$$
has the desired property.
\end{proof}
The lemma has the following consequence. 
\begin{lem} \label{asymptnormharm}
As $| \xi | \ra \infty $, we have the estimate 
$$
     \left\vert \left\Vert \hat{\sigma}(\xi) \right\Vert_{L^{2}}^{2} - 
     \left\Vert \hat{\pi}^{H}_{\xi} \hat{\sigma}(\xi) \right\Vert_{L^{2}}^{2} \right\vert \leq 
     K |\xi|^{-2\delta} \left\Vert \hat{\sigma}(\xi) \right\Vert_{L^{2}}^{2}
$$
with $K>0$ independent of $\xi$.
\end{lem}
\begin{proof}
It is sufficient to bound 
$$
     \left\Vert \hat{\sigma}(\xi) -\hat{\pi}^{H}_{\xi} \hat{\sigma}(\xi) \right\Vert_{L^{2}}^{2} 
$$
as in the lemma. The $\Dir_{\xi }^{\ast }$-harmonic representative $\hat{\pi}^{H}_{\xi} \hat{\sigma}(\xi)$ of 
$\hat{\sigma}(\xi)$ is given by the formula
$$
      (\Id - \Dir_{\xi } G_{\xi }\Dir_{\xi }^{\ast }) \hat{\sigma}(\xi)  ,  
$$
so the difference with $\hat{\sigma}(\xi)$ itself is 
$$
     \Dir_{\xi } G_{\xi }\Dir_{\xi }^{\ast } \hat{\sigma}(\xi). 
$$ 
Since for any positive spinor $\varphi $ the estimation 
$$
    \left\Vert \Dir_{\xi } \varphi \right\Vert_{L^{2}(\C )}^{2} \leq 
    K \left\Vert \varphi \right\Vert_{H^{1}(\C )}^{2} + K |\xi |^{2} \left\Vert \varphi \right\Vert_{L^{2}(\C )}^{2} 
$$
holds, we deduce that 
\begin{equation*}
     \left\Vert \Dir_{\xi } G_{\xi }\Dir_{\xi }^{\ast } \hat{\sigma}(\xi) \right\Vert_{L^{2}(\C )}^{2} \leq 
     K  \left\Vert G_{\xi }\Dir_{\xi }^{\ast } \hat{\sigma}(\xi) \right\Vert_{H^{1}(\C )}^{2} + 
     K |\xi |^{2} \left\Vert G_{\xi }\Dir_{\xi }^{\ast } \hat{\sigma}(\xi) \right\Vert_{L^{2}(\C )}^{2} .
\end{equation*}
Lemma \ref{estgreen} implies that both terms on the right-hand side are bounded from above by 
\begin{equation*}
     K  |\xi |^{-2} \left\Vert \Dir_{\xi }^{\ast } \hat{\sigma}(\xi) \right\Vert_{L^{2}(\C )}^{2}. 
\end{equation*}
We conclude by Lemma \ref{lemstep4}. 
\end{proof}

We can now finish the proof of Theorem \ref{parweightthminf}: as $|\xi|$ goes to infinity, 
by Lemma \ref{asymptnormharm}, we have 
$$
      \frac{\left\Vert \hat{\pi}^{H}_{\xi} \hat{\sigma}(\xi) \right\Vert_{L^{2}}^{2}}
      {\left\Vert \hat{\sigma}(\xi) \right\Vert_{L^{2}}^{2}} \lra 1. 
$$
In words, the norm of the harmonic representative of the spinor $\hat{\sigma}(z,\xi)$ 
is asymptotically equal to the norm of $\hat{\sigma}(z,\xi)$ itself. On the other hand, as it has already been 
remarked in the proof of Lemma \ref{lemstep4}, we have 
$$
      \frac{\Vert \hat{\sigma}(z,\xi) \Vert_{L^{2}(\C,h)}^{2}}{\Vert \hat{\sigma}^{\infty}_{k}(z,\xi) \Vert_{L^{2}(\C,h)}^{2}}
      \lra 1
$$
exponentially as $\xi \ra \infty$. Finally, by 
(\ref{diagharm}) the $L^{2}$-norm of the spinors $\hat{\sigma}^{\infty}_{k}(z,\xi)$ as measured by the harmonic 
metric $h$ satisfy 
\begin{equation}
      c|\xi  |^{2-2\alpha^{j}_{k}} \leq \Vert \hat{\sigma}^{\infty}_{k}(z,\xi) \Vert_{L^{2}}^{2} \leq C |\xi  |^{2-2\alpha^{j}_{k}} 
\end{equation}
for some $0<c<C$, where $\alpha^{j}_{k}$ is a parabolic weight of the original Higgs bundle at the point $p_{j}$. 
All this then implies (\ref{asymbehharmspin2}), so it follows that the parabolic weight of the transformed Higgs bundle 
on the given component is equal to $\alpha^{j}_{k}-1$, as it was stated in the theorem. 
\end{proof}

\subsection{The case of logarithmic singularities}
Next we compute the parabolic weights at a puncture $\xi_{l}$
corresponding to the extension of the holomorphic structure of $\hat{\E}$ given in Subsection \ref{logext}. 

Explicitly, here is the result we wish to show.
\begin{thm} \label{parweightthmlog}
The parabolic weight of the extension $^{i}\hat{\E}$ of the transformed Higgs bundle at the puncture $\xi_{l}$, 
restricted to the $-\lambda_{k}^{\infty }$-eigenspace of the residue of the transformed Higgs field (here 
$k \in \{ 1+a_{l},\ldots,a_{l+1 } \}$) is equal to $-1+\alpha^{\infty }_{k}$, where $\alpha^{\infty }_{k}$ is the parabolic weight of 
the original Higgs field at infinity, restricted to the $\xi_{l}$-eigenspace of the second-order term and the 
$\lambda_{k}^{\infty }$-eigenspace of the first-order term of the polar part of the Higgs field. 
\end{thm}
\begin{proof}
We follow the proof of Theorem \ref{parweightthminf}. 
Again, we divide the proof into two steps according to the number of distinct eigenvalues $\xi_{l}$ of 
the second order term of $D$ at infinity. Recall that some of the spectral points $q_{k}\in \Sigma_{\xi}$ 
converge to infinity as $\xi\ra \xi_{l}$, whereas others remain bounded. 
\paragraph{Step 1.} First we suppose that $n'=1$, that is to say $A$ is a simple diagonal matrix, 
and that in a global holomorphic basis $\{\sigma^{\infty}_{k}\}$ the Higgs field has is of the form 
\[
     \theta = \frac{\xi_{1}}{2}  \d z + diag(\lambda^{\infty}_{k})\frac{\d z}{z}
\]
with one regular singularity in $0$ and an irregular one at infinity, and finally the harmonic metric is 
\begin{equation}\label{modelherminf}
      h^{\infty}=diag(|z|^{-2\alpha^{\infty}_{k}})_{k=1}^{r}.
\end{equation}
This induces a parabolic structure on $\E$ with weights $-2\alpha^{\infty}_{k}$ at $0$ and $2\alpha^{\infty}_{k}$ at infinity. 
The deformed field is
\[
     \theta_{\xi} = \frac{\xi_{1}-\xi}{2} \d z + diag(\lambda^{\infty}_{k})\frac{\d z}{z},
\]
and the spectral points are 
\[
          \frac{2\lambda_{k}}{\xi - \xi_{1}}.
\]
Making the coordinate change 
\begin{align}\label{homothlog}
      h_{\xi }: \C & \lra \C \notag \\
              w & \mapsto z = \frac{w}{\xi-\xi_{1}}
\end{align} 
the field writes 
\begin{equation}\label{confhiggslog}
     \theta_{\xi} = -\frac{1}{2} \d w  + diag(\lambda^{\infty}_{k})\frac{\d w}{w}. 
\end{equation}
The Euclidean metric $|\d z|^{2}$ on the base and the fiber metric $h^{\infty}$ are transformed into 
\begin{align}
      |\xi-\xi_{1}|^{-2}&|\d w|^{2} \label{confeucllog} \\
      diag(|\xi - \xi_{1}|^{2\alpha^{\infty}_{k}}&|w|^{-2\alpha^{\infty}_{k}})_{k=1}^{r}\label{confhermlog}
\end{align}
and the position of the spectral points become simply
\[
     2\lambda_{k}, 
\]
independent of $\xi$. As in the case of the singularity at infinity, writing $h^{(w)}$ for the diagonal 
model metric 
$$
      diag(|w|^{-2\alpha^{\infty}_{k}})_{k=1}^{r}
$$ 
the coordinate changes induce tautological isomorphisms of Hermitian fiber bundles 
\begin{align} \label{confisomlog}
     ( h_{\xi }^{\ast } E,  h^{(w)}) & \lra  (E,  h^{\infty}) \\
     (h_{\xi }^{\ast } \sigma_{k})(w ) & \mapsto |\xi - \xi_{1}|^{-\alpha_{k}} \sigma_{k}(z) \notag .
\end{align}
Via this isomorphism the representatives $v_{k}(z,\xi)$ given in (\ref{vlog}) behave as follows: 
\[
    |\xi - \xi_{1}|^{-\alpha_{k}}v_{k}(z,\xi)=v_{k}(w),
\]
which is independent of $\xi$, or equivalently 
\[
    |\xi - \xi_{1}|^{-\alpha_{k}}v_{k}(z,\xi)(\xi - \xi_{1})\d z=v_{k}(w)\d w,
\]
independent of $\xi$. By the equation (\ref{tlog}), this implies 
\[
    |\xi - \xi_{1}|^{-\alpha_{k}}t_{k}(z,\xi)(\bar{\xi} - \bar{\xi}_{1})\d \bar{z}=t_{k}(w)\d \bar{w},
\]
independently of $\xi$. Exactly as in the case of the singularity at infinity, the Laplacian and the Green's
operator of $\Dir_{\xi}^{\ast}$ in the coordinate $w$ only depend on $\xi$ through a conformal factor 
$|\xi-\xi_{1}|^{-2}$ and $|\xi-\xi_{1}|^{2}$ respectively, so the pull-back $h_{\xi }^{\ast }\hat{\pi}_{\xi}$ of the projection 
onto $\Dir_{\xi}^{\ast}$-harmonic spinors is independent of $\xi$. We deduce using invariance of the $L^{2}$-norm of
$1$-forms by conformal coordinate change that for the $\Dir_{\xi}^{\ast}$-harmonic spinor $\hat{\sigma}_{k}(z,\xi)$ 
representing the cohomology class of $(v_{k}(z,\xi)\d z,t_{k}(z,\xi)\d \bar{z})$ we have 
\[
     \int_{\C}\left\vert \hat{\sigma}_{k}(z,\xi) \right\vert_{h^{\infty},|dz|^{2}}^{2} |\d z|^{2} = 
     |\xi - \xi_{1}|^{2\alpha_{k}-2} \int_{\C}\left\vert \hat{\sigma}_{k}(w) \right\vert_{h^{(w)},|dw|^{2}}^{2} |\d w|^{2}, 
\]
where $\hat{\sigma}_{k}(w)$ is the harmonic spinor representing $(v_{k}(w)\d w,t_{k}(w)\d \bar{w})$. We see 
also that the integral on the right-hand side is independent of $\xi$, hence we have the desired behavior 
giving parabolic weight $-1+\alpha_{k}$ on this component. 
\paragraph{Step 2.} 
We drop the assumption that the second-order term $A$ of the original Higgs field is a simple matrix. 
Let $\chi$ be a fixed cut-off function supported on the complementary $\C \setminus \Delta(0,1/ \varepsilon_{0})$ of a
large disk, equal to $1$ on $\C \setminus \Delta(0,2/ \varepsilon_{0})$. In $\C \setminus \Delta(0,1/ \varepsilon_{0})$,
the Higgs field is up to a perturbation
\[
    \theta^{\infty}=\frac{1}{2}A\d z + C \frac{\d z}{z}
\]
with $A$ and $C$ diagonal matrices as in (\ref{modelhiggsinf}), therefore decomposes into a direct sum of 
problems studied in Step 1. 
In particular, for each such model problem with eigenvalue of the second-order term $\xi_{l}$ we have 
harmonic spinors $\hat{\sigma}^{l}_{k}(z,\xi)$ where $k \in \{1+a_{l},\ldots,a_{l+1} \}$, such that 
$$
    \int_{\C}\left\vert\hat{\sigma}^{l}_{k}(z,\xi)\right\vert_{|dz|^{2},h^{\infty}} |\d z|^{2} = |\xi-\xi_{l}|^{-2+2\alpha^{\infty}_{k}}.
$$ 
Again, since the harmonic metric $h$ of the Higgs bundle $(\E,\theta)$ is mutually bounded 
with $h^{\infty}$ in a neighborhood of infinity and $\hat{\sigma}^{l}_{k}$ is supported there, this implies 
\begin{equation}\label{cKineqlog}
    c|\xi-\xi_{l}|^{-2+2\alpha^{\infty}_{k}} \leq \int_{\C}\left\vert\hat{\sigma}^{l}_{k}(z,\xi)\right\vert_{|dz|^{2},h} |\d z|^{2}
    \leq C |\xi-\xi_{l}|^{-2+2\alpha^{\infty}_{k}}
\end{equation}
for some $0<c<C$. The section 
\[
     \hat{\sigma}(z,\xi)= \chi(z) \hat{\sigma}^{l}_{k}(z,\xi)
\]
is well-defined because the local holomorphic trivialisation $\sigma^{\infty}_{k}$ of $\E$ is defined in 
$\C \setminus \Delta(0,1/ \varepsilon_{0})$ for $\varepsilon_{0}>0$ sufficiently small. 
The statement of the theorem will again follow if we prove 
\begin{equation}\label{asymbehharmlog}
     c|\xi-\xi_{l}|^{-2+2\alpha^{\infty}_{k}} \leq  \int_{\C}\left\vert \hat{\pi}_{\xi^{H}}\hat{\sigma}(z,\xi)\right\vert_{|dz|^{2},h} |\d z|^{2}
    \leq C |\xi-\xi_{l}|^{-2+2\alpha^{\infty}_{k}}
\end{equation}
where $\hat{\pi}_{\xi}^{H}\hat{\sigma}(z,\xi)$ is the harmonic representative of $\hat{\sigma}(z,\xi)$. 
As a first step in this direction, we prove: 
\begin{lem}\label{lemstep4log}
There exists $\delta>0$ and $K>0$ such that for $|\xi|$ sufficiently large the inequality 
$$
   \left\Vert \Dir_{\xi }^{\ast } \hat{\sigma}(z,\xi) \right\Vert_{L^{2}(\C )}^{2} 
   \leq K |\xi -\xi_{l}|^{2+2 \delta} \left\Vert  \hat{\sigma}(z,\xi) \right\Vert_{L^{2}(\C )}^{2}
$$
holds.
\end{lem}
\begin{proof}
We follow the proof of Lemma \ref{lemstep4}. We set $(D''_{\xi})^{\infty}=\dbar^{\E}+\theta^{\infty}$ and let $\Dir_{\xi}^{\infty}$ 
(respectively $(\Dir_{\xi}^{\infty})^{\ast}$) stand for its Dirac operator (respectively its adjoint).
By Lemma \ref{expdecord}, $\hat{\sigma}^{l}_{k}$ is supported in $L^{2}$-norm up to an exponentially decreasing factor
in $\xi$ in  $\C \setminus \Delta(0,1/ \varepsilon_{0})$. Therefore, the lemma reduces to the same estimation for
$\hat{\sigma}^{l}_{k}$. Moreover, by assumption we have 
\[
    (\Dir_{\xi}^{\infty})^{\ast}\hat{\sigma}^{l}_{k}(z,\xi)=0,
\]
so 
\[
    \Dir_{\xi }^{\ast } \hat{\sigma}^{l}_{k}(z,\xi) = [\Dir_{\xi }^{\ast } - (\Dir_{\xi}^{\infty})^{\ast}] \hat{\sigma}^{l}_{k}(z,\xi). 
\]
The difference on the right-hand side of this formula is bounded above by $K|z|^{-1-\delta}$ for some $K>0$
independent of $\xi$, because the two Dirac operators depend on $\xi$ in the same way, hence their difference 
does not depend on it at all. Introducing the coordinate $w=z(\xi -\xi_{l})$, this becomes 
$K|w|^{-1-\delta}|\xi-\xi_{l}|^{1+\delta}$. Therefore, it is sufficient to prove 
\begin{align*}
     \int_{\C \setminus \Delta(0,|\xi -\xi_{l}|/ \varepsilon_{0})}&
     |w|^{-2-2\delta}|\xi -\xi_{l}|^{2+2\delta} \absl \hat{\sigma}^{l}_{k}(z,\xi)\absr_{|dz|^{2},h}^{2} |\xi -\xi_{l}|^{-2} |\d w|^{2} \\ \leq 
     & K |\xi -\xi_{l}|^{2+2\delta} \int_{\C} \absl \hat{\sigma}^{l}_{k}(z,\xi)\absr_{|dz|^{2},h}^{2} |\xi -\xi_{l}|^{-2} |\d w|^{2},\notag
\end{align*}
for a suitable $K>0$, or more simply
\begin{align}\label{ineqxideltalog}
     \int_{\C \setminus \Delta(0,|\xi -\xi_{l}|/ \varepsilon_{0})}|w|^{-2-2\delta} &
     \absl \hat{\sigma}^{l}_{k}(z,\xi)\absr_{|dz|^{2},h}^{2} |\d w|^{2} \notag \\ 
     & \leq K \int_{\C} \absl \hat{\sigma}^{l}_{k}(z,\xi)\absr_{|dz|^{2},h}^{2} |\d w|^{2}.
\end{align}
This goes similarly to (\ref{perturbest}): because in the coordinate $w=h_{\xi}^{-1}z$ the spinor 
$|\xi -\xi_{l}|^{2-2\alpha^{\infty}_{k}}\hat{\sigma}^{l}_{k}(z,\xi)$ is independent of $\xi$ (see Step 1) and $h$ and $h^{\infty}$ are mutually
bounded, it boils down to  
\begin{align*}
    \int_{\C \setminus \Delta(0,|\xi -\xi_{l}|/ \varepsilon_{0})}|w|^{-2-2\delta} &
    \absl (h_{\xi}^{\ast}\hat{\sigma}^{l}_{k})(w)\absr_{|dz|^{2},h^{\infty}}^{2} |\d w|^{2} \\
    & \leq K \int_{\C} \absl (h_{\xi}^{\ast}\hat{\sigma}^{l}_{k})(w)\absr_{|dz|^{2},h^{\infty}}^{2} |\d w|^{2}.
\end{align*}
Now remark that $h_{\xi}^{\ast}\hat{\sigma}^{l}_{k}\in H^{1}(\C,|\d w|^{2},h^{\infty})$ implies in particular that 
$h_{\xi}^{\ast}\hat{\sigma}^{l}_{k}\in L^{2}(\C,|\d w|^{2},h^{\infty})$. Furthermore, near the origin 
$|w|^{-1-\delta} h_{\xi}^{\ast}\hat{\sigma}^{l}_{k}\in L^{2}_{loc}(|\d w|^{2},h^{\infty})$ provided that $\delta<\alpha^{\infty}_{k}$. 
Hence $|w|^{-1-\delta} h_{\xi}^{\ast}\hat{\sigma}^{l}_{k}\in L^{2}(\C,|\d w|^{2},h^{\infty})$, and 
$$
    K = 2 \frac{\left\Vert|w|^{-1-\delta} h_{\xi}^{\ast}\hat{\sigma}^{l}_{k}\right\Vert_{L^{2}(\C,|d w|^{2},h^{\infty})}^{2}}
    {\left\Vert h_{\xi}^{\ast}\hat{\sigma}^{l}_{k}\right\Vert_{L^{2}(\C,|d w|^{2},h^{\infty})}^{2}}
$$
has the desired property (\ref{ineqxideltalog}). 
\end{proof}
This has the following consequence.
\begin{lem}\label{asymptnormharmlog}
As $\xi \ra \xi_{l}$, we have the estimate
\[
     \left\vert \left\Vert \hat{\sigma}(z,\xi) \right\Vert_{L^{2}}^{2} - 
     \left\Vert \hat{\pi}^{H}_{\xi}  \hat{\sigma}(z,\xi)\right\Vert_{L^{2}}^{2} \right\vert \leq 
     K |\xi-\xi_{l}|^{2\delta} \left\Vert \hat{\sigma}(z,\xi) \right\Vert_{L^{2}}^{2}
\]
for some $K>0$ independent of $\xi$.
\end{lem}
\begin{proof}
Again as in Lemma \ref{asymptnormharm}, it is sufficient to bound 
$$
     \left\Vert \hat{\sigma}(z,\xi) -\hat{\pi}^{H}_{\xi} \hat{\sigma}(z,\xi) \right\Vert_{L^{2}}^{2} 
$$
as in the lemma, where 
$$
      \hat{\pi}^{H}_{\xi} \hat{\sigma}(z,\xi)=(\Id - \Dir_{\xi } G_{\xi }\Dir_{\xi }^{\ast }) \hat{\sigma}(z,\xi) 
$$
is the $\Dir_{\xi }^{\ast }$-harmonic representative of $\hat{\sigma}(\xi)$. Thus by Lemma \ref{estgreenlog} 
we have for the norm of the difference 
$$
    \left\Vert \Dir_{\xi} G_{\xi}\Dir_{\xi}^{\ast }\hat{\sigma}(z,\xi) \right\Vert_{L^{2}}^{2}
    \leq K|\xi-\xi_{l}|^{-2} \left\Vert \Dir_{\xi }^{\ast }\hat{\sigma}(z,\xi) \right\Vert_{L^{2}}^{2}
$$
and we conclude using Lemma \ref{lemstep4log}.
\end{proof}

We are now ready to finish the proof of Theorem \ref{parweightthmlog}: by Lemma \ref{asymptnormharmlog}, as 
$\xi\to\xi_{l}$ the norm of the harmonic representative of the spinor $\hat{\sigma}(z,\xi)$ verifies 
$$
    \frac{\left\Vert \hat{\pi}^{H}_{\xi} \hat{\sigma}(\xi) \right\Vert_{L^{2}}^{2}}
      {\left\Vert \hat{\sigma}(\xi) \right\Vert_{L^{2}}^{2}} \lra 1.
$$ 
On the other hand, since the support of $\chi$ in the coordinate $w$ is 
$\C \setminus \Delta(0,|\xi -\xi_{l}|/ \varepsilon_{0})$, and these sets exhaust $\C$ as $\xi\to\xi_{l}$, we have that 
$$
    \frac{\left\Vert \hat{\sigma}(\xi) \right\Vert_{L^{2}}^{2}}{\left\Vert \hat{\sigma}^{l}_{k}(\xi) \right\Vert_{L^{2}}^{2}}\lra 1.
$$
By (\ref{cKineqlog}) the $L^{2}$-norm of $\hat{\sigma}^{l}_{k}(z,\xi)$ as measured by
the harmonic metric $h$ satisfies 
$$
     c|\xi-\xi_{l}|^{-2+2\alpha^{\infty}_{k}} \leq \int_{\C}\left\vert\hat{\sigma}^{l}_{k}(z,\xi)\right\vert_{|dz|^{2},h}^{2} |\d z|^{2}
    \leq C |\xi-\xi_{l}|^{-2+2\alpha^{\infty}_{k}}.
$$
Putting together all this, we obtain (\ref{asymbehharmlog}), 
so that on the component of $\hat{E}$ near $\xi_{l}$ on which the transformed Higgs field has eigenvalue 
$-\lambda^{\infty}_{k}$, the parabolic weight of the induced extension is $-1+\alpha^{\infty}_{k}$. 
\end{proof}

\section{The topology of the transformed bundle} \label{secttop} 
In this section, we compute the topology of the underlying holomorphic
bundle $^{i}\hat{\E}$ of the transformed Higgs bundle (see (\ref{defholtrafo})) relative to its 
extension over the punctures given in Section \ref{ext}. We then deduce the topology 
of the transformed Higgs bundle relative to its transformed extension given by Definition \ref{transfext}. 
We recall that we have denoted
\begin{equation} 
    \hat{r} = \sum_{p \in P} rk(Res(\theta , p))). 
\end{equation} 
The result we wish to show is the following: 
\begin{thm} \label{thmtop}
The rank of $^{i}\hat{\E}$ is equal to $\hat{r}$, whereas its degree is equal to 
$\hat{r}+ deg(\E ) + r$, where $r$ and $deg(\E ) $ are the rank and degree of $\E$, respectively. 
\end{thm}
Notice that it gives in particular (\ref{i}) of Theorem \ref{mainthmhiggs}. 
\begin{proof}
Recall that we have denoted by $\E$ the sheaf of holomorphic
sections of the bundle $\E$ underlying the original Higgs bundle; $\F$ was defined as a sheaf of
meromorphic sections of $\E \otimes \Omega^{1,0}$ having singularities at $P \cup \{ \infty \}$
with singular parts in prescribed spaces (see Subsection \ref{sheaves});
and finally $\tilde{\F}=\pi_{1}^{\ast }\F \otimes \pi_{2}^{\ast } \mathcal{O}_{\CPt}(1)$. 
By hypothesis, $\theta$ (and so $\theta_{\eta } $ for any $\eta $) is holomorphic with respect to the holomorphic
structure $\dbar^{\E}$. Thus we may consider the holomorphic chain complex 
\[ 
\xymatrix{ 
              \E \ar[r] \ar[d]_{\Id} & 0 \ar[d] \\
              \E \ar[r]^{\theta_{\eta }} \ar[d] & \tilde{\F} \ar[d]^{\Id} \\
              0 \ar[r] & \tilde{\F}
}
\]
in $\eta \in \widehat{\mathbf{CP}}^1 $. The hypercohomology long exact sequence
associated to it yields the exact sequence of cohomology spaces
\begin{align} \label{les}
     0 \lra H^{0}({\CP} , \E)  & \xra{\theta_{\eta}} H^{0}({\CP} , \tilde{\F})
     \lra \H^{1}( \E \xra{\theta_{\eta }}\tilde{\F})\notag \\ 
     & \lra H^{1}({\CP},\E) \xra{\theta_{\eta}} H^{1}({\CP} , \F) \lra 0,  
\end{align}
since we have seen that $\H^{0}(\E \xra{\theta_{\eta }}\tilde{\F} ) = \H^{2}(\E
\xra{\theta_{\eta }}\tilde{\F} )=0$. 
All of the spaces in this exact sequence come with a natural holomorphic
structure over $\CPt$: 
\begin{itemize}
\item{the cohomology spaces of $\E$ because this latter is
trivial over $\CPt$}
\item{those of $\tilde{\F}$ because this latter is the tensor
product of a trivial vector bundle over $\CPt$ and $\mathcal{O}_{\CPt}(1)$}
\item{finally, $\H^{1}( \E \xra{\theta_{\bullet }}\tilde{\F} )=\hat{V}_{\bullet }$ has its
    holomorphic structure $\dbar^{\hat{\E}}$ induced by $\hat{\d}^{0,1}$, extended to 
    the singularities in Section \ref{ext} by the induced extension $^{i}\hat{\E}$.}
\end{itemize}
Moreover, all of the maps in the 
exact sequence (\ref{les}) vary holomorphically in $\eta \in \CPt$ with respect to
these structures and extensions: this follows from the definition of $\tilde{\F}$ and that of the induced
extension. Therefore, it induces an exact sequence of the sheaves over ${\CPt}$ 
of holomorphic sections of the corresponding cohomology spaces: 
\begin{align*}
     0 \lra  \mathcal{O}( H^{0}( \E)) & \xra{\theta_{\eta }} \mathcal{O}(H^{0}( \tilde{\F})) \lra \mathcal{O}(^{i}\hat{\E}) \\
      & \lra  \mathcal{O}(H^{1}(\E)) \lra \mathcal{O}(H^{1}( \F)) \lra 0, 
\end{align*}
where $\mathcal{O} $ stands to denote the sheaf of regular sections on $\CPt$ with respect
to the above mentioned holomorphic structures. 
By additivity of the Chern character, we deduce the equality 
\begin{align}
        ch(^{i}\hat{\E})= & ch(\mathcal{O}(\widehat{\mathbf{CP}}^1, H^{0}( \tilde{\F}))) -
        ch(\mathcal{O}(\widehat{\mathbf{CP}}^1, H^{1}( \tilde{\F}))) \label{chern1} \\
        & - ch(\mathcal{O}(\widehat{\mathbf{CP}}^1, H^{0}( \E))) + ch(\mathcal{O}(\widehat{\mathbf{CP}}^1, H^{1}( \E)))
        \label{chern2} 
\end{align}
in $H^{\ast }(\widehat{\mathbf{CP}}^1)$. Put $\pi = \pi_{2}$, the projection onto the second
factor in $\CP \times \widehat{\mathbf{CP}}^1$. One has direct image sheaves $\pi_{\ast }\E$ and 
$\pi_{\ast }\tilde{\F}$ on $\widehat{\mathbf{CP}}^1$ defined by 
\begin{align*}
   \pi_{\ast }\E |_{U} & = \mathcal{O}(U,H^{0}(\CP, \E )) \\
   \pi_{\ast }\tilde{\F}|_{U} & =   \mathcal{O}(U,H^{0}(\CP ,\tilde{\F})))=\mathcal{O}(U,H^{0}(\CP ,{\F})) \otimes
   \mathcal{O}_{\widehat{\mathbf{CP}}^1}(1)(U), 
\end{align*}
for any open set $U \in \CP$, and one can form the ``virtual'' sheaves 
\begin{align*}
        \pi_{!}\E|_{U} & =  \mathcal{O}(U,H^{0}(\CP , \E)) - \mathcal{O}(U,H^{1}(\CP , \E)) \\
        \pi_{!} \tilde{\F}|_{U} & =  \mathcal{O}(U,H^{0}(\CP , \tilde{\F})) - \mathcal{O}(U,H^{1}(\CP , \tilde{\F})).
\end{align*}
Again by additivity of the Chern character, the right-hand-side of
(\ref{chern1}) is equal to $ch(\pi_{!} \tilde{\F})$, which is in turn equal
to 
$$
    \pi_{\ast }(ch(\tilde{\F}) \cup Td(T_{\pi }) ), 
$$
by the Grothendieck-Riemann-Roch theorem, where 
$$
     T_{\pi } = T(\CP \times \widehat{\mathbf{CP}}^1) - \pi^{\ast }T\widehat{\mathbf{CP}}^1 = \pi_{1}^{\ast}T\CP
$$
is the relative tangent bundle  of $\pi $, and $Td$ stands for its Todd class. Moreover, $\pi_{\ast} $ is just evaluation
on the fundamental cycle of $\CP$. Similarly, we see that (\ref{chern2}) is just 
$$
   - ch(\pi_{!} \E ) = - \pi_{\ast }(ch(\E) \cup Td(T_{\pi }) ), 
$$
and thus we obtain 
\begin{equation} \label{grr}
     ch(^{i}\hat{\E})= [ ( ch(\tilde{\F})-ch(\E) ) \cup Td(\pi_{1}^{\ast} T \CP )  ]/[\CP ]. 
\end{equation}
Now we have 
\begin{align*}
       ch(\E ) & = r + c_{1}(\E) \\
       ch(\tilde{\F})  & = \left[ r + c_{1}(\E ) + h \sum_{p \in P} rk(Res(\theta ,
         p)) \right] (1+ \hat{h}) \\
       Td(T\CP) & = Td(\mathcal{O}_{\CP}(2))= 1+ h,  
\end{align*}
where $r$ is the rank of the bundle $\E$, $c_{1}(\E)$ is its first Chern
class, and $h$ and $\hat{h}$ are the hyper-plane classes of $\CP$ and
$\widehat{\mathbf{CP}}^1$ respectively. Putting all this together, we obtain 
\[ 
    ch(\tilde{\F}) - ch(\E ) =  \hat{r} h + [r+ c_{1}(\E)+ \hat{r}] \hat{h},
\]
and plugging this into (\ref{grr}), 
\begin{equation} \label{top}
       ch(^{i}\hat{\E}) = \hat{r}+ [r+ deg(\E)+ \hat{r}] \hat{h}, 
\end{equation}
as we wished. 
\end{proof}

We are now ready to pass back to the transformed extension of the Higgs bundle introduced in Definition
\ref{transfext}, hence establishing points (\ref{ii}), (\ref{v}) and (\ref{viii}) of Theorem
\ref{mainthmhiggs}. 

\begin{cor} \label{pardegcor} 
The parabolic weights of the transformed Higgs bundle endowed with its transformed extension are 
$\alpha_{k}^{\infty }$ at the logarithmic punctures (on the same subspace as in Theorem \ref{parweightthmlog}) 
and $\alpha_{k}^{j }$ at infinity (on the subspace in Theorem \ref{parweightthminf}). 
The degree of the transformed Higgs bundle $\hat{\E}$ with respect to its transformed extension is equal to 
the degree of ${\E}$. 
\end{cor}
\begin{proof}
Recall from Theorems \ref{parweightthmlog} and \ref{parweightthminf} that the parabolic weigths of 
the transformed Higgs bundle relative to the induced extensions considered in Subsections \ref{logext} and \ref{infext} 
are equal to $-1+\alpha_{k}^{\infty}$ at the logarithmic punctures and to $-1+\alpha_{k}^{j}$ at infinity. 
On the other hand, by Definition \ref{transfext}, the parabolic weights of the transformed Higgs bundle
with respect to its transformed extension are required to have parabolic weights between $0$ and $1$.  
This means that a local holomorphic trivialisation of the singular component of the transformed extension 
$\hat{\E}$ near the puncture $\xi_{l}$ is
\[
      (\xi -\xi_{l}) \hat{\sigma}^{l}_{k}(\xi ) , 
\]
where $\hat{\sigma}^{l}_{k}(\xi )$ is the local holomorphic section of the extension $^{i}\hat{\E}$ at $\xi_{l}$
defined in Subsection \ref{logext} and $k \in \{1+a_{l},\ldots,a_{l+1}\}$. On the regular component of
$\hat{\E}|_{\xi_{l}}$ the harmonic representatives have bounded norm, which gives $0$ parabolic weight. 
Therefore on this component one does not need to change the trivialisation. Similarly, a local holomorphic 
frame of $\hat{\E}$ near infinity can be expressed by 
\[
      \xi^{-1}  \hat{\sigma}_{k }^{\infty}(\xi) , 
\]
where $\hat{\sigma}^{\infty}_{k}$ is the local holomorphic section of the extension $^{i}\hat{\E}$ at infinity 
defined in Subsection \ref{infext} localized near $p_{j}$ for some $j\in \{ 1,\ldots,n\}$, and $k \in \{r_{j}+1,\ldots,r\}$. 
Clearly, this way we increased all non-vanishing parabolic weights by $1$. On the other hand, 
by Remark \ref{pardegvanish} even if the algebraic geometric degree of the bundle depends on the choice 
of extensions, the parabolic degree with respect to a fixed metric is independent of them, because it is always $0$. 
Recall from Definition \ref{parintcon} that 
\[
    deg_{par}(^{i}\hat{\E }) = deg(^{i}\hat{\E }) + \sum_{j \in \{1,\ldots,n, \infty \} } \sum_{k=r_{j}+1}^{r}  (-1+ \alpha^{j}_{k}). 
\]
This quantity is therefore equal to 
\begin{equation} \label{pardegEtr}
    deg_{par}(\hat{\E} ) = deg(\hat{\E} ) + \sum_{j \in \{1,\ldots,n, \infty \} } \sum_{k=r_{j}+1}^{r}  \alpha^{j}_{k}. 
\end{equation}
Putting these expressions together, we deduce that 
\[
    deg(\hat{\E} ) =  deg(^{i}\hat{\E }) - \hat{r} -r,
\]
where we recall again that we have defined 
\[
    \hat{r}= \sum_{j=1}^{n} rk (Res(\theta, p_{j} )).
\]
Using formula (\ref{top}) we get  
\begin{equation} \label{p-top}
       deg(\hat{\E })= deg(\E ).
\end{equation}
\end{proof}



\chapter{The inverse transform} \label{chinv}

In this chapter we construct the inverse of the transform introduced in the
previous chapters. In line with the properties of the ordinary Fourier
transform and its algebraic counterparts, the inverse is defined by a
formula which only differs from the transform in a sign. 

Recall from Section \ref{flattr} that the transformed flat connection on 
$\hat{E}_{\bullet }=L^{2}H^{1}(D_{\bullet })$ is defined by the $L^{2}$-orthogonal
projection of $\hat{\d}-z\d \xi \land$.  
For any parabolic vector bundle with integrable connection $ (F,D^{F},h^{F})$ on
$\hat{\C}$ satisfying the conditions of Section \ref{flatintro} (i.e. having
a finite number of simple poles in finite points and a second-order
pole at infinity, such that the eigenvalues and parabolic weights meet the conditions imposed 
in Theorem \ref{mainthm}), one can define the inverse transformed bundle with
integrable connection $(\check{F},\check{D}^{F},\check{h}^F)$ on $\C$ by a procedure similar to the one
defining $(\hat{E},\hat{D},\hat{h})$ starting from $(E,D,h)$: namely, consider the deformation 
\begin{equation} \label{invdeform}
     D^{F}_{z}=D^{F} + z \d \xi \land 
\end{equation}
of the connection parametrized by $z$ in $\C$ minus a finite set, 
and let $\check{F}_{z}$ be the first $L^{2}$-cohomology of 
$$ 
     F \xra{D^{F}_{z}} \Omega^{1}_{\hat{C}} \otimes F \xra{D^{F}_{z}} \Omega^{2}_{\hat{C}} \otimes F.
$$ 
These vector spaces are of the same dimension and form a
smooth vector bundle over $\C$ minus a finite number of points. The
critical points are easily seen to be the opposites of the eigenvalues of the second-order 
term of ${D}^{F}$ at infinity. The proof goes similarly to the case
of the direct transform. We also define the Hilbert bundle $\check{H}$ over
$\C $, the $L^{2}$-metric $\check{h}$ and the orthogonal projection 
$\check{\pi}_{z}: \check{H}_{z} \ra  \check{F}_{z} $ in an analogous manner as in Section
\ref{flattr}. 
Next, let the inverse transformed integrable connection
$\check{D}^{F}$ be defined by the parallel sections 
$\check{\pi}_{z}(e^{(z_{0}-z)\xi } \phi_{z_{0}} (\xi ))$ for any harmonic section 
$\phi_{z_{0}} (\xi ) \in \check{F}_{z_{0}}$. Equivalently, denoting by $\check{\d}$
the trivial connection with respect to $w$ in the trivial Hilbert bundle
$\check{H}$, the inverse transformed flat connection can be given by the 
formula 
\begin{equation} \label{flinvexpl}
        \check{D}^{F} = \check{\pi}_{z} (\check{\d} + z \d \xi ),
\end{equation}
as it can be seen by the argument given in Section \ref{flattr}, changing
signs. Finally, we define the inverse transformed metric $\check{h}^F$ on the fiber
$\check{F}_{z_{0}}$ again as the $L^{2}$-norm on $\hat{\C}$ of a ${D}^{F}_{z_{0}}$-harmonic representative. 
We can now state the 
\begin{thm} \label{invtr} 
The inverse transform of $\Nahm : (E,D,h)\mapsto (\hat{E},\hat{D}, \hat{h})$ is 
$\Nahm^{-1} :(F,D^{F},h^{F}) \mapsto (\check{F},\check{D}^{F},\check{h}^{F})$. In different terms, for any
bundle with integrable connection and harmonic metric $(E,D,h)$ satisfying
the conditions of Section \ref{flatintro} and the ones imposed in Theorem \ref{mainthm}, there exists a canonical
Hermitian bundle isomorphism $\omega $ between $\check{\hat{E}}$ and $E$ such
that $\omega^{\ast } D = \check{\hat{D}} $.
\end{thm}
\begin{rk}
As one can check using the transform on the level of singularity parameters described in Theorem
\ref{mainthm}, the assumptions (1) and (2) of that theorem are symmetric, in the sense that 
if they are fulfilled by $(E,D)$ than the same is true for $(\hat{E},\hat{D})$. Therefore, the 
transform $\check{}$ can be applied to this latter, so the affirmation of the theorem has 
a meaning. 
\end{rk}
\begin{proof}
The proof is done in four steps: 
first, we prove that the fibers over $0 \in \C$ of $E$ and $\check{\hat{E}}$ are
canonically isomorphic. Next we show the same thing for the other fibers.
Then we prove that the integrable connections are the same, and finally 
we establish equality of the harmonic metrics and parabolic structures.  

\paragraph{Step 1.} 
Consider the product manifold $\C \times \hat{\C}$, and let $\pi_{1}$ and $\pi_{2}$
be the projection to the first and second factor, respectively. 
Denote by $E$ the pull-back vector bundle $\pi_{1}^{\ast }E$ on the product, and
define the connection $\D =\pi_{1}^{\ast }D - \xi \d z - z \d \xi $. Notice that on the 
fiber $\C \times \{ \xi_{0} \}$ this just gives the deformation $D_{\xi_{0}}$. Now
form the double complex 
$$ 
    \Dou^{p,q} = \Omega_{{\C}}^{p} \otimes \Omega_{\hat{\C}}^{q}(E),
$$ 
where $\Omega_{{\C}}^{p}$ (respectively $\Omega_{\hat{\C}}^{q}$) denote smooth
$p$-forms (smooth $q$-forms) on $\C$ ($\hat{\C}$); and 
with differentials $\d_{1}=D_{\xi }, \d_{2}=\hat{\d}-z \d \xi \land $. 
Remark that these differentials commute (in the graded sense), and their 
sum is just $\D$. The desired
isomorphism will result from the study of the spectral sequences
corresponding to the two different filtrations of this double complex. 

Namely, consider the first filtration of $\Dou$: the first page of the
corresponding spectral sequence $\S_{1}^{ \bullet , \bullet  }$ is 
\begin{eqnarray} \label{D1}
      \xymatrix{ 0 & \Omega_{\hat{\C}}^{2} \otimes \hat{E} & 0 \\
                 0 & \Omega_{\hat{\C}}^{1} \otimes \hat{E} \ar[u]^{\d_{2}^{\sharp }} & 0 \\ 
                 0 & \hat{E} \ar[u]^{\d_{2}^{\sharp }} & 0 }
\end{eqnarray}
where $\d_{2}^{\sharp }$ stands for the operator induced by $\d_{2}$. More
precisely, this operator is obtained as follows. Consider for example a
local section of $\hat{E}$: if $B(\xi_{0} )$ is an open ball in $\hat{\C}$,
it is given by cohomology classes $[\phi_{\xi }]$ in
$L^{2}H^{1}(D_{\xi })$ changing smoothly with $\xi \in B(\xi_{0} )$. Here $\phi_{\xi }=\phi_{\xi }(z)$ is a global
$L^{2}$-section of  $E$ over $\C$, in the kernel of $\Dir_{\xi }^{\ast }$. In
particular, $D_{\xi } \phi_{\xi }=0$, and since the two differentials commute, we
then have $D_{\xi } \circ \d_{2}  \phi_{\xi }=0$. In other words, $\d_{2}\phi_{\xi}$
is a $\d_{1}$-closed section of $\Dou^{1,1}$ on $\C \times B(\xi_{0} )$; hence we
may consider its cohomology class with respect to $\d_{1}$, and letting $\xi$
vary these give a section of $\Omega^{1} \otimes \hat{E}$ over $B(\xi_{0} )$, which is
by definition $\d_{2}^{\sharp } [\phi_{\xi }]$. Now remark that under the
isomorphism of the first $L^{2}$-cohomology of the elliptic complex (\ref{ltwocompl})
and the space of $\Dir_{\xi }$-harmonic sections given in Theorem
\ref{laplker}, this induced connection goes over to $\hat{D}$
defined in Section \ref{flattr}; 
in other words, under these identifications $\d_{2}^{\sharp }=\hat{D}$.  
Moreover, the connection
$\hat{D}$ also satisfies the conditions of Section \ref{flatintro}. 
Therefore, by Chapter \ref{chFr} and Section \ref{l2coh}
the $L^{2}$-cohomology of $\hat{D}=\hat{D}_{0}$ is non-trivial
only in degree $1$, and so when passing to the second page $\S_{2}^{ \bullet , \bullet }$ of the
spectral sequence, we obtain by definition
$\S_{2}^{1,1}=\check{\hat{E}}_{0}$ and all other terms equal to $0$. In
particular, the spectral sequence collapses at the second page, and the
total cohomology of the double complex is canonically isomorphic to $\check{\hat{E}}_{0}$ in degree $2$
and vanishes in all other degrees. 

Consider now the second filtration of $\Dou$: in order to form the first
page $\tilde{\S}_{1}^{ \bullet , \bullet }$ of the corresponding spectral sequence, we first take
cohomology on each column of the double complex with respect to
$\d_{2}=\hat{\d} -z \d \xi $, and so it is equal to 
\begin{eqnarray} \label{D2}
     \xymatrix{ 0 & 0 & 0 \\
                0 & 0 & 0 \\ 
                 L^{2}(\C,E) e^{z\xi } \ar[r]^(.42){\d_{1}^{\sharp }} &
                 L^{2}(\C,\Omega_{\C}^{1} \otimes  E) e^{z\xi } \ar[r]^(.48){\d_{1}^{\sharp }} & 
                 L^{2}(\C,\Omega_{\C}^{2} \otimes  E) e^{z\xi }. }
\end{eqnarray}
In words: for example, the $(0,0)$-term consists of $L^{2}$-sections of $E$
on $\C \times \hat{\C}$ which are a product of an arbitrary section of $E$ on
$\C$ and the function $e^{z\xi }$. Now notice that the only
possibility for a non-zero section of this form to be in $L^{2}$ on $\{z \} \times  \hat{\C}$ is for $z=0$.
Put another way, the cohomology along the slices $\{z\} \times \hat{\C}$ vanishes
for all $z \neq 0$. Hence we may replace the double complex $\Dou$ without changing the spectral
sequence associated with this filtration (and so the total cohomology), 
by the double complex $\germ$ whose component of bidegree $(p,q)$ is the
space of $L^{2}$-forms with values in $E$ of bidegree $(p,q)$ defined on 
$V_{0} \times \hat{\C}$  for any neighborhood $V_{0}$ of $0 \in \C$, and where we
identify such forms if they coincide on an arbitrary neighborhood of 
$\{0\} \times \hat{\C}$. Of course, the differentials of this new double complex
are induced by those of $\Dou$ in a trivial way. 

The idea now is to consider the spectral sequence $\Sg$ corresponding to the
first filtration of $\germ$: by the general theory of spectral sequences,
this will then abut to the total cohomology of $\germ$, which is, as we
saw in the previous paragraph, equal to that of $\Dou$, that is to
$\check{\hat{E}}_{0}$. First trivialize $E$ in $V_{0}$: this just means
that we identify the total space of the bundle with $V_{0} \times E_{0}$. 
Since the vector bundle $E$ on $\C \times \hat{\C}$ is just the pull-back of $E$
on $\C$, this also gives an identification of $E \ra V_{0} \times \hat{\C}$ with 
the trivial bundle $ (V_{0} \times \hat{\C}) \times  E_{0}$. 
Without loss of generality we may assume $0 \notin P$, so for $V_{0}$ sufficiently small 
the connection $D$ can also be taken by a gauge transformation $\tilde{g}$ to the 
trivial one. Thus in this trivialisation and gauge we have 
$\d_{1}=\d -\xi \d z$ where $\d$ stands for the trivial connection in the $z$
direction. The first page $\Sg_{1}^{ \bullet , \bullet }$ is then equal to
the cohomology spaces with respect to this differential: 
\begin{eqnarray} 
      \label{germD1} \xymatrix{ \Omega_{\hat{\C}}^{2} \otimes L^{2}(\hat{\C},E_{0}) e^{z\xi } & 0  & 0\\
                \Omega_{\hat{\C}}^{1} \otimes L^{2}(\hat{\C},E_{0}) e^{z\xi } \ar[u]^{\d_{2}^{\sharp }} & 0 & 0 \\ 
                  L^{2}(\hat{\C},E_{0}) e^{z\xi } \ar[u]^{\d_{2}^{\sharp }} & 0 & 0, }
\end{eqnarray}
where, as before, $ L^{2}(\hat{\C},E_{0}) e^{z\xi }$ stands to denote
functions with values in $E_{0}$ of
the form $\gamma(\xi ) e^{z \xi }$ but this time on $V_{0} \times \hat{\C}$, and the
$L^{2}$ condition now only implies that $\gamma $ must be rapidly decreasing as $|\xi| \ra \infty $. 
The next remark is that since we only have terms in degree $p=0$, the
differential induced by $\d_{2}$ is just itself: indeed, it is by
definition $\d_{2}$ modulo the image of $\d_{1}$, but this latter vanishes
for $p=0$. Thus, in order to obtain the second page $\Sg_{2}^{ \bullet , \bullet}$ 
of the spectral sequence, we take cohomology with respect to
$\d_{2}=\hat{\d} -z \d \xi \land$. Notice that via the gauge transformation
$e^{-z \xi }$ the whole picture can be rephrased as the de Rham cohomology of
rapidly decreasing sections $\sigma $ on $\hat{\C}$ with values in $E_{0}$, 
which is similar to compactly supported de Rham cohomology. 
Therefore in $\Sg_{2}^{ \bullet , \bullet }$ all elements except for the one
corresponding to bidegree $(0,2)$ vanish, and this latter is canonically
isomorphic to $E_{0}$ via mapping an element $\gamma_{0} \in E_{0}$ into the germ 
$$
        [\gamma_{0} \chi(\xi ) e^{z\xi } \d \xi \land \d \bar{\xi }],
$$
where $\chi$ is a fixed exponentially decreasing bump-function on $\hat{\C}$ with
integral (with respect to the volume form $|\d \xi |^{2}$) equal to $1$, and 
$[.]$ stands to denote the de Rham cohomology class of exponentially decreasing
forms on $\hat{\C} $ with values in $E_{0}$. 
Conversely, for an arbitrary class $ [\gamma_{0}( \xi ) e^{z \xi } \d \xi \land \d \bar{\xi }]$ where
$\gamma_{0}( \xi ) e^{z \xi }$ is a germ of exponentially decreasing functions on
$\hat{\C}$ with values in $E_{0}$ and in the kernel of $\d_{1}=(\d -\xi \d z)$, we may define 
\begin{align}\label{explisom}
       [\gamma_{0}( \xi )  e^{z \xi } \d \xi \land \d \bar{\xi }] \mapsto & 
       eval_{z=0}\int_{\hat{\C}} \gamma_{0}( \xi )  e^{z \xi } |\d \xi |^{2} \notag \\ 
       & = \int_{\hat{\C}} \gamma_{0}(\xi ) |\d \xi |^{2} \in E_{0}
\end{align}
and verify readily that it is independent of the section representing a cohomology
class. 
The fact that $E_{0}$ and $\check{\hat{E}}_{0}$ are canonically isomorphic
now follows from the fact that they are both canonically isomorphic to
(different gradings of) the total cohomology of the double complex $\Dou$. 

\paragraph{Step 2.} The first thing to do is to describe explicitly the
isomorphism obtained above. Let $\left[\check{\hat{\delta}}_{0} \right]$ be an element in
$\check{\hat{E}}_{0}$: it is a class in the cohomology space 
$\S_{2}^{1,1}$ in the spectral sequence corresponding to the first
filtration of $\Dou$. Hence it is represented by a $(1,1)$-form
$\check{\hat{\delta}}_{0}(z;\xi )$ over $\C \times \hat{\C}$ such that 
\begin{enumerate}
\item{ $(D-\xi \d z \land) \check{\hat{\delta}}_{0}(z;\xi ) =0$}
\item{ $(\hat{\d}-z \d \xi \land)^{\sharp }\check{\hat{\delta}}_{0}(z;\xi )=0$; in other words, there exists a
$(0,2)$-form $\gamma_{0}(z;\xi ) $ over $\C \times \hat{\C}$
satisfying 
$$
      D_{\xi } \gamma_{0}(z;\xi )  = (\hat{\d}-z \d \xi \land)\check{\hat{\delta}}_{0}.
$$}
\end{enumerate}
Concatenating the map 
$$ 
     \left[ \check{\hat{\delta}}_{0} \right] \mapsto \gamma_{0}(z;\xi ) 
$$ 
with an analog of (\ref{explisom}), namely 
\begin{equation} \label{eval0}
     [\gamma_{0}( z;\xi ) ] \mapsto eval_{z=0}
     \int_{\hat{\C}} \gamma_{0}( z;\xi )  
\end{equation}
we get the canonical isomorphism 
$$
     \omega_{0}: \left[ \check{\hat{\delta}}_{0} \right] \mapsto \delta_{0} =  eval_{z=0}
     \int_{\hat{\C}} \gamma_{0}( z;\xi )  
$$
between $\check{\hat{E}}_{0}$ and $ E_{0}$ provided by the previous step.

Fix now an arbitrary $z_{0} \in \C$, and consider the double complex
$\Dou_{z_{0}}$ having the same $(p,q)$-components as $\Dou$, but with
differentials $\d_{1}=D_{\xi }, \d_{2}=\hat{\d}-(z-z_{0}) \d \xi \land $. In order
to obtain the components of the first page $(\S_{z_{0}})_{1}^{\bullet , \bullet}$ of the
spectral sequence corresponding to the first filtration of $\Dou_{z_{0}}$,
we need to take cohomology with respect to $\d_{1}$, hence these 
will be the same as those of $\Dou$ in (\ref{D1}), and the
differentials  will be induced by $\d_{2}$. Now since $z_{0}$ is a
constant, observe that for any local section 
$\phi_{\xi }(z) \in Ker \Dir_{\xi}^{\ast}$ in $\xi $ of harmonic sections over $\C$ the relation 
$$
     \d_{2}^{\sharp } \phi_{\xi } = [(\hat{\d}-(z-z_{0}) \d \xi \land) \phi_{\xi }] = 
     [(\hat{\d}-z \d \xi \land) \phi_{\xi }] + z_{0} \d \xi \land \phi_{\xi }=
     \hat{D}_{z_{0}}(\phi_{\xi }), 
$$
holds, where $\hat{D}_{z_{0}}$ is the deformation of $\hat{D}$ introduced
in (\ref{invdeform}). To get the second page of the spectral sequence, we
take cohomology with respect to $\d_{2}^{\sharp } = \hat{D}_{z_{0}}$, and
therefore if $z_{0}$ does not belong to the set of opposites of eigenvalues of the
leading term of $\hat{D}$ then this is a finite-dimensional space, equal by
definition to $\check{\hat{E}}_{z_{0}}$. Notice that by the results of
Subsection \ref{sing}, the set of $z_{0}$ where this does not hold is
exactly $P$, the set of singularities (at finite points) of $E$. 
Similarly, the second filtration of $\Dou_{z_{0}}$ gives rise to a spectral
sequence whose first page is (analogously to (\ref{D2})) 
\begin{eqnarray*} \label{Dz2}
     \xymatrix{ 0 & 0 & 0 \\
                0 & 0 & 0 \\ 
                 L^{2}(\C,E) e^{(z-z_{0})\xi } \ar[r]^(.45){\d_{1}^{\sharp }} &
                 L^{2}(\C,\Omega_{\C}^{1} \otimes  E) e^{(z-z_{0})\xi } \ar[r]^{\d_{1}^{\sharp }} & 
                 L^{2}(\C,\Omega_{\C}^{2} \otimes  E) e^{(z-z_{0})\xi }. }
\end{eqnarray*}
Hence the only fiber $\{ z \} \times \hat{\C}$ over which these spaces are
non-trivial is for $z=z_{0}$, so we may consider the double complex 
$\germz$ whose components are germs of forms in a neighborhood
$V_{z_{0}}\times\hat{\C}$ of the fiber $\{ z_{0} \} \times \hat{\C}$, two such germs
being identified if they coincide in any such neighborhood, and with
differentials coming from those of $\Dou_{z_{0}}$. As before, the spectral
sequences corresponding to the second filtration of these double complexes 
agree starting from the first page, so in particular their total
cohomologies are the same. Now, we pass back again to the first filtration
and compute the spectral sequence of $\germz$ with respect to it: in
a convenient trivialisation of $E$ in $V_{0}$ and gauge, the first page is
equal to 
\begin{eqnarray} 
      \label{germDz1} \xymatrix{ \Omega_{\hat{\C}}^{2} \otimes L^{2}(\hat{\C},E_{z_{0}}) e^{z\xi } & 0  & 0\\
                \Omega_{\hat{\C}}^{1} \otimes L^{2}(\hat{\C},E_{z_{0}}) e^{z\xi } \ar[u]^{\d_{2}^{\sharp }} & 0 & 0 \\ 
                  L^{2}(\hat{\C},E_{z_{0}}) e^{z\xi } \ar[u]^{\d_{2}^{\sharp }} & 0 & 0, }
\end{eqnarray}
with differentials given by $\d_{2}=\hat{\d}-(z-z_{0}) \d \xi \land $. 
As in step 1, the second page therefore contains only one non-vanishing
component: the one corresponding to bidegree $(0,2)$, and it is canonically
isomorphic to the vector space $E_{z_{0}}$; this proves that the vector
spaces $E_{z_{0}}$ and $\check{\hat{E}}_{z_{0}}$ are canonically
isomorphic to each other. Again, an element 
$\left[ \check{\hat{\delta}}_{z_{0}} \right]$ of
$\check{\hat{E}}_{z_{0}}$ is represented by a $(1,1)$-form
$\check{\hat{\delta}}_{z_{0}}(z;\xi )$ over $\C \times \hat{\C} $
satisfying 
$(\hat{\d}-(z-z_{0}) \d \xi )^{\sharp }\check{\hat{\delta}}_{z_{0}}(z;\xi )=0$, i.e. there exists a 
$(0,2)$-form $\gamma_{z_{0}}(z;\xi )  $ over $\C \times \hat{\C}$
with
$$
      D_{\xi } (\gamma_{z_{0}}(z;\xi ) ) = (\hat{\d}- (z-z_{0}) 
      \d \xi \land )\check{\hat{\delta}}_{z_{0}}(z;\xi ), 
$$
and an explicit way of describing the obtained isomorphism is given by
\begin{equation} \label{evalz}
     \omega_{z_{0}} : \left[\check{\hat{\delta}}_{z_{0}} \right] \mapsto \delta_{z_{0}} = eval_{z=z_{0}}
     \int_{\hat{\C}} \gamma_{z_{0}}( z;\xi )   
\end{equation}

\paragraph{Step 3.} By the previous points, we have that the bundle 
$\check{\hat{E}}$ is isomorphic to $E$ via the isomorphisms $\omega_{\bullet }$.
Now we prove that the integrable connection 
$\check{\hat{D}}$ on $\check{\hat{E}}$ is carried into $D$ on $E$ by this 
bundle isomorphism: for this, it is clearly sufficient to prove that any
local parallel section for $\check{\hat{D}}$ is carried into a parallel
section for $D$. 
For simplicity, we shall consider a local section near $w=0$, but we 
will see that the proof does not use this. 

For this purpose, we need to work on the product $\C \times \hat{\C} \times \C$,
parametrized by $(z, \xi , w)$; we keep on writing the variable $w$ in lower index. 
We shall consider  $E$ as being a bundle over this space by pull-back,
without writing it out explicitly. Let $\left[\check{\hat{\delta}}_{w} \right]$ be a
$\check{\hat{D}}$-parallel local section of $\check{\hat{E}}$. 
As in Step 2, such a section is represented by giving a global
$(1,1)$-form $\check{\hat{\delta}}_{w}(z;\xi )$ of $E$ on $\C \times \hat{\C}$ for each $w$ in a neighborhood
$V_{0}$ of $0 \in \C$, verifying 
\begin{enumerate}
\item{$D_{\xi_{0}}\check{\hat{\delta}}_{w}(z;\xi )=0$ for all fixed $w_{0} \in V_{0}$
    and $\xi_{0} \in \hat{\C}$} 
\item{$(\d_{2} -(z- w_{0}) \d \xi \land)^{\sharp} \check{\hat{\delta}}_{w}(z;\xi )=0$  for
    all fixed $w_{0} \in V_{0}$}
\item{the section in $w$ of the cohomology classes of the above elements 
    is $\check{\hat{D}}$-parallel.} 
\end{enumerate}
By Hodge theory, we may suppose that $\check{\hat{\delta}}_{w_{0}}(z;\xi_{0} )$ is
the $D_{\xi_{0}}$-harmonic representative of 
$\left[ \check{\hat{\delta}}_{w_{0}}|_{\C \times \{ \xi_{0} \}} \right]$ 
and also that $\check{\hat{\delta}}_{w_{0}}(z;\xi )$ is the
$\hat{D}_{w_{0}}$-harmonic representative of $\left[\check{\hat{\delta}}_{w_{0}} \right]$. 
This way we rephrase the above conditions as 
\begin{enumerate}
\item{for all fixed $w_{0} \in V_{0}$ and $\xi_{0} \in \hat{\C}$
      its restriction to the fiber $\C \times \{ \xi_{0} \} \times \{ w_{0} \}$ is 
      in $\hat{E}_{\xi_{0}}$, that is 
      $\Dir_{\xi_{0}}^{\ast } \check{\hat{\delta}}_{w_{0}}(z;\xi_{0})=0$} \label{cond1}
\item{for all fixed $w_{0} \in V_{0}$ the global section in $\xi $ of the above
    elements of $\hat{E}_{\xi}$ is in $\check{\hat{E}}_{w_{0}}$, in different
    terms 
    $\hat{\Dir}_{w_{0}}^{\ast }\check{\hat{\delta}}_{w_{0}}(z;\xi )=0 $}\label{cond2}
\item{and for all $w_{0} \in V_{0}$, 
  $\check{\pi}_{w}\circ (\check{\d}+\xi \d w \land)\check{\hat{\delta}}_{w}(z; \xi)|_{w=w_{0}}=0$.} \label{cond3}
\end{enumerate} 
As before, (\ref{cond2}) means that for all $w \in V_{0}$ there exists
$\gamma_{w}(z;\xi ) \in \Gamma(\C \times \hat{\C}, \Omega^{2,0} \otimes E)$ such that 
\begin{equation} \label{gammadelta1}
    D_{\xi } \gamma_{w}(z;\xi ) = (\hat{\d}-(z-w) \d \xi \land) \check{\hat{\delta}}_{w}(z;\xi );
\end{equation}
and by Hodge theory, such a section can be defined by the formula 
\begin{equation} \label{gammadelta2}
      \gamma_{w}(z;\xi ) = G_{\xi } D_{\xi }^{\ast } (\hat{\d}-(z-w) \d \xi \land) \check{\hat{\delta}}_{w}(z;\xi ),
\end{equation}
where $G_{\xi }$ is the Green's operator of $\Dir_{\xi }^{\ast } \Dir_{\xi }$. (Here
we used that $G_{\xi }$ is diagonal with respect to the
decomposition $\Omega_{\C}^{0} \oplus \Omega_{\C}^{2}$, a standard consequence of the fact
that $\Dir_{\xi }^{\ast } \Dir_{\xi }$ is diagonal with respect to the same
decomposition, which comes immediately from harmonicity of the metric.)
Now by (\ref{evalz}) and (\ref{gammadelta1}) we have 
\begin{align*}
    D \delta (w)|_{w=w_{0}} = &  D \left( eval_{z=w} \int_{\hat{\C}} \gamma_{w}( z;\xi )   \right)_{|{w=w_{0}}}\\
     = & \int_{\hat{\C}} D \gamma_{w_{0}}( z;\xi )|_{z=w_{0}} + \check{\d}  \gamma_{w}( w_{0};\xi )|_{w=w_{0}} \\
     = & \int_{\hat{\C}} \xi \d z \land \gamma_{w_{0}}( w_{0};\xi ) \\
      & +(\hat{\d}-(w_{0}-w_{0}) \d \xi \land) \check{\hat{\delta}}_{w_{0}}(w_{0};\xi ) 
      + \check{\d}  \gamma_{w}( {w_{0}};\xi )|_{w=w_{0}} 
\end{align*}
(remember that $\check{\d}$ stands for
the trivial connection with respect to $w$ in the trivial Hilbert bundle
$\check{H}$, whereas $\hat{\d}$ is the trivial connection with respect to
$\xi $ in the trivial Hilbert bundle $\hat{H}$). 
The integral of the middle term in this last formula vanishes by Stokes's
theorem. Furthermore, on the diagonal $z=w$ of $\C \times \C$ we have $\d z =\d w$, 
so we are left with 
$$
      \int_{\hat{\C}} (\check{\d} + \xi \d w \land ) \gamma_{w_{0}}( {w_{0}};\xi ).
$$
Applying to this quantity (\ref{gammadelta2}) and the commutation relations 
\begin{align} \label{comm1}
    [ \check{\d} + \xi \d w \land , \hat{\d}-(z-w) \d \xi \land ] = 0 & &
    [ \check{\d} + \xi \d w \land , D_{\xi }] = 0
\end{align}
we obtain 
\begin{equation} \label{int1}
    \int_{\hat{\C}} G_{\xi } D_{\xi }^{\ast } (\hat{\d}-(z-w) \d \xi \land) (\check{\d}+ \xi \d w \land)  
    \check{\hat{\delta}}_{w_{0}}(w_{0};\xi ).
\end{equation}
Consider now condition (\ref{cond3}) above: denoting by $\hat{\Dir}_{w}$
and $\hat{\Dir}_{w}^{\ast }$ the positive and negative Dirac operators of 
the deformation $\hat{D} + w \d \xi $, moreover by $\hat{G}_{w}$ the Green's
operator of $\hat{\Dir}_{w}^{\ast } \hat{\Dir}_{w}$, it can be rewritten as 
\begin{equation*}
     (\Id - \hat{\Dir}_{w}  \hat{G}_{w}  \hat{\Dir}_{w}^{\ast }) (\check{\d}+
     \xi \d w \land)  \check{\hat{\delta}}_{w}(z;\xi ) = 0.    
\end{equation*}
In order to finish the proof, it is sufficient to prove the commutation
relation 
\begin{equation} \label{comm2}   
     [\check{\d}+ \xi \d w \land , \hat{\Dir}_{w}] = 0.
\end{equation}
Indeed, this then implies 
\begin{align*}
     [\check{\d}+ \xi \d w \land , \hat{\Dir}_{w}^{\ast }] = 0 & & [\check{\d}+ \xi
     \d w \land , \hat{G}_{w}] = 0,
\end{align*}
and interchanging $\check{\d}+ \xi \d w \land $ turn by turn with
$\hat{\Dir}_{w_{0}}^{\ast }$, $\hat{G}_{w_{0}}$ and $\hat{\Dir}_{w_{0}}$ 
using each time condition (\ref{cond2}), we get  
\begin{align*} 
      (\check{\d}+ \xi \d w \land ) \check{\hat{\delta}}_{w_{0}}(w_{0};\xi ) & = 
      (\check{\d}+ \xi \d w \land ) (\Id - \hat{\Dir}_{w_{0}}  \hat{G}_{w_{0}}
      \hat{\Dir}_{w_{0}}^{\ast }) \check{\hat{\delta}}_{w_{0}}(w_{0};\xi ) \\
       & = (\Id - \hat{\Dir}_{w_{0}}  \hat{G}_{w_{0}}  \hat{\Dir}_{w_{0}}^{\ast }) 
       (\check{\d}+ \xi \d w \land)  \check{\hat{\delta}}_{w_{0}}(w_{0};\xi ) \\
       & = 0,
\end{align*}
and so (\ref{int1}) is equal to $0$; but on the other hand it is just the
expression for $D \delta (w)|_{w=w_{0}}$, and this shows that $\delta (w)$ is parallel in
$w_{0}$. There remains to show (\ref{comm2}): recall that $\hat{\Dir}_{w} =
\hat{D}_{w}-\hat{D}_{w}^{\ast }$, with 
$$
     \hat{D}_{w} = \hat{\pi}_{\xi } (\hat{\d} - (z-w) \d \xi ). 
$$
Now the first relation in (\ref{comm1}) and 
$\hat{\pi}_{\xi } = (\Id - \Dir_{\xi} G_{\xi } \Dir_{\xi}^{\ast } ) $ combined with the
second relation in (\ref{comm1}) show that 
$$
     [\check{\d}+ \xi \d w \land , \hat{D}_{w}] = 0, 
$$
and we conclude.

\paragraph{Step 4.} Here we wish to show that the double transformed metric
$\check{\hat{h}}$ is equal to $h$. In Step 3 we have already shown that the flat connections $D$ and 
$\check{\hat{D}}$ agree. On the other hand, using the results of Section \ref{sectharm} twice, we see that 
$\check{\hat{h}}$ is a harmonic metric for $\check{\hat{D}}=D$. 
Therefore by uniqueness (up to a constant) of the harmonic metric corresponding to an
integrable connection, we get that $\check{\hat{h}}=h$.  

An equivalent way of deducing the same assertion would be as follows: using again the already proved 
equality $\check{\hat{D}}=D$ and uniqueness of the harmonic metric, we will be done if we can prove  
that the unitary part $\check{\hat{D}}^{+}$ (with respect to $\check{\hat{h}}$) 
of the double transformed flat connection $\check{\hat{D}}$ is equal to 
$D^{+} $, the unitary part of $D$ with respect to $h$. This can be done in a completely analogous way to Steps
1-3. The changes we have to make are the following: consider the double complex $\Dou^{H}_{z_{0}}$ having the
same components as $\Dou_{z_{0}}$, but with differentials $\d_{1}=D_{\xi }^{H}$ and 
$\d_{2}= \hat{\d} - z/2 \d \xi \land -  \bar{z}/2 \d \bar{\xi} \land $. One establishes that these operators commute,
therefore $\Dou^{H}_{z_{0}}$ really forms a double complex. We then see from (\ref{hdecomp}) that the deformation 
$$
     \hat{D}^{H}_{w} = \hat{D}^{H} + \frac{1}{2} w \d \xi \land +  \frac{1}{2} \bar{w} \d \bar{\xi} \land 
$$ 
induced from the differential 
$$
     \hat{\d} - \frac{1}{2} (z-w) \d \xi \land - \frac{1}{2} (\bar{z} -  \bar{w}) \d \bar{\xi} \land 
$$
is the natural deformation of the Higgs-bundle structure induced by the deformation $\hat{D}_{w}$. 
In concrete terms, they are related by the gauge transformation $g^{-1}$. Therefore the double transformed
bundle $\check{\hat{E}}^{H}$ is isomorphic to $g^{-1}gE=E$, and the unitary connection 
$$
   \check{\hat{D}}^{+} = \check{\pi}_{w} \circ \left( \check{\d} + \frac{\xi }{2} \d w \land + 
                            \frac{\bar{\xi} }{2} \d \bar{w} \land \right)
$$
is identified to $D^{+}$ just as $\check{\hat{D}}$ with $D$, using the commutation relations 
\begin{align*}
         \left[\check{\d} + \frac{\xi }{2} \d w \land +  \frac{\bar{\xi} }{2} \d \bar{w} \land  , 
           \hat{\d} - \frac{1}{2} (z-w) \d \xi \land - \frac{1}{2} (\bar{z} -  \bar{z}) \d \bar{\xi} \land \right] =0, 
      \\ \left[\check{\d} + \frac{\xi }{2} \d w \land +  \frac{\bar{\xi} }{2} \d \bar{w} \land  , 
           D_{\xi }^{H} \right] =0 
\end{align*}
instead of (\ref{comm1}), which together imply the analog 
$$
       \left[\check{\d} + \frac{\xi }{2} \d w \land +  \frac{\bar{\xi} }{2} \d \bar{w} \land  , 
       \hat{\Dir}_{w}^{H} \right] = 0
$$
of (\ref{comm2}) for the deformed Dirac operator 
$$
     \hat{\Dir}_{w}^{H} = D_{w}^{H} - (D_{w}^{H})^{\ast }.   
$$
This then allows us to conclude equality of the unitary connections. 

Since the Hermitian bundles $(\check{\hat{E}}, \check{\hat{h}})$ and $(E,h)$ coincide, so do the flags
of their parabolic structures in the singular points; as well as the
parabolic weights, because they are supposed to be between $0$ and $1$, and there is a unique way of choosing
holomorphic sections with such behaviors. 
\end{proof}

\end{document}